\documentclass{tac1}

\usepackage{xy}

\input diagxy

\usepackage[colorlinks=true]{hyperref}
\hypersetup{allcolors=[rgb]{0.1,0.1,0.4}}

\usepackage{bbm}				    
\usepackage{eufrak}					
                                    
\usepackage{graphicx}				
\usepackage[makeroom]{cancel}       

\author{Camilo Angulo}

\thanks{The author would like to acknowledge the financial support of CAPES and the \emph{National Council for Scientific and Technological Development} - CNPq during his Ph.D. studies at the University of S\~ao Paulo, where this research was carried out, as well as the FAPERJ grant number E-26/202.439/2019.}

\address{Departamento de Matem\'atica Aplicada, Instituto de Matem\'atica da Universidade Federal Fluminense\\
 Rua Professor Marcos Waldemar de Freitas Reis, s/n, Campus do Gragoat\'a, Niter\'oi, RJ, Brazil, 24210-201.\\
}

\title{Towards a new cohomology theory for strict Lie $2$-groups}

\copyrightyear{2020}

\keywords{Cohomology, higher Lie theory}
\amsclass{17B56, 18D05, 22A22}

\eaddress{ca.angulo951@gmail.com}

\mathcaldef[G]{G}												   
\mathcaldef[Lie]{L}					 						       
\newcommand{\ggl}{\mathfrak{gl}}                                   
\renewcommand{\gg}{\mathfrak{g}}            				       
\newcommand{\hh}{\mathfrak{h}}            						   
\newcommand{\Rr}{\mathbbm{R}}    		 						   
\renewcommand{\Join}{\bowtie}									   
\newcommand{\vJoin}{\mathbin{\rotatebox[origin=c]{90}{$\bowtie$}}} 
\newcommand{\Tp}[1]{\overline{\underline{#1}}}       			   
\newcommand{\Lf}[1]{\vert{#1}\vert}                			       
\newcommand{\Rg}[1]{{#1}_V}                         			   
\newcommand{\ad}{\textnormal{ad}}           					   

\newcommand{\ltimes}{\mathrel{\times\mkern-12.5mu\raise0.3ex\hbox{\scalebox{0.5}{$\vert$}}}~} 
\newcommand{\rtimes}{\mathbin{\rotatebox[origin=c]{180}{$\ltimes$}}} 

\mathchardef\gt="313E 
\mathchardef\lt="313C 

\newtheorem{thm}{Theorem}
 
\newtheorem{lem}{Lemma}

\newtheoremrm{rem}{Remark}

\let\pf\proof
\let\epf\endproof
\mathrmdef{Hom}
\mathbfdef{Set}

\begin{document}

\maketitle
\begin{abstract}
In this article, we introduce the first degrees of a cochain complex 
associated to a strict Lie 2-group whose cohomology is shown to extend 
the classical cohomology theory of Lie groups. In particular, we show 
that the second cohomology group classifies an appropriate type of 
extensions. We conclude putting forward evidence that this complex can 
be extended to arbitrary degrees.
\end{abstract}


\section{Introduction}\label{sec-Introduction}

Crossed modules enjoyed a lively research activity in the late 1940s after 
Whitehead introduced them to classify $2$-types \cite{Whitehead}. Said 
activity was arguably revived after Loday \cite{Loday} justified the 
classification of all types using the so-called $n$-cat Groups. This 1990s 
revival, spearheaded by \cite{Norrie}, aimed at studying crossed modules as 
algebraic objects in their own right. This paper is dedicated to crossed 
modules having a smooth structure, thought of as the global counterpart to 
Lie $2$-algebras \cite{Lie2Alg}. Recall the definition of a 
\emph{Crossed module} ({\em e.g.}, \cite{Lie2Gps}).

\definition\label{crossMod} 
A {\em crossed module of Lie groups} is a Lie group homomorphism $G\to^i H$ 
together with a right action of $H$ on $G$ by Lie group automorphisms satisfying
$$i(g^h)=h^{-1}i(g)h, $$
$$g_1^{i(g_2)}=g_2^{-1}g_1g_2, $$
for all $g,g_1,g_2\in G$ and $h\in H$, where we write $g^h$ for $h$ acting on $g$. 
Following the convention in the literature, we refer to these equations 
respectively as {\bf equivariance} and {\bf Peiffer}. \enddefinition

Throughout, we use the equivalence between the category of crossed modules 
with that of strict Lie $2$-groups (see, {\em e.g.}, \cite{Ellis}); henceforth, 
we refer to the latter simply as \emph{Lie 2-groups}.

\definition\label{Lie2-gp} 
A {\em Lie $2$-group} is a groupoid object internal to the category of Lie groups. \enddefinition

Though several categorifications of the group structure have made their 
way into the literature \cite{Lie2Gps,2Coh,IntInfDim,StackyInt} proving 
specially useful in applications to mathematical physics (see \cite{Lie2Gps} 
for a list of references) and to Lie theoretical problems (\cite{IntInfDim,StackyInt}), 
we restrict ourselves to this seemingly restrictive subclass. We do so 
because, in either presentation, one can naturally differentiate the 
structure of a Lie $2$-group to get a Lie $2$-algebra. 

Using Lie algebra paths, it was proven in \cite{ZhuInt2Alg} that all 
finite-dimensional Lie $2$-algebras arise this way, by differentiating a 
Lie $2$-group structure. This paper sprout out of attempting to find a 
cohomological proof of the integrability of strict Lie $2$-algebras to 
strict Lie $2$-groups following along the lines of the strategy devised by 
van Est \cite{VanEst,VanEstC}. This approach still works in infinite 
dimensions and was historically used to construct the first example of a 
non-integrable Lie algebra \cite{NonInt}; thus, it bears the potential to 
improve our current understanding of the Lie theory of other categorified 
objects (see, {\em e.g.}, \cite{BCD,LucaPhD,Neeb}).

Van Est's strategy could be roughly summarized as follows: 
Given a Lie algebra $\gg$, one uses its adjoint representation to 
recast it as an extension of the Lie subalgebra $\ad(\gg)\leq\ggl(\gg)$ 
by the center $\mathfrak{z}(\gg)$. There exists a Lie algebra cohomology 
class $[\omega]\in H^2(\ad(\gg),\mathfrak{z}(\gg))$ that classifies this 
extension. Since linear Lie algebras are integrable and one can always 
pick a $2$-connected integration, the van Est Theorem says there is a 
group cohomology class whose associated extension is a Lie group 
integrating $\gg$. Implicitly, the preliminary step to carry out van 
Est's strategy is to have cohomology theories that respectively classify 
abelian extensions of Lie groups and Lie algebras, or rather of the 
global and infinitesimal counterparts, as well as a van Est map relating 
them (see, {\em e.g.}, \cite{AAC} for the extensions and \cite{VanEstC,CamFlor} 
for the van Est map in the Lie groupoid/algebroid case). 

This paper is the sequel to \cite{Lie2AlgCoh}, where the author introduced 
a cohomology theory for Lie 2-algebras that suitably classifies abelian 
extensions. In the present article, we give part of a complex that serves 
as the global counterpart to the one in \cite{Lie2AlgCoh}. There currently 
are several cohomologies associated with a crossed module that classify 
specific types of extensions in the literature \cite{Ellis,Baues1,Vietes,(Co)Ho}; 
however, it is worth pointing out that most of these disregard the topology. 
For instance, \cite{(Co)Ho} is based on a research program that rests 
upon the tripeability of the underlying set functor from the category of 
crossed modules to that of sets. The left adjoint to the referred functor 
is a composition of the free group functor and taking free products of 
a group with itself, both of which have little meaning in the smooth 
category, for either one loses all topological data, or ends up 
with an infinite-dimensional object instead. Such is the general 
drawback impeding applying the machinery developed after \cite{Norrie} to  
the smooth category, where, as another instance, the second homology is 
computed in terms of {\em generators and relations}.  

The well-established cohomology of categories does not suffice either. 
Indeed, taking values in a natural system, its second cohomology relates 
to the so-called linear extensions. The latter, however, are defined to 
have the same space of objects as the category under consideration. For 
the application we have in mind, more general coefficients are needed. 
The complex is prescribed to take coefficients on {\bf $2$-vector spaces}, 
{\it i.e.}, on flat abelian Lie $2$-groups or, equivalently on $2$-term 
complexes of vector spaces.

\definition\label{defRep} 
A \emph{representation} of a Lie $2$-group $\G$ on $W\to^\phi V$ is a 
morphism of Lie $2$-groups
\begin{equation}\label{2rep}
\rho :\G \to<350> GL(\phi ).
\end{equation} \enddefinition

We remark that, by a morphism of Lie $2$-groups, we mean a na\"ive smooth 
functor respecting the Lie group structures, as opposed to more general 
types of morphisms (such as bibundles of Lie groupoids). We refer to 
a map (\ref{2rep}) as a $2$-representation. The co-domain of $\rho$ in 
Definition \ref{defRep} is a category of linear automorphisms and linear 
natural transformations which happens to be a Lie $2$-group (see 
Proposition \ref{GLPhi}). A $2$-representation taking values in either 
$A\to 1$ or $1\to A$ coincides with the modules of \cite{Ellis,Baues1} 
if the abelian group $A$ is simply connected. 

In \cite{Ellis}, it is explained that the groupoid nature of a Lie 
$2$-group can be exploited to build a bi-bisimplicial object whose 
associated double complex ({\em cf.} Proposition \ref{DoubleGpdDoubleCx}) 
with values in an abelian Lie group $A$ classifies extensions by the unit 
Lie groupoid $A\two A$. The novel complex we are to outline is an extension 
of the referred double complex, though allowing the space of values to be 
an arbitrary $2$-vector space.   

We define the first degrees of the complex of Lie $2$-group cochains of 
$\G$ with values in the $2$-representation (\ref{2rep}) as the total 
complex of a triple complex of sorts. Assuming the convention that 
$\G_0=H$, $\G_1=\G$ and $\G_p$ is the space of $p$-composable arrows 
usually noted $\G^{(p)}=\G\times_{H}...\times_{H}\G$, we define
\begin{eqnarray}\label{3dimLat}
C^{p,q}_r(\G,\phi):=C(\G^q_p\times G^r,W)
\end{eqnarray}
for $r\neq 0$, and 
\begin{eqnarray}\label{3dimLat0}
C^{p,q}_0(\G,\phi):=C(\G^q_p ,V),
\end{eqnarray}
where $C(\G_p^q\times G^r,A)$ is the vector space of $A$-valued smooth 
functions.

This three dimensional lattice of vector spaces can be enhanced to a  
{\em grid} of complexes of Lie groupoid cochains in each direction by 
placing the complex  associated to certain Lie group bundles in the 
$r$-direction; the complex associated to certain action groupoids in the 
$q$ direction; and the complex associated to the product between powers 
of the Lie $2$-group and powers of the unit groupoid $G\two G$ in the 
$p$-direction. All these complexes take values in different representations 
on either the trivial vector bundle with fibre $V$ or the trivial vector 
bundle with fibre $W$ (see Section \ref{sec-tetrahedralCx} for a detailed 
explanations of the groupoids and the representations involved). 

We refrain from calling this grid a triple complex because not all these 
differentials pairwise commute; thus failing to build a complex by taking 
their alternated sum. The page $r=0$ commutes only up to isomorphism in 
the $2$-vector space while the successive $r$-{\em pages} commute only up 
to homotopy (see Proposition \ref{FrontPage-upToIso}). Adding these 
homotopies that we call \textit{difference maps} to the total differential 
takes care of the pair of non-commuting differentials; however, in order, 
it makes other coordinates in the square of the total differential not 
necessarily vanish. We show that, in the lowest degrees, the difference 
maps are subject to higher relations encoded by homotopies that we call 
accordingly \textit{higher difference maps}. By adding these higher 
difference maps to the total differential we thus ensure that it squares 
to zero.

Summing up, grading by counter-diagonal planes and letting $\delta_{(1)}$, 
$\delta$ and $\partial$ be the differentials of the complexes in the $r$, 
$q$ and $p$ directions respectively and denoting the difference maps by 
$\Delta$ and the second difference maps by $\Delta_{a,b}$, we get the 
following result.

\thm\label{The2GpCx}
Associated to a Lie $2$-group $\G$ together with a $2$-representation $\rho$, 
there is a (truncated) complex $(C^{\leq 3}_{tot}(\G,\phi),\nabla)$ with 
$$C^{n}_{tot}(\G,\phi)=\bigoplus_{p+q+r=n}C^{p,q}_r(\G,\phi) $$
and
\begin{equation}\label{preDiff3}
\nabla=(-1)^p\Big{(}\delta_{(1)}+\partial+\Delta+\Delta_{1,2}+(-1)^r(\delta+\Delta_{2,1})\Big{)}
\end{equation} \endthm

For the purpose of this note, the main application of Theorem \ref{The2GpCx} 
is that, as desired, the second cohomology of the defined complex 
classifies abelian extensions (see Theorem \ref{H2Gp}).

The main advantage of the complex of Theorem \ref{The2GpCx} is that it is 
built out of complexes of Lie groupoids; hence, one can directly generalize 
the van Est map by twice-assembling usual van Est maps to land in the 
complex of \cite{Lie2AlgCoh}. This will be explored in the separate paper 
\cite{2VanEst}.
 
This paper is organized as follows. In Section \ref{sec-Pre}, we recall 
some basic facts and convene notation. We motivate the emergence 
of the complex and its differential by recalling the canonical double 
complex associated to a Lie $2$-group whose second cohomology 
classifies extensions of a particular type and the complex of \cite{Ellis}. 
Then, we recall the definition of the general linear Lie $2$-group, 
spell out the definition of $2$-representations, provide examples and 
some associated constructions. We conclude the section by proving that 
$2$-representations are indeed the kind of actions induced by extensions 
(see Proposition \ref{Ind2GpRep}). In Section \ref{sec-tetrahedralCx}, we 
carefully define the three dimensional grid and the difference maps out 
of which we get the truncated $(4\times 4\times 4)$-tetrahedral complex 
of Lie $2$-group cochains with values in a $2$-representation, and study 
its cohomology. In particular, we show that the equations that define a 
$2$-cocycle are equivalent to the equations defining an abstract 
extension. In Section \ref{sec-Ink}, we prove the general relations that 
the differentials in the background grid verify and heal the non-commuting 
part by introducing the general formula for the difference maps. We 
conclude by discussing what is needed to fully extend the grid to a complex.
We add an appendix with the general formula for the second difference 
maps, as well as the necessary maps to extend the complex of Theorem 
\ref{The2GpCx} to degree $5$ by taking a $(6\times 6\times 6)$-tetrahedral 
slice of the grid of Section \ref{sec-Ink}.

\section{Preliminaries}\label{sec-Pre}

In this section, we establish the notation conventions used throughout. 
As a motivation, we recall the complex of \cite{Ellis} and study its 
cohomology. We also recall the notions of general linear Lie $2$-group 
and $2$-representation.

\rem\label{Equiv2Gp:Mods} From here on out, we make no distinction 
between a Lie $2$-group and its associated crossed module. For future 
reference, we outline the equivalence between the category of Lie 
$2$-groups and crossed modules of Lie groups at the level of objects (see, 
{\em e.g.}, \cite{Lie2Gps} and \cite{Loday} for details). 

We write a generic Lie $2$-group as
$$\G\times_H\G\to^m\iloop\G(dr,dl)_\iota\two^s_t H\to^u\G . $$
In order to make clear the difference between the group operation and the 
groupoid operation in $\G$, we assume the following convention: 
$$g_1\vJoin g_2 \qquad\qquad g_3\Join g_4, $$
stand respectively for the group multiplication and the groupoid 
multiplication whenever $(g_1,g_2)\in\G^2$ and $(g_3,g_4)\in\G\times_H\G$. 
This notation intends to reflect that we think of the group multiplication 
as being ``vertical'', whereas the groupoid multiplication as  being 
``horizontal''.

Given a crossed module $G\to^i H$ as in Definition \ref{crossMod}, the 
space of arrows of its associated Lie $2$-group $\G$ is defined to be the 
semi-direct product $G\rtimes H$ with respect to the $H$-action, whose 
product is explicitly given by 
\begin{equation}\label{ArrowsProduct}
\pmatrix{g_1\cr h_1}\vJoin\pmatrix{g_2\cr h_2}=\pmatrix{g_1^{h_2}g_2\cr h_1 h_2},
\end{equation}
for $(g_1,h_1),(g_2,h_2)\in G\rtimes H$. The structural maps are given 
by
$$s\pmatrix{g\cr h}=h \qquad\qquad t\pmatrix{g\cr h}=hi(g) \qquad\qquad \iota\pmatrix{g\cr h}=\pmatrix{g^{-1}\cr hi(g)} \qquad\qquad u(h)=\pmatrix{1\cr h} $$
\begin{equation}\label{strMaps}
\pmatrix{g'\cr hi(g)}\Join\pmatrix{g\cr h}:=\pmatrix{gg'\cr h}
\end{equation} 
for $h\in H$ and $g,g'\in G$.

Conversely, given a Lie $2$-group $\G\two H$, let $G$ to be the Lie 
subgroup $\ker s\leq\G$. The associated crossed module is given by 
$G\to^{t\vert_{\ker s}}H$ together with the right action given by 
conjugation by units in the group $\G$
$$g^h :=u(h)^{-1}\vJoin g\vJoin u(h), $$
for $h\in H$ and $g\in G$. We stress that the $-1$ power stands for the 
inverse of the group multiplication $\vJoin$. 

Notice that the isomorphism of vector spaces $\G\cong G\rtimes H$ is 
canonical because the unit map provides a natural splitting. \endrem

\subsection{Lie groupoid cohomology}\label{GpdCoh}
Let $G\two M$ be a Lie groupoid. There is a simplicial structure on the 
nerve of $G$ whose maps are given by 
\begin{equation}\label{faceMaps}
\partial _k (g_0,...,g_p)=\cases{(g_1,...,g_p)&if $k=0$\cr(g_0,...,g_{k-1}g_k,...,g_p)&if $0\lt k\leq p$\cr(g_0,...,g_{p-1})&if $k=p+1$,}
\end{equation} 
for a given element $(g_0,...,g_p)\in G^{(p+1)}$. With these, one builds 
the complex of Lie groupoid cochains $C^p(G):=C^\infty(G^{(p)})$ whose 
differential
$\bfig
\morphism<800,0>[\partial :C^\bullet(G)`C^{\bullet +1}(G);]
\efig$
is defined by the formula
\begin{equation}\label{trivDiff}
\partial\varphi=\sum_{k=0}^{p+1}(-1)^k\partial_k^*\varphi,
\end{equation} 
for $\varphi\in C^p(G)$. 

Thus defined, $(C^\bullet(G),\partial)$ is referred to as the 
{\em groupoid complex} of $G$, and its cohomology is called 
{\em differentiable cohomology} of $G$. 

\rem\label{2GpNerve} Under the isomorphism of Remark \ref{Equiv2Gp:Mods}, 
the space of $p$-composable arrows 
$$\G_p=\G\times_H ...\times_H\G =\lbrace(\gamma_1,...,\gamma_p)\in\G^p :s(\gamma_j)=t(\gamma_{j+1}),\quad 1\leq j\lt p\rbrace $$
corresponds to $G^p\times H$; again, hereafter, we consider this 
isomorphism to be fixed and treat it as an equality, when necessary. For 
each coordinate $\gamma_j$ of $\gamma\in\G_p$, there is a corresponding 
$(g_j,h_j)\in G\rtimes H$. The defining relation for $\G_p$ then reads 
$h_j=h_{j+1}i(g_{j+1})$ (see Eq. (\ref{strMaps})), thus making the map 
$\gamma\to/|->/(g_1,...,g_p;h_p)$ an isomorphism with inverse 
$(g_1,...,g_p;h) \to/|->/(g_1,hi(g_p...g_2);...;g_{p-1},hi(g_p);g_p,h)$. 
Under this isomorphism, we rewrite the face maps to be
\begin{equation}\label{2GpFaceMaps}
\partial_k(g_0,...,g_p;h)=\cases{(g_1,...,g_p;h)&if $k=0$\cr(g_0,...,g_{k-2},g_kg_{k-1},g_{k+1},...,g_p;h)&if $0\lt k\leq p$\cr(g_0,...,g_{p-1};hi(g_p))&if $k=p+1$.}
\end{equation} \endrem

\subsubsection{Representations and cohomology with values} A 
{\em left representation} of $G\two M$ is a vector bundle $E$ over $M$, 
together with a left action 
$$G_s\times_M E\to<350>E:(g,e)\to/|->/<350>\Delta_g e$$
along the projection of the vector bundle. Let 
$t_p:G^{(p)}\to M:(g_1,...g_p)\to/|->/t(g_1)$ be the map that returns the 
final target of a $p$-tuple of composable arrows.

\rem\label{2GpFinalTarget} Under the isomorphism of Remark \ref{2GpNerve}, 
the final target map gets rewritten as
$$\bfig
\morphism<600,0>[t_p:G^p\times H`H:;]
\morphism(950,0)/|->/<1150,0>[(g_1,...,g_p;h)`hi\Big{(}\prod_{j=0}^{p-1}g_{p-j}\Big{)}=hi(g_p...g_1).;]
\efig$$ 
Observe that the final target map is a composition of face maps and hence 
a group homomorphism. \endrem
 
The complex of Lie groupoid cochains with values on the left 
representation $E$ is defined by
$$C^p(G;E):=\Gamma(t_p^*E)$$ 
together with the differential 
$\partial :C^\bullet (G;E)\to C^{\bullet +1}(G;E) $ whose formula is 
essentially (\ref{trivDiff}), though modifying the first term, so that 
all terms lie on the same fibre and the sum can be performed. More 
specifically, for $\varphi\in C^p(G;E)$ and $(g_0,...,g_p)\in G^{(p+1)}$,
\begin{equation}\label{LGpdDifferential}
(\partial\varphi)(g_0,...,g_p) := \Delta_{g_0}\partial_0^*\varphi(g_0,...,g_p)+\sum_{k=1}^{p+1}(-1)^k\partial_k^*\varphi(g_0,...,g_p). 
\end{equation} 
Similarly, a {\em right representation} of $G\two M$ is a vector bundle 
$E$ over $M$, together with a right action 
$$E_M\times_t G\to<350>E:(e,g)\to/|->/<350>\Delta_g e$$
along the projection of the vector bundle for which we use the same 
notation. Replacing each instance of the target map by the source map in 
the preceding discussion, one defines the complex of Lie groupoid 
cochains with values on the right representation $E$, whose differential 
is 
\begin{equation}\label{RGpdDifferential}
(\partial\varphi)(g_0,...,g_p) :=\sum_{k=0}^{p}(-1)^k\partial_k^*\varphi(g_0,...,g_p)+(-1)^{p+1}\Delta_{g_p}\partial_{p+1}^*\varphi(g_0,...,g_p). 
\end{equation}

Both left and right representations can be pulled-back along 
homomorphisms. If $E$ is a left (resp. right) representation of $G\two M$ 
and         
$$\bfig
\square/>`@{>}@<3pt>`@{>}@<3pt>`>/[H`G`N`M;\varphi```f]
\square/`@{>}@<-3pt>`@{>}@<-3pt>`/[``N`M;```]
\efig$$
is a Lie groupoid homomorphism, then there is a left (resp. right) action 
of $H\two N$ on the pull-back bundle $f^*E=N\times_M E$, where $h\in H$ 
acts on $e\in E_{f(s(h))}$ (resp. $E_{f(t(h))}$) by $\Delta_{\varphi(h)}e$.

\subsection{The canonical simplicial object} - In this subsection, we 
recall that, as it is explained in \cite{Ellis}, given a Lie $2$-group 
$\G\two H$, its nerve $\G_\bullet$ is a simplicial group. For each $p$, 
$\G_p$ is a Lie subgroup of $\G^p$ and one can thus consider its nerve.  
Considering simultaneously the nerve of all $\G_p$'s yields a 
bisimplicial set $\G_\bullet^\bullet$. In particular, the face maps of 
the two simplicial structures always commute with one another; hence, 
{\em dualizing}, one gets the double complex:
  
$$\bfig
\node 1a(0,1350)[\vdots]
\node 1b(550,1350)[\vdots]
\node 1c(1100,1350)[\vdots]
\node 1d(1500,1350)[]
\node 2a(0,1000)[C(H^3)]
\node 2b(550,1000)[C(\G^3)]
\node 2c(1100,1000)[C(\G_2^3)]
\node 2d(1500,1000)[\cdots]
\node 3a(0,500)[C(H^2)]
\node 3b(550,500)[C(\G^2)]
\node 3c(1100,500)[C(\G_2^2)]
\node 3d(1500,500)[\cdots]
\node 4a(0,0)[C(H)]
\node 4b(550,0)[C(\G)]
\node 4c(1100,0)[C(\G_2)]
\node 4d(1500,0)[\cdots]
\arrow|l|[4a`3a;\delta]
\arrow[4a`4b;\partial]
\arrow|l|[3a`2a;\delta]
\arrow[3a`3b;\partial]
\arrow|l|[4b`3b;\delta]
\arrow[4b`4c;\partial]
\arrow[2a`1a;]
\arrow[2a`2b;\partial]
\arrow|l|[3b`2b;\delta]
\arrow[3b`3c;\partial]
\arrow|l|[4c`3c;\delta]
\arrow[4c`4d;]
\arrow[2b`1b;]
\arrow[2b`2c;\partial]
\arrow|l|[3c`2c;\delta]
\arrow[3c`3d;]
\arrow[2c`1c;]
\arrow[2c`2d;]
\efig$$

whose columns are complexes of Lie group cochains, and whose $q$th row is 
the groupoid complex of $\G^q\two H^q$. For future reference, we recall  
the total complex of this double complex
$$\Omega ^k_{tot}(\G)=\bigoplus _{p+q=k}C(\G _p^q), $$
with differential $d=(-1)^p(\partial+\delta)$. 	

In Section \ref{pq-pag}, we use a generalization of this construction to 
double Lie groupoids. Recall that a double Lie groupoid is a groupoid 
object internal to the category of Lie groupoids. More explicitly, a 
{\em double Lie groupoid} consists of a square
$$\bfig
\square/@{>}@<3pt>`@{>}@<3pt>`@{>}@<3pt>`@{>}@<3pt>/[D`V`H`M;```]
\square/@{>}@<-3pt>`@{>}@<-3pt>`@{>}@<-3pt>`@{>}@<-3pt>/[D`V`H`M;```]
\efig$$
where each side is a Lie groupoid and the structural maps are smooth 
functors. The elements in $D$ can be thought of as being squares whose 
vertical edges are arrows in $V$, whose horizontal edges are arrows in 
$H$ and all of whose vertices are points in $M$. In order to recognize 
the structural maps then, we adopt the following mnemonic device. We write 
$\Tp{s}$, $\Tp{t}$, etc. for the structural maps of the top groupoid; 
$\Lf{s}$, $\Lf{t}$, etc. for the left vertical groupoid; 
$\Rg{s}$, $\Rg{t}$, etc. for the right vertical groupoid; and concludingly, 
the usual $s$, $t$, etc. for the bottom groupoid. Thus, a given element 
$d\in D$, thought of as a square, has the following edges
$$\bfig
\square/<-`<-`<-`<-/[\bullet`\bullet`\bullet`\bullet;\Lf{t}(d)`\Tp{t}(d)`\Tp{s}(d)`\Lf{s}(d)]
\efig$$
Additionally, as we did for Lie $2$-groups, we use the shorthand 
$$d_1\vJoin d_2:=\Lf{m}(d_1,d_2) \qquad\qquad d_3\Join d_4:=\Tp{m}(d_3,d_4), $$
whenever $(d_1,d_2)\in D\times_H D$ and $(d_3,d_4)\in D\times_V D$ to 
reflect the fact that $d_1\vJoin d_2$ and $d_3\Join d_4$ represent 
respectively
$$\bfig
\iiixiii|aalramalraba|/<-``<-`<-``<-``<-`<-``<-`/{0}[\bullet`\bullet``\bullet`\bullet``\bullet`\bullet`;\Lf{t}(d_1)``\Tp{t}(d_1)`\Tp{s}(d_1)``\frac{\Lf{s}(d_1)}{\Lf{t}(d_2)}``\Tp{t}(d_2)`\Tp{s}(d_2)``\Lf{s}(d_2)`]
\iiixii(1000,250)|aalmrbb|/<-`<-`<-`<-`<-`<-`<-/{0}[\bullet`\bullet`\bullet`\bullet`\bullet`\bullet;\Lf{t}(d_3)`\Lf{t}(d_4)`\Tp{t}(d_3)`\Tp{s}(d_3)=\Tp{t}(d_4)`\Tp{s}(d_4)`\Lf{s}(d_3)`\Lf{s}(d_4)]
\efig$$
With this notation, the fact that the multiplication in either groupoid 
is a groupoid homomorphism yields the formula
$$ (d_1\Join d_2)\vJoin(d_3\Join d_4)=(d_1\vJoin d_3)\Join(d_2\vJoin d_4), $$ 
whenever it makes sense. We call this formula the interchange law.

It is also customary to add the axiom that the double source map
$$\bfig
\morphism<900,0>[\mathbbm{S}:=(\Lf{s},\Tp{s}):D`H_s\times _{\Rg{s}}V;] 
\efig$$
is a submersion, though it is immaterial for our purposes.

If $D$ is a double Lie groupoid, there are two simplicial structures for 
$D$ given by its vertical and its horizontal groupoid structures. The 
commutativity of the structural maps of the square diagram representing a 
double Lie groupoid is but a shadow of the general interaction between 
these two simplicial structures. In what follows, let 
$\vec{d}=(d_1,...,d_q)\in M_{q\times p}(D)$ be represented by the matrix 
$\pmatrix{d_{11} & d_{12} & ... & d_{1p} \cr
d_{21} & d_{22} & ... & d_{2p} \cr
\vdots & \vdots & ... & \vdots \cr
d_{q1} & d_{q2} & ... & d_{qp} \cr
}$, where each $d_{mn}\in D$. In order to distinguish the vertical 
and the horizontal simplicial structures, we do as before and write 
$\partial_k$ for the face maps associated to the horizontal groupoid and 
$\delta_j$ for the face maps associated to the vertical groupoid.

The complexes of groupoid complex of the horizontal and vertical groupoids 
fit into a double complex. Although this is an expected relation, we could 
not find a reference in the literature.
\proposition\label{DoubleGpdDoubleCx}
Given a double Lie groupoid $D$ with our conventions,
$$\bfig
\node 1a(0,1350)[\vdots]
\node 1b(850,1350)[\vdots]
\node 1c(1650,1350)[\vdots]
\node 1d(2150,1350)[]
\node 2a(0,1000)[C(V^{(2)})]
\node 2b(850,1000)[C(D\times_H D)]
\node 2c(1650,1000)[C(D_2^2)]
\node 2d(2150,1000)[\cdots]
\node 3a(0,500)[C(V)]
\node 3b(850,500)[C(D)]
\node 3c(1650,500)[C(D\times_V D)]
\node 3d(2150,500)[\cdots]
\node 4a(0,0)[C(M)]
\node 4b(850,0)[C(H)]
\node 4c(1650,0)[C(H^{(2)})]
\node 4d(2150,0)[\cdots]
\arrow|l|[4a`3a;\delta]
\arrow[4a`4b;\partial]
\arrow|l|[3a`2a;\delta]
\arrow[3a`3b;\partial]
\arrow|l|[4b`3b;\delta]
\arrow[4b`4c;\partial]
\arrow[2a`1a;]
\arrow[2a`2b;\partial]
\arrow|l|[3b`2b;\delta]
\arrow[3b`3c;\partial]
\arrow|l|[4c`3c;\delta]
\arrow[4c`4d;]
\arrow[2b`1b;]
\arrow[2b`2c;\partial]
\arrow|l|[3c`2c;\delta]
\arrow[3c`3d;]
\arrow[2c`1c;]
\arrow[2c`2d;]
\efig$$
is a double complex, where 
$$\mld D_p^q:=\lbrace\vec{d}\in M_{q\times p}(D): &  
                \quad\Tp{s}(d_{mn})=\Tp{t}(d_{mn+1}),\quad\Lf{s}(d_{mn})=\Lf{t}(d_{m+1n}), \\
                \quad\Tp{t}(d_{mn})=\Tp{s}(d_{mn-1}),\quad\Lf{t}(d_{mn})=\Lf{s}(d_{m-1n})\rbrace . $$\endproposition

In the sequel, we refer to this object as the 
{\em double complex associated to} $D$, and the cohomology of its total 
complex $(C_{tot}^\bullet(D),d)$,
$$C_{tot}^k(D)=\bigoplus_{p+q=k}C(D_p^q)\qquad d=(-1)^p(\partial+\delta), $$
as the \textit{double groupoid cohomology of} $D$. 

\subsubsection{Cohomology with trivial coefficients}\label{GpDcx}

We interpret $H^2_{tot}(\G)$ for a Lie $2$-group $\G$ as classifying 
certain type of extensions. Though this is not a remarkable observation, 
it does not follow from neither \cite{Ellis, Baues1} because, in the 
extension, the space of objects is modified.

A $2$-cocycle consists of a pair of functions $(F,f)\in C(H^2)\oplus C(\G)$ 
such that:
\begin{enumerate}
\item $\delta F=0$, i.e. 
$F(h_1 ,h_2)+F(h_0 ,h_1 h_2)=F(h_0 h_1 ,h_2)+F(h_0 ,h_1 )$ for all triples 
$h_0,h_1,h_2\in H$ .
\item $\partial f=0$, i.e.  
$f(\gamma _1 \bowtie\gamma _2 )=f(\gamma _1 )+f(\gamma _2 )$ for all 
$(\gamma _1 ,\gamma _2 )\in\G _2$. Using the isomorphism of Remark 
\ref{Equiv2Gp:Mods}, this is equivalent to 
$f(g_2g_1,h)=f(g_2,h)+f(g_1,hi(g_2))$ for all $h\in H$ and $g_1,g_2\in G$.
\item $\partial F-\delta f=0$, i.e. 
$F(s(\gamma _0),s(\gamma _1))-F(t(\gamma _0),t(\gamma _1))=f(\gamma _1)-f(\gamma _0\vJoin\gamma _1)+f(\gamma _0)$ 
for all pairs $\gamma _0,\gamma _1\in\G$. Again, under the isomorphism of 
Remark \ref{Equiv2Gp:Mods}, this equation can be rewritten as
\begin{equation}\label{mixed}
F(h_0,h_1)-F(h_0i(g_0),h_1i(g_1))=f\pmatrix{g_1\cr h_1}-f\pmatrix{g_0^{h_1}g_1\cr h_0h_1}+f\pmatrix{g_0\cr h_0},
\end{equation}
where $(g_1,h_1),(g_2,h_2)\in G\rtimes H$.
\end{enumerate}
We point out that making $h_0=h_1=1$ in Eq. (\ref{mixed}) yields 
\begin{equation}\label{singledOut}
f\pmatrix{g_0g_1\cr 1}=f\pmatrix{g_1\cr 1}+f\pmatrix{g_0\cr 1}+F(i(g_0),i(g_1)).
\end{equation}
Also, putting $g_2=1$, $\partial f=0$ implies that $f(1,h)=0$ for all 
$h\in H$. 

Since $\delta F=0$, $F$ induces a (central) extension of $H$, 
$$\bfig
 \node a(0,0)[1]
 \node b(400,0)[\Rr]
 \node c(1000,0)[H\ltimes^F\Rr]
 \node d(1500,0)[H]
 \node e(2000,0)[1,]
 \arrow[a`b;]
 \arrow[b`c;\bar{1}\times I]
 \arrow[c`d;pr_1]
 \arrow[d`e;]
\efig$$
where $H\ltimes^F\Rr$ is the twisted semi-direct product, whose 
multiplication is given by the formula
$$(h_0,\lambda_0)\odot_F(h_1,\lambda_1):=(h_0h_1,\lambda_0+\lambda_1+F(h_0,h_1)), $$
where $(h_0,\lambda_0),(h_1,\lambda_1)\in H\times\Rr$.

\lem\label{easyExt}
If $d(F,f)=0$, then 
$$\bfig
 \morphism(0,0)<800,0>[\psi_f:G`H\ltimes^F\Rr;]
 \morphism(1100,0)/{|->}/<800,0>[:g`(i(g),f(g,1));]
\efig$$
defines a crossed module for the action $g^{(h,\lambda)}:=g^h$. \endlem

\pf
$\psi_f$ is a Lie group homomorphism: 
$$\mld \psi_f(g_0)\odot_F\psi_f(g_1)& =(i(g_0),f(g_0,1))\odot_F(i(g_1),f(g_1,1))\\=(i(g_0)i(g_1),f(g_0,1)+f(g_1,1)+F(i(g_0),i(g_1)))\\=(i(g_0g_1),f(g_0g_1,1))=\psi_f(g_0g_1), $$
where the third equality follows from Eq. (\ref{singledOut}). Due to the 
independence of the variable in $\Rr$, thus defined, the action is still 
a right action by automorphisms and verifies the Peiffer identity. As for 
the equivariance of $\psi_f$, on the one hand we have got
$$\psi_f(g^{(h,\lambda)})=(i(g^h),f(g^h,1)), $$
while on the other, 
$$\mld (h,\lambda)^{-1}\odot_F\psi_f(g)\odot_F(h,\lambda)& =(h^{-1},-\lambda-F(h^{-1},h))\odot_F(i(g),f(g,1))\odot_F(h,\lambda) \\=(h^{-1}i(g),-\lambda-F(h^{-1},h)+f(g,1)+F(h^{-1},i(g)))\odot_F(h,\lambda)\\=(h^{-1}i(g)h,-F(h^{-1},h)+f(g,1)+F(h^{-1},i(g))+F(h^{-1}i(g),h)) $$
The first entries coincide because $i$ is the structural morphism of a 
crossed module. Evaluating $((1,h^{-1}),(g,1))$ and $((g,h^{-1}),(1,h))$ 
in Eq. (\ref{mixed}), one gets respectively
$$f(g,1)+F(h^{-1},i(g))=f(g,h^{-1}), $$
and
$$f(g,h^{-1})-F(h^{-1},h)+F(h^{-1}i(g),h)=f(g^h,1); $$
which combined, imply the result.

\epf

As a consequence Lemma \ref{easyExt}, there is a short exact sequence of 
Lie $2$-groups that we write using their associated crossed modules
$$\bfig
\square/>``>`>/<525,500>[1`1`1`\Rr;```]
\square(525,0)/>``>`>/<750,500>[1`G`\Rr`H\ltimes^F\Rr;```]
\square(1275,0)/>``>`>/<750,500>[G`G`H\ltimes^F\Rr`H;Id_G`\psi_f`i`pr_1]
\square(2025,0)/>```>/<525,500>[G`1`H`1;```]
\efig$$

\lem\label{isoEasyExt}
Let $(F,f),(F',f')\in\Omega^2_{tot}(\G)$ be a pair of cohomologous 
$2$-cocycles. Then the induced extensions of Lemma \ref{easyExt} are 
isomorphic. \endlem

\pf
Let $\phi\in\Omega^1_{tot}(\G)=C(H)$ be such that $(F,f)-(F',f')=d\phi=(\delta\phi,\partial\phi)$. In particular, $F$ and $F'$ are cohomologous cocycles; thus, the \emph{object} extensions are isomorphic via
$$\bfig
 \morphism(0,0)<800,0>[\alpha:H\ltimes^F\Rr`H\ltimes^{F'}\Rr;]
 \morphism(1200,0)/{|->}/<800,0>[:(h,\lambda)`(h,\lambda+\phi(h)).;]
\efig$$
We claim that $\alpha$, together with the identity of $G$ induce the 
claimed isomorphism between the extensions. Indeed, using the notation of 
Lemma \ref{easyExt},
$$\mld \alpha(\psi_f(g))& =\alpha(i(g),f(g,1))=(i(g),f(g,1)+\phi(i(g)))\\=(i(g),f'(g,1))=\psi_{f'}(g). $$
Also, trivially, $Id_G(g^{(h,\lambda)})=Id_G(g)^{\alpha(h,\lambda)}$, 
thus finishing the proof.

\epf

Lemmas \ref{easyExt} and \ref{isoEasyExt} should be taken as motivation 
to look for a complex whose cohomology classifies extensions starting 
from the bisimplicial structure naturally associated with a Lie $2$-group. 
They should be interpreted as the ``trivial coefficients'' case, thus 
prompting us to define a representation of a Lie $2$-group in hopes to 
classify extensions by more general $2$-vector spaces.

\subsection{Representations of Lie 2-group}\label{subsec-rep}
The General Linear Lie $2$-group \cite{IntSubLin2} is the Lie $2$-group 
which plays the r\^ole of space of automorphisms of a $2$-vector space. 
This $2$-group can be traced back at least to Norrie's thesis \cite{Norrie}, 
where it is called the {\em actor crossed module} of the $2$-vector space
regarded as an (abelian) $2$-group. The domain of the crossed module of the 
actor is the space of regular derivations -referred to as the Whitehead 
group- and the codomain is the space of automorphisms of the $2$-group. 
The associated $2$-group via Remark \ref{Equiv2Gp:Mods} can be identified 
with the category of linear invertible functors and natural homomorphic 
transformations further endowed with the horizontal composition of natural 
transformations which yields a group operation. We recall its structure 
for reference: Let $W\to^\phi V$ be a $2$-vector space. Then, the space 
of objects of its General Linear Lie $2$-group is the subgroup of 
invertible self functors
$$GL(\phi)_0 =\lbrace (F,f)\in GL(W)\times GL(V) : \phi\circ F=f\circ\phi\rbrace . $$
The Whitehead group of $W\to^\phi V$ is given by 
$$GL(\phi)_1 =\lbrace A\in Hom(V,W):(I+A\phi ,I+\phi A)\in GL(W)\times GL(V)\rbrace , $$
endowed with the operation
\begin{equation}\label{gpStr}
A_1\odot A_2 := A_1+A_2+A_1\phi A_2,
\end{equation} 
for which the identity element is the $0$ map, and inverses are given by 
either
$$A^{\dagger}=-A(I+\phi A)^{-1}=-(I+A\phi )^{-1}A. $$
We write $\dagger$ instead of $-1$ to avoid any possible overlap of 
notation with the actual inverse of a matrix. The crossed module map
$$\bfig
\morphism<650,0>[GL(\phi)_1`GL(\phi)_0.;\Delta]
\efig$$
is given by
$$\Delta A=(I+A\phi ,I+\phi A) $$
for $A\in GL(\phi)_1$. This is well defined since by definition it takes 
values in $GL(W)\times GL(V)$, and 
$$\phi(I+A\phi)= \phi +\phi A\phi = (I+\phi A)\phi . $$
Concluding, the right action of $GL(\phi)_0$ on $GL(\phi)_1$ is given by
\begin{equation}\label{act}
A^{(F,f)}=F^{-1}Af. 
\end{equation}

\proposition\label{GLPhi}\cite{Norrie, IntSubLin2}
Along with the group structure (\ref{gpStr}) and the action (\ref{act}),
$$\bfig
\morphism<650,0>[GL(\phi)_1`GL(\phi)_0;\Delta]
\efig$$ 
is a crossed module of Lie groups. \endproposition

In the sequel, we write $GL(\phi)$ for the crossed module of 
Proposition \ref{GLPhi} and its associated Lie $2$-group as well. 

Definition \ref{defRep} is an instance of the actions of \cite{Norrie} in 
the particular case when the $2$-group that is being acted on is abelian 
and simply connected. By definition, if $G\to^i H$ is the crossed module 
associated with $\G$ via Remark \ref{Equiv2Gp:Mods}, the $2$-representation 
(\ref{2rep}) consists of a Lie group homomorphism 
$$\rho_0 :H\to<350> GL(\phi)_0\leq GL(W)\times GL(V), $$
at the level of objects, which amounts to two Lie group representations 
$\rho_0^1 :H\to GL(W)$, $\rho_0^0 :H\to GL(V)$ fitting in
$$\bfig
\square[W`W`V`V,;(\rho_0^1)_h`\phi`\phi`(\rho_0^0)_h]
\efig$$ 
for all $h\in H$, {\em i.e.}, 
\begin{eqnarray}\label{0}
\rho_0^0(h)\circ\phi=\phi\circ\rho_0^1(h).
\end{eqnarray}
At the level of arrows, 
$$\rho_1 :G\to<350> GL(\phi)_1\leq Hom(V,W), $$
is a Lie group homomorphism if and only if
\begin{eqnarray}\label{1}
\rho_1(g_0g_1)=\rho_1(g_0)+\rho_1(g_1)+\rho_1(g_0)\circ\phi\circ\rho_1(g_1) 
\end{eqnarray}
for all $g_0,g_1\in G$. The compatibility between the rest of the 
crossed module structures is encoded in the following relations: 
\begin{eqnarray}\label{2}
\rho^0_0(i(g))=I+\phi\circ\rho_1(g), & \quad\rho^1_0(i(g))=I+\rho_1(g)\circ\phi 
\end{eqnarray}
for all $g\in G$, and 
\begin{eqnarray}\label{3}
\rho_1(g^h)=\rho_0^1(h)^{-1}\rho_1(g)\rho_0^0(h) 
\end{eqnarray}
for all $h\in H$, $g\in G$.

\example\label{trivRep}
Trivial representations. If $(\rho_1,\rho_0)\equiv(0,I)$, the defining 
equations for a $2$-representation are trivially satisfied. \endexample

\example\label{unitGpRep}
Usual Lie group representations. Letting $W=(0)$, a $2$-representation is 
ultimately equivalent to a representation of $H/i(G)$ on $V$. More 
precisely, the $2$-representation is defined by a single representation 
of $H$ on $V$ that vanishes along $i(G)$. In particular, a Lie group 
representation defines a $2$-representation of a unit Lie $2$-group on a 
unit $2$-vector space. 

Analogously, if $V=(0)$, a $2$-representation is ultimately equivalent to 
a representation of the orbit space on $W$. \endexample

\rem\label{OldLitRep}
The class of examples in Example \ref{unitGpRep} are referred to in the 
literature ({\em e.g.} \cite{GrandJLadra}) as $\G$-{\em modules}, though 
generalized to allow other abelian groups besides vector spaces. 
\endrem

\example\label{AdRep}
The adjoint representation: Let $\gg\to^\mu\hh$ be the Lie $2$-algebra of 
the Lie $2$-group $G\to^iH$, then we have 
$$\bfig
\square[G`GL(\mu)_1`H`GL(\mu)_0;Ad_1```Ad_0]
\efig$$
where $(Ad_0)_h=((-)^{h^{-1}},Ad_{h^{-1}})$, and $(Ad_1)_g=d_e({}_g\wedge)$ with
$${}_g\wedge : H\to G:h\to/|->/ g(g^{-1})^{h^{-1}}. $$ \endexample 

\rem\label{AdAppces}
Example \ref{AdRep} appears in \cite{IntSubLin2} and is also the derivative 
of the {\em canonical} action of \cite{Norrie} which is a generalized 
conjugation.
\endrem

Already in \cite{Norrie}, it is explained that given a $2$-representation 
$\rho$ of $\G$ on $\mathbbm{V}$, one can build a semi-direct product 
$2$-group $\G_\rho\ltimes\mathbbm{V}$. With the notation conventions of  
this section, the crossed module of the semi-direct product is 
\begin{equation}
G{}_{\rho_0^1\circ i}\ltimes W \to<350>^{i\times\phi} H{}_{\rho_0^0}\ltimes V,
\end{equation}\label{semiDirProd}
together with the right action given by
$$(g,w)^{(h,w)}=(g^h,\rho_0^1(h)^{-1}(w+\rho_1(g)v)) $$
for $(h,w)\in H\times V$ and $(g,w)\in G\times W$. 

The Lie group of arrows of $\G_\rho\ltimes\mathbbm{V}$,
$$(G{}_{\rho_0^1\circ i}\ltimes W)\rtimes(H{}_{\rho_0^0}\ltimes V), $$
is isomorphic to a semi-direct product of the Lie group $\G$ and 
$W\oplus V$ with respect to the honest representation that is the content 
of the following proposition. 

\proposition\label{honestGpRep}
Given a representation $2$-representation $\rho :\G\to GL(\phi )$, there 
is an honest representation 
$$\bar{\rho}:G\rtimes H\to<350> GL(W\oplus V):(g,h) \to/|->/<350>\pmatrix{
    \rho_0^1(hi(g)) & \rho_0^1(h)\rho_1(g) \cr
    0               & \rho_0^0(h) 
} $$ \endproposition

\pf
Consider the product
$$\mld \bar{\rho}(g_0,h_0)\bar{\rho}(g_1,h_1) & = \pmatrix{
                                               \rho_0^1(h_0i(g_0)) & \rho_0^1(h_0)\rho_1(g_0) \cr
                                               0                   & \rho_0^0(h_0) 
                                           }\pmatrix{
                                               \rho_0^1(h_1i(g_1)) & \rho_0^1(h_1)\rho_1(g_1) \cr
                                               0                   & \rho_0^0(h_1) 
                                           } \\ =\pmatrix{
 \rho_0^1(h_0i(g_0))\rho_0^1(h_1i(g_1)) & \rho_0^1(h_0i(g_0))\rho_0^1(h_1)\rho_1(g_1) +\rho_0^1(h_0)\rho_1(g_0)\rho_0^0(h_1) \cr
                               0            & \rho_0^0(h_0)\rho_0^0(h_1) 
                                           }. $$
The bottom row agrees with the bottom row of 
$\bar{\rho}((g_0,h_0)\vJoin (g_1,h_1))=\bar{\rho}(g_0^{h_1}g_1,h_0h_1)$ 
because $\rho_0^0$ is a group homomorphism. The first entries of the top 
row coincide too, put simply, because the target is a group homomorphism 
as well:
$$\rho_0^1(h_0h_1i(g_0^{h_1}g_1))=\rho_0^1(h_0h_1h_1^{-1}i(g_0)h_1i(g_1))=\rho_0^1(h_0i(g_0))\rho_0^1(h_1i(g_1)). $$
For the remaining entries, we use Eq.'s (\ref{1}) and (\ref{3}) to 
compute
$$\mld \rho_0^1(h_0h_1)\rho_1(g_0^{h_1}g_1) & = \rho_0^1(h_0h_1)(\rho_1(g_0^{h_1})(I+\phi\rho_1(g_1))+\rho_1(g_1)) \\
                                      = \rho_0^1(h_0h_1)\rho_0^1(h_1)^{-1}\rho_1(g_0)\rho_0^0(h_1)(I+\phi\rho_1(g_1))+\rho_0^1(h_0h_1)\rho_1(g_1) \\
                                      = \rho_0^1(h_0)\rho_1(g_0)\rho_0^0(h_1)+\rho_0^1(h_0)\rho_1(g_0)\phi\rho_0^1(h_1)\rho_1(g_1)+\rho_0^1(h_0)\rho_0^1(h_1)\rho_1(g_1) \\
                                      = \rho_0^1(h_0)\rho_1(g_0)\rho_0^0(h_1)+\rho_0^1(h_0)(\rho_1(g_0)\phi +I)\rho_0^1(h_1)\rho_1(g_1), $$
and the result follows from the second relation in Eq. (\ref{2}).

\epf

\rem\label{repUpToH}
Forgetting for the time being the Lie group structure, 
$$\bfig
\square/@{>}@<3pt>`>`>`@{>}@<3pt>/<650,500>[\G{}_{\bar{\rho}}\ltimes\mathbbm{V}`H{}_{\rho_0^0}\ltimes V`\G`H;```]
\square/@{>}@<-3pt>```@{>}@<-3pt>/<650,500>[\G{}_{\bar{\rho}}\ltimes\mathbbm{V}`H{}_{\rho_0^0}\ltimes V`\G`H;```]
\efig$$
has got the structure of a VB-groupoid \cite{BCD}. Thus, there is an 
associated representation up to homotopy (see \cite{VB&Reps}). In fact, 
since all vector bundles are trivial, there is an obvious splitting of 
the core sequence 
$$(0)\to<350> \G\times W\to<350>\G{}_{\bar{\rho}}\ltimes\mathbbm{V}\to<350>\G\times V\to<350> (0) $$
given by
$$\sigma:\G\times V\to<350>\G{}_{\bar{\rho}}\ltimes\mathbbm{V}:(g,h;v)\to/|->/<350> \sigma_{(g,h)}(h,v):=(g,h;0,v), $$
which verifies 
$\sigma_{u(h)}(h,v)=\sigma_{(1,h)}(h,v)=(1,h;0,v)=\hat{u}(h,v)$. Here, we 
use $\hat{\cdot}$ to refer to the structural maps of the top groupoid. 
$\sigma$ defines a canonical representation up to homotopy 
$(\varrho,\Delta^V,\Delta^W,\Omega)$ associated to the $2$-representation, 
where
$$\varrho :H\times W\to<250> H\times V:(h,w)\to/|->/<250>\hat{t}(1,h;w,0)=(h,\rho_0^0(h)\phi(w)), $$
the quasi-actions of $\G\cong G\rtimes H$ on $H\times V$ and $H\times W$ 
are respectively
$$\Delta^V_{(g,h)}(h,v) = \hat{t}(\sigma_{(g,h)}(h,v))=(hi(g),v), $$
$$\Delta^W_{(g,h)}(h,w) = \sigma_{(g,h)}(\varrho(h,w))\hat{\Join}(1,h;w,0)\hat{\Join}(g^{-1},hi(g);0,0)=(hi(g);\rho_0^1(i(g))^{-1}w), $$
and the curvature form 
$\Omega\in\Gamma(s_2^*(H\times V^*)\otimes t_2^*(H\times W))$ at 
$(g_1,g_2,h)\in G^2\times H\cong \G_2$ is 
$$\Omega_{(g_1,g_2,h)}(v)=\Big{(}\sigma_{(g_2g_1,h)}(h,v)-\sigma_{(g_1,hi(g_2))}\big{(}\Delta^V_{(g_2,h)}(h,v)\big{)}\hat{\Join}\sigma_{(g_2,h)}(h,v)\Big{)}\hat{\Join} 0_{(g_2g_1,h)^{-1}}=0_{hi(g_2g_1)}. $$
Since $\Omega$ is identically zero, the quasi-actions define actual 
representations of the Lie $2$-group $\G$ on the corresponding vector bundles. \endrem

We close this section by proving that splitting an abstract extension of 
Lie $2$-groups induces a $2$-representation in the sense of Definition 
\ref{defRep}.

\definition\label{2-extensions}
An {\em extension} of the Lie $2$-group $G\to^iH$ by the $2$-vector space 
$W\to^\phi V$ is a Lie $2$-group $E_1\to^\epsilon E_0$ that fits in 
$$\bfig
\square/>``>`>/<525,500>[1`W`1`V;```]
\square(525,0)/>``>`>/<750,500>[W`E_1`V`E_0;j_1`\phi``j_0]
\square(1275,0)/>``>`>/<750,500>[E_1`G`E_0`H;\pi_1`\epsilon`i`\pi_0]
\square(2025,0)/>```>/<525,500>[G`1`H`1,;```]
\efig$$
{\em i.e.}, where the top and bottom rows are short exact sequences and 
the squares are maps of Lie $2$-groups. \enddefinition

\proposition\label{Ind2GpRep}
Given an extension of the Lie $2$-group $G\to^iH$ by the $2$-vector space 
$W\to^\phi V$ and a smooth splitting $\sigma$,
$$\bfig
\square/>``>`>/<525,500>[1`W`1`V;```]
\square(525,0)/>``>`>/<750,500>[W`E_1`V`E_0;j_1`\phi``j_0]
\square(1275,0)|blrb|/>``>`>/<750,500>[E_1`G`E_0`H;\pi_1`\epsilon`i`\pi_0]
\square(2025,0)/>```>/<525,500>[G`1`H`1,;```]
\square(1275,0)|alra|/{@{<-}@/^1em/}```{@{<-}@/^1em/}/<750,500>[E_1`G`E_0`H;\sigma_1```\sigma_0]
\efig$$
there is an induced $2$-representation $\rho^\epsilon_\sigma:\G\to GL(\phi)$ 
given by
$$\rho^0_0(h)v=\sigma_0(h)v\sigma_0(h)^{-1} \qquad\qquad \rho^1_0(h)w=w^{\sigma_0(h)^{-1}} $$
$$\rho_1(g)v =\sigma_1(g)^v\sigma_1(g)^{-1}, $$
for $h\in H$, $v\in V$, $w\in W$ and $g\in G$. \endproposition
 
The proof of Proposition \ref{Ind2GpRep} follows easily after using the 
splitting to write $E_1$ and $E_0$ as semi-direct products. We thus 
postpone it to the end of the section in order to introduce the necessary 
notation.

We regard injective maps as inclusions; in so, $\phi=\epsilon\vert_W$. 
Given an extension as in the statement of Proposition \ref{Ind2GpRep}, 
one uses the splitting to get the diffeomorphisms $H\times V\cong E_0$ 
and $G\times W \cong E_1$ given respectively by
$$(z,a)\to/|->/<350> a\sigma_k (z) $$
with inverse 
$$e\to/|->/<350> (\pi_k (e),e\sigma_k (\pi_k (e))^{-1}), $$
for $k\in\lbrace 0,1\rbrace$. We recall that one can use these 
diffeomorphisms to transfer the group structure, thus getting 
\begin{equation}\label{twstdProd}
(h_0,v_0)\cdot(h_1,v_1):=(h_0h_1,v_0+\rho_0^0(h_0)v_1+\omega_0(h_0,h_1)),
\end{equation}  
where $(h_0,v_0),(h_1,v_1)\in H\times V$, $\rho_0^0$ is defined as 
in Proposition \ref{Ind2GpRep} and 
$\omega_0(h_0,h_1):=\sigma_0(h_0)\sigma_0(h_1)\sigma_0(h_0h_1)^{-1}\in V$. 
This is the usual twisted semi-direct product $H{}_{\rho^0_0}\ltimes^{\omega_0} V$ 
from the theory of Lie group extensions. Conversely, the operation defined 
by Eq. (\ref{twstdProd}) with $\omega_0\in C(H^2,V)$ is an associative 
product if and only if $\omega_0$ is a $2$-cocycle for the Lie group 
cohomology of $H$ with values in $\rho^0_0$ (see Eq. (\ref{LGpdDifferential})). 

An identical reasoning implies there is an isomorphism of Lie groups 
$E_1\cong G{}_{\rho^1_0\circ i}\ltimes^{\omega_1}W$, with $\rho_0^1$ 
defined as in Proposition \ref{Ind2GpRep} and 
$\omega_1(g_0,g_1)=\sigma_1(g_0)\sigma_1(g_1)\sigma_1(g_0g_1)^{-1}$. 

We rewrite the rest of the crossed module structure of the extension 
using these trivializations. The homomorphism $\epsilon$ gets rewritten 
as
$$G\times W\to<250>H\times V:(g,w) \to/|->/<250>(i(x),\phi(w)+\varphi(g)), $$
where $\varphi(g):=\epsilon(\sigma_1(g))\sigma_0(i(g))^{-1}$, and the 
action of $(h,v)\in H\times V$ on $(g,w)\in G\times W$ as
$$(g,w)^{(h,v)}:=(g^h,\rho^1_0(h)^{-1}(w+\rho_1(g)v)+\alpha(h;g)), $$
where $\rho_0^1$ and $\rho_1$ are defined as in Proposition \ref{Ind2GpRep}, 
and $\alpha(h;g):=\sigma_1(g)^{\sigma_0(h)}\sigma_1\big{(}g^h\big{)}^{-1}$.

\rem\label{vanishingCocycle} Notice that, in case one can take the 
splitting $\sigma$ in Proposition \ref{Ind2GpRep} to be a crossed module 
map, $\omega_0$, $\omega_1$, $\varphi$ and $\alpha$ vanish identically 
and one recovers the semi-direct product structure (\ref{semiDirProd}). \endrem 

\pf[of Proposition \ref{Ind2GpRep}]
We make the necessary computations to prove that 
$\rho^{\epsilon}_{\sigma}$ is a $2$-representation.
\begin{itemize}
\item Well-defined: We use the exactness of the sequences to see that the 
maps land where they are supposed to. Let $h\in H$, $v\in V$, $w\in W$ 
and $g\in G$, then
$$\pi_0\big{(}\sigma_0(h)v\sigma_0(h)^{-1}\big{)}=\pi_0(\sigma_0(h))\pi_0(v)\pi_0\big{(}\sigma_0(h)^{-1}\big{)}=h1h^{-1}=1 \Longrightarrow \rho_0^0(h)v\in V, $$
$$\pi_1\big{(}w^{\sigma_0(h)^{-1}}\big{)}=\pi_1(w)^{\pi_0\big{(}\sigma_0(h)^{-1}\big{)}}=1^{h^{-1}}=1 \Longrightarrow \rho_0^1(h)w\in W, $$
$$\pi_1\big{(}\sigma_1(g)^v\sigma_1(g)^{-1}\big{)}=\pi_1(\sigma_1(g))^{\pi_1(v)}\pi_1\big{(}\sigma_1(g)^{-1}\big{)}=g^1g^{-1}=1  \Longrightarrow \rho_1(g)v\in W. $$
These components are smooth and linear. Further, thus defined, 
$\rho_0^0(h)\circ\phi=\phi\circ\rho_0^1(h)$ for each $h\in H$. Indeed, 
let $w\in W$, then
$$\rho_0^0(h)(\phi(w))=\sigma_0(h)\epsilon(w)\sigma_0(h)^{-1}=\epsilon\big{(}w^{\sigma_0(h)^{-1}}\big{)} =\phi(\rho_0^1(h)w). $$
Thus, for all $h\in H$, $\rho_0(h):=(\rho_0^0(h),\rho_0^1(h))\in GL(\phi)_0$.

\item $\rho^0_0$ is a Lie group representation: Let $h_1,h_2\in H$ and 
$v\in V$, then
$$\mld \rho^0_0(h_1)\rho^0_0(h_2)v & =\sigma_0(h_1)\sigma_0(h_2)v\sigma_0(h_2)^{-1}\sigma_0(h_1)^{-1} \\
							         =\omega_0(h_1,h_2)(\rho^0_0(h_1h_2)v)\omega_0(h_1,h_2)^{-1}. $$
Since both $\omega_0(h_1,h_2),\rho^0_0(h_1h_2)v\in V$, the conjugation is 
trivial and the claim follows.
\item $\rho_0^1$ is a Lie group representation: Let $h_1,h_2\in H$ and 
$w\in W$, then
$$\mld \rho^1_0(h_1)\rho^1_0(h_2)w & =w^{\sigma_0(h_2)^{-1}\sigma_0(h_1)^{-1}} \\
							 		 =(\rho^1_0(h_1h_2)w)^{\omega_0(h_1,h_2)^{-1}}. $$
Since $\omega_0(h_1,h_2)\in V$ and $\rho^1_0(h_1h_2)w\in W$, the action 
is trivial and the claim follows.
\item $\rho_1$ is a Lie group homomorphism: Let $g_1,g_2\in G$ and 
$v\in V$, then
$$\mld \rho_1(g_1)\odot\rho_1(g_2)v & =(\rho_1(g_1)+\rho_1(g_2)+\rho_1(g_1)\phi\rho_1(g_2))v \\
							  =\sigma_1(g_1)^{v}\sigma_1(g_1)^{-1}\sigma_1(g_2)^{v}\sigma_1(g_2)^{-1}\rho_1(g_1)\epsilon(\sigma_1(g_2)^{v}\sigma_1(g_2)^{-1}) \\
							  =\sigma_1(g_1)^{v}\sigma_1(g_1)^{-1}\sigma_1(g_2)^{v}\sigma_1(g_2)^{-1}\sigma_1(g_1)^{\epsilon\big{(}\sigma_1(g_2)^{v}\sigma_1(g_2)^{-1}\big{)}}\sigma_1(g_1)^{-1}. $$
Using Peiffer equation, we get 
$$\sigma_1(g_1)^{\epsilon\big{(}\sigma_1(g_2)^{v}\sigma_1(g_2)^{-1}\big{)}}=\big{(}\sigma_1(g_2)^{v}\sigma_1(g_2)^{-1}\big{)}^{-1}\sigma_1(g_1)\big{(}\sigma_1(g_2)^{v}\sigma_1(g_2)^{-1}\big{)}; $$
thus yielding,
$$\mld \rho_1(g_1)\odot\rho_1(g_2)v & =\sigma_1(g_1)^{v}\sigma_1(g_2)^{v}\sigma_1(g_2)^{-1}\sigma_1(g_1)^{-1} \\
							  =\sigma_1(g_1)^{v}\sigma_1(g_2)^{v}\big{(}\sigma_1(g_1g_2)^{v}\big{)}^{-1}\sigma_1(g_1g_2)^{v}\sigma_1(g_1g_2)^{-1}\sigma_1(g_1g_2)\sigma_1(g_2)^{-1}\sigma_1(g_1)^{-1} \\
							  =(\omega_1(g_1,g_2))^v(\rho_1(g_1g_2)v)\omega_1(g_1,g_2)^{-1}. $$
Since both $\omega_1(g_1,g_2),\rho_1(g_1g_2)v\in W$, both the action and 
the conjugation are trivial and the claim follows.       
\item $\rho_0\circ i=\Delta\circ\rho_1$: For $g\in G$, this equation 
breaks into two components,
$$\rho_0^0(i(g))=I+\phi\circ\rho_1(g)\in GL(V) \qquad\textnormal{and}\qquad \rho_0^1(i(g))=I+\rho_1(g)\circ\phi\in GL(W). $$
Let $v\in V$ and $w\in W$, then
$$\mld (I+\phi\rho_1(g))v & =v\epsilon\big{(}\sigma_1(g)^v\sigma_1(g)^{-1}\big{)} =\epsilon(\sigma_1(g))v\epsilon(\sigma_1(g))^{-1} \\      
                    =\varphi(g)(\rho_0^0(i(g))v)\varphi(g)^{-1}=\rho_0^0(i(g))v, $$
and
$$\mld (I+\rho_1(g)\phi)w & =w\sigma_1(g)^{\epsilon(w)}\sigma_1(g)^{-1}=\sigma_1(g)w\sigma_1(g)^{-1} \\
                    =w^{\epsilon\big{(}\sigma_1(g)^{-1}\big{)}} =(\rho_0^1(i(g))w)^{\varphi(g)^{-1}}=\rho_0^1(i(g))w, $$
where the last equality in each sequence follows from $\varphi(g)\in V$.
\item $\rho_1$ respects the actions: Let $g\in G$, $h\in H$ and $v\in V$, 
then
$$\mld \rho_1(g)^{\rho_0(h)}v & =\rho_0^1(h)^{-1}\rho_1(g)\rho_0^0(h)v =\big{(}\sigma_1(g)^{\sigma_0(h)v\sigma_0(h)^{-1}}\sigma_1(g)^{-1}\big{)}^{\sigma_0(h)} \\  
                        =\big{(}\sigma_1(g)^{\sigma_0(h)}\sigma_1(g^h)^{-1}\big{)}^v\sigma_1\big{(}g^h\big{)}^v\big{(}\sigma_1(g)^{\sigma_0(h)}\big{)}^{-1}=\alpha(h;g)^v\big{(}\rho_1(g^h)v\big{)}\alpha(h;g)^{-1}. $$
Since both $\alpha(h;g),\rho_1(g^h)v\in W$, both the action and the 
conjugation are trivial and the claim follows.
\end{itemize}
\epf

\section{The grid, the truncated complex and its cohomology}\label{sec-tetrahedralCx}

In this section, inspired by the triple complex associated to a Lie 
$2$-algebra \cite{Lie2AlgCoh} and by the concurrence of the double 
cohomology with the cohomology of \cite{Ellis}, we introduce the 
{\em grid} of complexes of Lie $2$-group cochains with values in a 
$2$-representation. We establish notation and specify all groupoids and 
representations appearing in the three dimensional grid. Since we define 
the complex in the $r$-direction to start differently in $0$th degree, we 
prove that this fits into a complex accordingly. We expose how this grid 
fails to be a triple complex. Specifically, we study the square of the 
total differential degree by degree and define the {\em difference maps} 
as the commutator of the non-commuting differentials, thus producing a 
complex up to degree $3$. We study the cohomology of the resulting 
truncated complex and show that its second cohomology classifies 
abelian extensions.

\subsection{The background grid}\label{Grid}

Throughout, let $\G$ be a Lie $2$-group with associated crossed module 
$G\to^iH$ and let $\rho$ be a $2$-representation of $\G$ on the $2$-vector 
space $W\to^\phi V$. Using the notation conventions laid down in Subection 
\ref{subsec-rep}, we think of $\rho$ as a triple 
$(\rho_0^0,\rho_0^1;\rho_1)$. 
Also, we take the isomorphisms of Remarks \ref{Equiv2Gp:Mods} and 
\ref{2GpNerve} to be fixed and we abuse notation and often treat them as 
equalities.

\rem\label{initialSce}
In the sequel, we define the grid using a series of groupoids and 
groupoid representations that involve taking pull-backs along the final 
target map (cf. Remark \ref{2GpFinalTarget}). One could also define an 
equivalent structure by pulling-back along the ``initial source'' map 
$s_p:\G_p\to H:(\gamma_1,...,\gamma_p)\to/|->/s(\gamma_p)$. There is no 
economy in working with either one, one necessarily pays a computational 
price somewhere. \endrem

Let $C^{p,q}_r(\G,\phi)$ be defined by Eq.'s (\ref{3dimLat}) and 
(\ref{3dimLat0}). This three dimensional lattice of vector spaces comes 
together with a grid of complexes of groupoid cochains (cf. 
(\ref{LGpdDifferential}) and (\ref{RGpdDifferential})). 

\subsubsection{The $p$-direction}\label{p-dir} When $q=0$, one has got 
the trivial complexes
$$C^{0,0}_0(\G,\phi)=V\to<350>^{\partial=0} C^{1,0}_0(\G,\phi)=V\to<350>^{\partial=Id_V} C^{2,0}_0(\G,\phi)=V\to<350>^{\partial=0} C^{3,0}_0(\G,\phi)=V\to\cdots $$
when $r=0$, and 
$$C(G^r,W)\to<350>^{\partial=0}C(G^r,W)\to<350>^{\partial=Id_{C(G^r,W)}}C(G^r,W)\to<350>^{\partial=0} C(G^r,W)\to<350>\cdots $$
otherwise. When $q\neq 0$ and $r=0$, the complex
$$C(H^q,V)\to<350>^{\partial} C(\G^q,V)\to<350>^{\partial} C(\G_2^q,V)\to<350>^{\partial} C(\G_3^q,V)\to<350>\cdots $$
is the cochain complex of the product groupoid 
$$\G^q\two H^q$$
with respect to the trivial representation on the vector bundle 
$H^q\times V\to^{pr_1} H^q$. For any other value of $r$, the complex
$$C(H^q\times G^r,W)\to<350>^{\partial} C(\G^q\times G^r,W)\to<350>^{\partial} C(\G_2^q\times G^r,W)\to<350>^{\partial} C(\G_3^q\times G^r,W)\to<250>\cdots $$
is the cochain complex of the product groupoid 
$$\G^q\times G^r\two H^q\times G^r$$
with respect to the left representation on the trivial 
bundle $H^q\times G^r\times W\to^{pr_1}H^q\times G^r$,
\begin{equation}\label{p-rep}
(\gamma_1,...,\gamma_q;\vec{f})\cdot(h_1,...,h_q;\vec{f},w):=(h_1i(g_1),...,h_qi(g_q);\vec{f},\rho_0^1(i(pr_G(\gamma_1\vJoin ...\vJoin\gamma_q)))^{-1}w), 
\end{equation}
where $\gamma_k=(g_k,h_k)\in\G$ and $\vec{f}\in G^r$.

\rem\label{repUpToHom-onVandW} 
Observe that for $q=1$, the representations defining the complexes in the 
$p$-direction are those of the representation up to homotopy induced by 
the $2$-representation (see Remark \ref{repUpToH}). Indeed, for $r=0$, 
the representation coincides with $\Delta^V$, and for $r>0$, the 
representation coincides with the pull-back of $\Delta^W$ along the 
projection onto $\G\two H$. \endrem

\lem\label{p-cx}
Eq. (\ref{p-rep}) defines a representation. \endlem

Notice that this lemma is not straightforward, because the projection 
$pr_G$ is not in general a group homomorphism. In order to make cleaner 
computations, we introduce the following auxiliary straightforward lemma.

\lem\label{multiprods}
Let $\gamma_1,...,\gamma_q\in\G$. If $\gamma_k=(g_k,h_k)\in G\rtimes H$, 
then
$$\gamma_1\vJoin ...\vJoin\gamma_q=(g_1^{h_2...h_q}g_2^{h_3...h_q}...g_{q-2}^{h_{q-1}h_q}g_{q-1}^{h_q}g_q,h_1...h_q). $$ \endlem

\pf
By induction on $q$, for $q=2$ the formula is nothing but the definition 
of the product in $\G\cong G\rtimes H$. Now, suppose the equation holds 
for $q-1$ elements, then
$$\mld \gamma_1\vJoin ...\vJoin\gamma_q & =(\gamma_1\vJoin ...\vJoin\gamma_{q-1})\vJoin\gamma_q \\
 =(g_1^{h_2...h_{q-1}}g_2^{h_3...h_{q-1}}...g_{q-3}^{h_{q-2}h_{q-1}}g_{q-2}^{h_{q-1}}g_{q-1},h_1...h_{q-1})\vJoin (g_q,h_q) \\
 =((g_1^{h_2...h_{q-1}}g_2^{h_3...h_{q-1}}...g_{q-3}^{h_{q-2}h_{q-1}}g_{q-2}^{h_{q-1}}g_{q-1})^{h_q}g_q,h_1...h_{q-1}h_q), $$
and the result follows since the action of $H$ is by automorphisms.

\epf

\pf[of Lemma \ref{p-cx}]
Observe that the first coordinates in the right hand side of Eq. 
(\ref{p-rep}) are given by the target; hence, we just need to focus on 
the $W$ coordinate. 

Units act trivially: Let $w\in W$ and $\gamma_1,...,\gamma_q\in\G$ be 
such that $\gamma_k=(1,h_k)\in G\rtimes H$ for each $k$. Then, the $W$ 
coordinate in the right hand side of Eq. (\ref{p-rep}) reads
$$\mld \rho_0^1(i(pr_G(\gamma_1\vJoin ...\vJoin\gamma_q)))^{-1}w & =\rho_0^1(i(1^{h_2...h_q}1^{h_3...h_q}...1^{h_{q-1}h_q}1^{h_q}1))^{-1}w \\
	= \rho_0^1(i(1...1))^{-1}w =w. $$
	
Groupoid multiplication: The arrows 
$(\gamma'_1,...,\gamma'_q;(\vec{f})'),(\gamma_1,...,\gamma_q;\vec{f})\in\G^q\times G^r$ 
are composable if and only if $(\vec{f})'=\vec{f}$ and 
$(\gamma'_k,\gamma_k)\in\G^{(2)}$ for each $k$, or equivalently, if 
$\gamma'_k=(g'_k,h_ki(g_k))\in G\rtimes H$. For such a pair, the groupoid 
multiplication is given by 
$$(\gamma'_1,...,\gamma'_q;\vec{f})\cdot(\gamma_1,...,\gamma_q;\vec{f}):=\Big{(}\pmatrix{g_1g'_1 & ... & g_qg'_q \cr
h_1     & ... & h_q};\vec{f}\Big{)}. $$
The compatibility follows from the formula
$$\mld (g_1g'_1)^{h_2...h_q} & (g_2g'_2)^{h_3...h_q}...(g_{q-1}g'_{q-1})^{h_q}g_qg'_q = \\
               g_1^{h_2...h_q}g_2^{h_3...h_q}...g_{q-1}^{h_q}g_q(g'_1)^{h_2i(g_2)...h_qi(g_q)}(g'_2)^{h_3i(g_3)...h_qi(g_q)}...(g'_{q-1})^{h_qi(g_q)}g'_q . $$
We proceed by induction on $q$. For $q=2$,
$$(g_1g'_1)^{h_2}g_2g'_2 = g_1^{h_2}(g'_1)^{h_2}g_2g'_2 = g_1^{h_2}g_2g_2^{-1}(g'_1)^{h_2}g_2g'_2 = g_1^{h_2}g_2(g'_1)^{h_2i(g_2)}g'_2 . $$
Suppose now that the equation holds for $q-1$, then
$$\mld \prod_{k=1}^q (g_kg'_k)^{\prod_{j=k+1}^q h_j} & =\Bigg{(}\prod_{k=1}^{q-1} (g_kg'_k)^{\Big{(}\prod_{j=k+1}^{q-1} h_j\Big{)}h_q}\Bigg{)}g_qg'_q \\
	=\Bigg{(}\prod_{k=1}^{q-1} (g_kg'_k)^{\prod_{j=k+1}^{q-1} h_j}\Bigg{)}^{h_q}g_qg'_q \\
	=^{(I.H.)}\Bigg{(}\prod_{k=1}^{q-1} g_k^{\prod_{j=k+1}^{q-1} h_j}\prod_{k=1}^{q-1}(g'_k)^{\prod_{j=k+1}^{q-1} h_ji(g_j)}\Bigg{)}^{h_q}g_qg'_q \\
	=\Bigg{(}\prod_{k=1}^{q-1} g_k^{\prod_{j=k+1}^{q-1} h_j}\Bigg{)}^{h_q}g_qg_q^{-1}\Bigg{(}\prod_{k=1}^{q-1}(g'_k)^{\prod_{j=k+1}^{q-1} h_ji(g_j)}\Bigg{)}^{h_q}g_qg'_q \\
	=\prod_{k=1}^{q} g_k^{\prod_{j=k+1}^{q} h_j}\Bigg{(}\prod_{k=1}^{q-1}(g'_k)^{\prod_{j=k+1}^{q-1} h_ji(g_j)}\Bigg{)}^{h_qi(g_q)}g'_q , $$
which is precisely what we wanted.

\epf

\subsubsection{The $q$ direction}\label{q-dir} When $r=0$, the complex
$$V\to<350>^{\delta} C(\G_p,V)\to<350>^{\delta} C(\G_p^2,V)\to<350>^{\delta} C(\G_p^3,V)\to<350>\cdots $$
is the group complex of $\G_p$ with values in the pull-back of the 
representation $\rho_0^0$ along the final target map $t_p$; when 
$r\neq 0$, the complex 
$$C(G^r,W)\to<350>^{\delta} C(\G_p\times G^r,W)\to<350>^{\delta} C(\G_p^2\times G^r,W)\to<350>^{\delta} C(\G_p^3\times G^r,W)\to<350>\cdots $$
is the cochain complex of the (right!) transformation groupoid 
$$\G_p\ltimes G^r\two G^r$$
with respect to the right representation
\begin{equation}\label{q-Rep}
(g_1,...,g_r;w)\cdot(\gamma;g_1,...,g_r):=(g_1^{t_p(\gamma)},...,g_r^{t_p(\gamma)};\rho_0^1(t_p(\gamma))^{-1}w)
\end{equation}
on the trivial vector bundle $G^r\times W\to^{pr_1}G^r$, where 
$g_1,...,g_r\in G$, $\gamma\in\G_p$ and $w\in W$. It is obvious that Eq. 
(\ref{q-Rep}) defines a representation, as, on the one hand, is defined 
explicitly using the source map of the transformation groupoid and, on 
the other, it is given fibre-wise by the pull-back of a representation 
along a homomorphism.

When writing the groupoid differential, we use the shorthand 
$\rho^r_{\G_p}(\gamma;g)w$ instead of the lengthier Eq. (\ref{q-Rep}).
 
\subsubsection{The $r$-direction}\label{r-dir} When $q=0$, the complex
$$V\to<350>^{\delta'} C(G,W)\to<350>^{\delta_{(1)}} C(G^2,W)\to<350>^{\delta_{(1)}} C(G^3,W)\to<350>\cdots $$
is the group complex of $G$ with values in the pull-back of the 
representation $\rho_0^1$ along the crossed module homomorphism $i$, but 
for the $0$th degree; when $q\neq 0$, the complex 
$$C(\G_p^q,V)\to<350>^{\delta'} C(\G_p^q\times G,W)\to<350>^{\delta_{(1)}} C(\G_p^q\times G^2,W)\to<350>^{\delta_{(1)}} C(\G_p^q\times G^3,W)\to<350>\cdots $$
is, again except for the $0$th degree, the cochain complex of the Lie 
group bundle
$$\G_p^q\times G\two \G_p^q$$
with respect to the left representation
\begin{equation}\label{r-Rep}
(\gamma_1,...,\gamma_q;g)\cdot(\gamma_1,...,\gamma_q;w):=(\gamma_1,...,\gamma_q;\rho_0^1(i(g^{t_p(\gamma_1)...t_p(\gamma_q)}))w)
\end{equation}
on the trivial vector bundle $G_p^q\times W\to^{pr_1}G_p^q$, where 
$\gamma_1,...,\gamma_q\in\G_p$, $g\in G_p$ and $w\in W$. Eq. (\ref{r-Rep}) 
clearly defines a representation, as it is given fibre-wise by the 
pull-back of a representation along a homomorphism, namely the composition 
of the crossed module homomorphism $i$ with the crossed module action. 
Notice that right and left representations of a Lie group bundle coincide; 
hence, though Eq. (\ref{r-Rep}) could be taken as a right representation, 
we emphasize that it is a left representation as the formula for the 
differential of right and left representations differ 
(cf. (\ref{LGpdDifferential}) and (\ref{RGpdDifferential})). 

The missing maps $\delta':V\to C(G,W)$ and 
$\delta':C(\G_p^q,V)\to C(\G_p^q\times G,W)$ are defined respectively by
\begin{equation}\label{q=01st-r}
(\delta'v)(g):=\rho_1(g)v,
\end{equation}
for $v\in V$ and $g\in G$, and by
\begin{equation}\label{1st-r}
\delta'\omega(\gamma_1,...,\gamma_q;g)=\rho_0^1(t_p(\gamma_1)...t_p(\gamma_q))^{-1}\rho_1(g)\omega(\gamma_1,...,\gamma_q),
\end{equation}
for $\omega\in C(\G_p^q,V)$, $\gamma_1,...,\gamma_q\in\G_p$ and $g\in G$. 

The next two lemmas justify how, in spite of the replacements in $0$th 
degree, the complexes in the $r$-direction remain complexes.

\lem\label{q=0r-cx}
$$V\to<350>^{\delta'} C(G,W)\to<350>^{\delta_{(1)}} C(G^2,W), $$
where $\delta'$ is defined by Eq. (\ref{q=01st-r}), is a complex. \endlem

\pf
We prove that, for $v\in V$, $\delta_{(1)}\delta'v=0$. Let $g_0,g_1\in G$, 
then
$$\mld \delta_{(1)}(\delta'v)(g_0,g_1) & =\rho_0^1(i(g_0))(\delta'v)(g_1)-(\delta'v)(g_0g_1)+(\delta'v)(g_0) \\
									 =\rho_0^1(i(g_0))\rho_1(g_1)v-\rho_1(g_0g_1)v+\rho_1(g_0)v \\
									 =(I+\rho_1(g_0)\phi)\rho_1(g_1)v-\rho_1(g_0g_1)v+\rho_1(g_0)v, $$
which is zero due to Eq. (\ref{1}).

\epf

\lem\label{r-cx}
$$C(\G_p^q,V)\to<350>^{\delta'} C(\G_p^q\times G,W)\to<350>^{\delta_{(1)}} C(\G_p^q\times G^2,W), $$
where $\delta'$ is defined by Eq. (\ref{1st-r}), is a complex. \endlem

\pf
We prove that, for $\omega\in C(\G_p^q,V)$, $\delta_{(1)}\delta'\omega=0$. 
Let $\gamma_1,...,\gamma_q\in\G_p$ and $g_0,g_1\in G$, then
$$\mld \delta_{(1)} & \delta'\omega(\gamma_1,...,\gamma_q;g_0,g_1) \\
      =(\gamma_1,...,\gamma_q;g_0)\cdot\delta'\omega(\gamma_1,...,\gamma_q;g_1)-\delta'\omega(\gamma_1,...,\gamma_q;g_0g_1)+\delta'\omega(\gamma_1,...,\gamma_q;g_0) \\
      =\rho_0^1(i(g_0^{t_p(\gamma_1)...t_p(\gamma_q)}))\rho_0^1(t_p(\gamma_1)...t_p(\gamma_q))^{-1}\rho_1(g_1)\omega(\gamma_1,...,\gamma_q)+ \\
	  \quad-\rho_0^1(t_p(\gamma_1)...t_p(\gamma_q))^{-1}\rho_1(g_0g_1)\omega(\gamma_1,...,\gamma_q)+\rho_0^1(t_p(\gamma_1)...t_p(\gamma_q))^{-1}\rho_1(g_0)\omega(\gamma_1,...,\gamma_q) \\
	  =\rho_0^1(t_p(\gamma_1)...t_p(\gamma_q))^{-1}\Big{(}\rho_0^1(i(g_0))\rho_1(g_1)\omega(\gamma_1,...,\gamma_q)-\rho_1(g_0g_1)\omega(\gamma_1,...,\gamma_q)+\rho_1(g_0)\omega(\gamma_1,...,\gamma_q)\Big{)}, $$
which is zero due to Eq.'s (\ref{1}) and (\ref{2}).

\epf

Let 
\begin{equation}\label{totCx}
C_{tot}^{n}(\G,\phi)=\bigoplus_{p+q+r=n}C^{p,q}_r(\G,\phi).
\end{equation}
For expository purposes, we preliminarily define the total differential 
\begin{equation}\label{preDiff}
\nabla:=(-1)^p(\delta_{(1)}+\partial+(-1)^r\delta), 
\end{equation}
as though the grid were a triple complex. In the course of the remainder 
of this section, we study $\nabla^2$ degree by degree and conclude that, 
despite $\partial$ and $\delta$ fail to commute, one can add corrections 
to have a complex
\begin{equation}\label{trunc3}
C^0_{tot}(\G,\phi)\to<350>^\nabla C^1_{tot}(\G,\phi)\to<350>^\nabla C^2_{tot}(\G,\phi)\to<350>^\nabla C^3_{tot}(\G,\phi).
\end{equation}
Then, we move on to study the cohomology of (\ref{trunc3}). We postpone a 
more general study of the relations among the differentials in the grid 
until the next section.

\subsection{Degree 0}\label{H^0} By definition 
$(C_{tot}^{n}(\G,\phi),\nabla)$ is concentrated in nonnegative degrees, 
thereby defining a complex in degree $0$. Let 
$v\in C_{tot}^0(\G,\phi)=C^{0,0}_0(\G,\phi)=V$ and consider its 
differential
$$\nabla v=(\cancelto{0}{\partial v},\delta v,\delta'v)\in C^1(\G,\phi)=C^{1,0}_0(\G,\phi)\oplus C^{0,1}_0(\G,\phi)\oplus C^{0,0}_1(\G,\phi). $$ 
If $v$ is a $0$-cocycle, then 
$$\rho_0^0(h)v=v \qquad\textnormal{and}\qquad \rho_1(g)v=0 $$
for all $h\in H$ and all $g\in G$; therefore, 
$$H^0_\nabla(\G,\phi)=V^{\G}:=\lbrace v\in V:\bar{\rho}_{(g,h)}(0,v)=(0,v),\quad\forall (g,h)\in G\rtimes H\cong\G\rbrace , $$
where $\bar{\rho}$ is the honest representation of Proposition 
\ref{honestGpRep}.

\subsection{Degree 1}\label{H^1} A $1$-cochain $\lambda$ is a triple 
$(v,\lambda_0,\lambda_1)\in C_{tot}^1(\G,\phi)=V\oplus C(H,V)\oplus C(G,W)$ 
whose differential has six entries. Adopting the convention that 
$(\nabla\lambda)^{p,q}_r\in C^{p,q}_r(\G,\phi)$, 
$$(\nabla\lambda)^{0,2}_0 = \delta\lambda_0 $$
$$(\nabla\lambda)^{1,1}_0 = \partial\lambda_0-\delta v \qquad\qquad (\nabla\lambda)^{0,1}_1 = \delta'\lambda_0-\delta\lambda_1 $$
$$(\nabla\lambda)^{2,0}_0 =-\partial v=-v \qquad\qquad (\nabla\lambda)^{1,0}_1 =\partial\lambda _1-\delta'v=-\delta'v \qquad\qquad (\nabla\lambda)^{0,0}_2 = \delta_{(1)}\lambda_1 . $$

Let $v\in V$ and put $\lambda=\nabla v$. With the exception 
$(\nabla^2v)^{0,1}_1$ and $(\nabla^2v)^{1,1}_0$, all components of 
$\nabla^2v$ vanish by definition. In the next lemma, we prove that 
$(\nabla^2v)^{0,1}_1$ also vanishes.

\lem\label{1stSquare}
$$\bfig
\square/>`<-`<-`>/<1000,450>[C(H,V)`C(H\times G,W)`V`C(G,W),;\delta'`\delta`\delta`\delta']
\efig$$
commutes. \endlem

\pf
Let $v\in V$ and $(h;g)\in H\times G$, then
$$\mld (\delta\delta'v)(h;g) & =(\delta'v)(g^{h})-\rho_{H}^1(h;g)(\delta'v)(g) 
							   =\rho_1(g^{h})v-\rho_0^1(h)^{-1}\rho_1(g)v \\
							   =\rho_0^1(h)^{-1}\rho_1(g)(\rho_0^0(h)v-v) 
							   =\rho_0^1(h)^{-1}\rho_1(g)(\delta v)(h) = (\delta'\delta v)(h;g). $$
\epf

In fact, since for $p>0$ the action of $\G_p$ on $G$ and the right 
representation $\rho_{\G_p}^1$ are respectively pull-backs along $t_p$ 
of the action of $H$ on $G$ and the right representation $\rho_{H}^1$, 
the proof of Lemma \ref{1stSquare} implies the following corollary.

\corollary\label{1st pSquare} 
$$\bfig
\square/>`<-`<-`>/<1000,450>[C(\G_p,V)`C(\G_p\times G,W)`V`C(G,W),;\delta'`\delta`\delta`\delta']
\efig$$
commutes. \endcorollary

$(\nabla^2v)^{1,1}_0$, as it is, does not vanish. Let $\gamma\in\G$ and 
let $(g,h)\in G\rtimes H$ be its image under the isomorphism of Remark 
\ref{Equiv2Gp:Mods}, then
$$\mld (\nabla^2v)^{1,1}_0(\gamma) & =(\partial\delta v)(\gamma)=(\delta v)(h)-(\delta v)(hi(g)) \\
								     =\rho_0^0(h)(v-\rho_0^0(i(g))v)=-\rho_0^0(h)\big{(}\phi\circ\rho_1(g)v\big{)}, $$
where the last equality is the first part of Eq. (\ref{2}).	Let 
$$\bfig
\morphism<900,0>[\Delta:C^{0,0}_1(\G,\phi)`C^{1,1}_{0}(\G,\phi);]
\efig$$				   
be defined by
\begin{equation}\label{1.Delta}
\Delta\omega(\gamma):=\rho_0^0(h)\circ\phi(\omega(g)),
\end{equation} 
for $\omega\in C(G,W)$ and $\gamma\in\G$. Here, $(g,h)\in G\rtimes H$ 
is the image of $\gamma$ under the isomorphism of Remark \ref{Equiv2Gp:Mods}. 

\rem\label{upToHomSuite}
Observe that $\Delta$ is related to the structural map $\varrho$ of the representation up to homotopy induced by the $2$-representation (see Remark \ref{repUpToH}): $\omega\in C(G,W)$ defines the bundle map $\bar{\omega}:H\times G\to H\times W$ , $\bar{\omega}(h;g)=(h;\omega(g))$. Then, correctly interpreted, $\Delta\omega=\varrho\circ\bar{\omega}$. \endrem

Adding $\Delta$ to Eq. (\ref{preDiff}), makes $(\nabla^2v)^{1,1}_0=0$ as
$$(\Delta\delta'v)(\gamma)=\rho_0^0(h)\circ\phi((\delta'v)(g))=\rho_0^0(h)\big{(}\phi\circ\rho_1(g)v\big{)}, $$
for all $\gamma=(g,h)\in G\rtimes H\cong\G$.

Schematically, the updated differential of the $1$-cochain $\lambda$ is
$$\bfig
\Atriangle(0,500)/.`.`/[\delta\lambda_0`\partial\lambda_0-\delta v+\Delta\lambda_1`\delta'\lambda_0-\delta\lambda_1;``]
\Atriangle(-500,0)/.``./[`-v`-\delta'v;``]
\Atriangle(500,0)/`.`./[`-\delta'v`\delta_{(1)}\lambda_1,;``]
\morphism(500,700)/|->/<0,300>[\lambda_0`;]
\morphism(500,700)/|->/<-300,-150>[\lambda_0`;]
\morphism(500,700)/|->/<300,-150>[\lambda_0`;]
\morphism(0,200)/|->/<0,300>[v`;]
\morphism(0,200)/|->/<-300,-150>[v`;]
\morphism(0,200)/|->/<300,-150>[v`;]
\morphism(1000,200)/|->/<0,300>[\lambda_1`;]
\morphism(1000,200)/|->/<-300,-150>[\lambda_1`;]
\morphism(1000,200)/|->/<300,-150>[\lambda_1`;]
\morphism(1000,200)/=>/<-700,250>[\lambda_1`;]
\morphism(500,430)/./<0,250>[`;]
\morphism(500,430)/./<-400,-200>[`;]
\morphism(500,430)/./<400,-200>[`;]
\efig$$
where the solid arrows represent the grid differentials and the double 
arrow represents $\Delta$. 

If $(v,\lambda_0,\lambda_1)\in C^1_{tot}(\G,\phi)$ is a $1$-cocylce, then 
$v=0$, $\lambda_0$ is a crossed homomorphism of $H$ into $V$ with respect 
to $\rho_0^0$, and $\lambda_1$ is a crossed homomorphisms of $G$ into $W$ 
with respect to $\rho_0^1\circ i$. In symbols,
\begin{equation}\label{lambda0CrHom}
\lambda_0(h_0h_1) = \lambda_0(h_0)+\rho_0^0(h_0)\lambda_0(h_1),\qquad\qquad \forall h_0,h_1\in H,
\end{equation}
\begin{equation}\label{lambda1CrHom}
\lambda_1(g_0g_1) = \lambda_1(g_0)+\rho_0^1(i(g_0))\lambda_1(g_1),\qquad\quad \forall g_0,g_1\in G.
\end{equation}
Additionally, the following relations hold for every $\gamma\in\G$ and 
all $(h;g)\in H\times G$:
$$(\partial\lambda_0+\Delta\lambda_1)(\gamma) = 0, \qquad\qquad (\delta'\lambda_0-\delta\lambda_1)(h;g) = 0.$$
If $(g,h)\in G\rtimes H$ is the image of $\gamma$ under the isomorphism 
of Remark \ref{Equiv2Gp:Mods}, these relations are respectively
\begin{equation}\label{respSandT}
\lambda_0(h)+\rho_0^0(h)\circ\phi(\lambda_1(g)) = \lambda_0(hi(g)),
\end{equation} 
and
\begin{equation}\label{CrHomo}
\rho_0^1(h)^{-1}\rho_1(g)\lambda_0(h)+\rho_0^1(h)^{-1}\lambda_1(g) = \lambda_1(g^h).
\end{equation}
Eq. (\ref{respSandT}) implies that the map
$$\bar{\lambda}:\G\to<350> W\oplus V:(g,h)\to/->/<350>(\rho_0^1(h)\lambda_1(g),\lambda_0(h)) $$
respects both the source and the target; indeed, when $h=1$, it implies
the commutativity of 
$$\bfig
\square[G`W`H`V.;\lambda_1`i`\phi`\lambda_0]
\efig$$
In fact, $\bar{\lambda}$ is a functor. Let $(\gamma_0,\gamma_1)\in\G_2$ 
and let $(g_0,hi(g_1)),(g_1,h)\in G\rtimes H$ be their respective images 
under the isomorphism of Remark \ref{Equiv2Gp:Mods}. 
$\bar{\lambda}(\gamma_0)$ and $\bar{\lambda}(\gamma_1)$ are composable 
in the $2$-vector space; indeed, combining Eq.'s (\ref{respSandT}) and 
(\ref{0}), one gets
$$\lambda_0(hi(g))=\lambda_0(h)+\rho_0^0(h)\phi(\lambda_1(g))=\lambda_0(h)+\phi\big{(}\rho_0^1(h)\lambda_1(g)\big{)}. $$
Using Eq. (\ref{lambda1CrHom}), one computes the composition of 
$\bar{\lambda}(\gamma_0)$ and $\bar{\lambda}(\gamma_1)$ to be
$$\mld \bar{\lambda}(\gamma_0)\bar{\lambda}(\gamma_1) & =\big{(}\rho_0^1(h)\big{(}\lambda_1(g_1)+\rho_0^1(i(g_1))\lambda_1(g_0)\big{)},\lambda_0(h)\big{)} \\ =(\rho_0^1(h)\lambda_1(g_1g_0),\lambda_0(h))=\bar{\lambda}(g_1g_0,h)=\bar{\lambda}(\gamma_0\Join\gamma_1), $$
yielding the claim. Further using Eq. (\ref{CrHomo}), one shows that 
$\bar{\lambda}$ is a crossed homomorphism into $W\oplus V$ with respect 
to $\bar{\rho}$:
$$\mld \bar{\lambda}(\gamma_0\vJoin\gamma_1) & =\bar{\lambda}(g_0^{h_1}g_1,h_0h_1) 
        =(\rho_0^1(h_0h_1)\lambda_1(g_0^{h_1}g_1),\lambda_0(h_0h_1)) \\  
        =\big{(}\rho_0^1(h_0h_1)\big{(}\lambda_1(g_0^{h_1})+\rho_0^1(i(g_0^{h_1}))\lambda_1(g_1)\big{)},\lambda_0(h_0)+\rho_0^0(h_0)\lambda_0(h_1)\big{)} \\ 
        =\big{(}\rho_0^1(h_0)\big{(}\rho_1(g_0)\lambda_0(h_1)+\lambda_1(g_0) +\rho_0^1(i(g_0)h_1)\lambda_1(g_1)\big{)},\lambda_0(h_0)+\rho_0^0(h_0)\lambda_0(h_1)\big{)} \\
        =\bar{\lambda}(\gamma_0)+\bar{\rho}_{(g_0,h_0)}\bar{\lambda}(\gamma_1) $$

A coboundary $\nabla v$, will induce a crossed homomorphism-functor that 
can also be seen as $\bar{\rho}_{(g,h)}(0,v)$. We could not find a 
terminology for these in the literature; hence, we introduce the 
following definitions by analogy. We use the notation 
conventions of this section. 

\definition\label{crHomo}
The \textit{space of crossed functors} of a Lie $2$-group $\G$ 
with respect to a $2$-representation $\rho$ on the $2$-vector 
$\mathbbm{V}=W\to^\phi V$ is defined to be
$$\mld Cr\Hom(\G,\phi):=\lbrace \bar{\lambda}\in \Hom_{Gpd}(\G,\mathbbm{V}) & :\bar{\lambda}(g,h)=(\rho_0^1(h)\lambda_1(g),\lambda_0(h)) \\
  \qquad\textnormal{ is a crossed homomorphism with respect to }\bar{\rho}\rbrace . $$
The \textit{space of principal crossed functors} is defined to be
$$PCr\Hom(\G,\phi):=\lbrace \bar{\lambda}\in CrHom(\G,\phi):\bar{\lambda}(g,h)=\bar{\rho}_{(g,h)}(0,v)\textnormal{ for some }v\in V\rbrace . $$ \enddefinition

With these definitions, 
$$H^1_\nabla(\G,\phi)=Cr\Hom(\G,\phi)/PCr\Hom(\G,\phi). $$

\subsection{Degree 2}\label{H^2} A $2$-cochain $\vec{\omega}$ is a 
$6$-tuple $(v,\varphi,\omega_0,\lambda,\alpha,\omega_1)$, where
$$\omega_0\in C(H^2,V) $$
$$\varphi\in C(\G,V) \qquad\qquad \alpha\in C(H\times G,W) $$
$$v\in V \qquad\qquad \lambda\in C(G,W) \qquad\qquad \omega_1\in C(G^2,W). $$
The coordinates of the differential $\nabla\vec{\omega}$ are 
$$(\nabla\vec{\omega})^{3,0}_0=\partial v=0 \qquad\qquad (\nabla\vec{\omega})^{0,3}_0=\delta\omega_0 \qquad\qquad (\nabla\vec{\omega})^{0,0}_3=\delta_{(1)}\omega_1 $$
$$(\nabla\vec{\omega})^{1,2}_0=\partial\omega_0-\delta\varphi \qquad\qquad (\nabla\vec{\omega})^{0,2}_1=\delta'\omega_0-\delta\alpha $$
$$(\nabla\vec{\omega})^{2,1}_0=\delta v-\partial\varphi \quad (\nabla\vec{\omega})^{1,1}_1=\partial\alpha+\delta\lambda-\delta'\varphi \quad (\nabla\vec{\omega})^{0,1}_2=\delta_{(1)}\alpha+\delta\omega_1 $$
$$(\nabla\vec{\omega})^{2,0}_1=\delta'v-\cancelto{\lambda}{\partial\lambda} \qquad\qquad (\nabla\vec{\omega})^{1,0}_2=\cancelto{0}{\partial\omega_1}-\delta_{(1)}\lambda . $$

Let $\lambda=(v,\lambda_0,\lambda_1)\in C_{tot}^1(\G,\phi)$. We study 
$\nabla^2$ applied to each coordinate of $\lambda$ separately:
 
In $\nabla^2 v$,  
$(\nabla^2 v)^{3,0}_0=(\nabla^2 v)^{1,2}_0=0$,
because $\partial$ and $\delta$ are differentials. Moreover, 
$(\nabla^2 v)^{1,0}_2=0$ due to Lemma \ref{q=0r-cx}, and 
$(\nabla^2 v)^{1,1}_1=0$ due to Corollary \ref{1st pSquare}. By 
definition, $\partial:C^{1,0}_r(\G,\phi)\to C^{2,0}_r(\G,\phi)$ 
is the identity; hence, $(\nabla^2 v)^{2,0}_1=0$ trivially.  
$(\nabla^2 v)^{2,1}_0$, as it is, does not vanish. Let 
$(g_1,g_2,h)\in G^2\times H\cong\G_2$, then
$$\mld (\nabla^2 v)^{2,1}_0(g_1,g_2,h) & =(\partial\delta v-\delta\partial v)(g_1,g_2,h) 
										  = \rho_0^0(hi(g_2))v-v-\big{(}\rho_0^0(hi(g_2g_1))v-v\big{)} \\
										  = \rho_0^0(hi(g_2))\big{(}v-\rho_0^0(i(g_1))v\big{)}
										  = -\rho_0^0(hi(g_2))\big{(}\phi\circ\rho_1(g_1)v\big{)}, $$
where the last equality is the first part of Eq. (\ref{2}). Let 
$$\bfig
\morphism<900,0>[\Delta:C^{1,0}_1(\G,\phi)`C^{2,1}_{0}(\G,\phi);]
\efig$$				   
be defined by 
\begin{equation}\label{2.Delta}
\Delta\omega(g_1,g_2,h):=\rho_0^0(hi(g_2))\circ\phi(\omega(g_1)),
\end{equation}
for $\omega\in C(G,W)$ and $(g_1,g_2,h)\in G^2\times H\cong\G_2$ 
under the isomorphism of Remark \ref{2GpNerve}. Subtracting $\Delta$ 
from Eq. (\ref{preDiff}), makes $(\nabla^2 v)^{2,1}_0=0$ as
$$-\Delta(-\delta'v)(g_1,g_2,h)=\rho_0^0(hi(g_2))\circ\phi((\delta'v)(g))=\rho_0^0(hi(g_2))\big{(}\phi\circ\rho_1(g_1)v\big{)}, $$
for all $(g_1,g_2,h)\in\G_2$. All other components vanish from the onset.

In $\nabla^2\lambda_0$,  
$(\nabla^2\lambda_0)^{0,3}_0=(\nabla^2\lambda_0)^{2,1}_0=0$,
because $\delta$ and $\partial$ are differentials. Moreover, 
$(\nabla^2\lambda_0)^{0,1}_2=0$ due to Lemma \ref{q=0r-cx}. Let 
$(h_1,h_2;g)\in H^2\times G$, then
$$\mld (\delta\delta'\lambda_0)(h_1,h_2;g) & =(\delta'\lambda_0)(h_2;g^{h_1})-(\delta'\lambda_0)(h_1h_2;g)+\rho_0^1(h_2)^{-1}(\delta'\lambda_0)(h_1;g) \\
											 = \rho_0^1(h_2)^{-1}\rho_1(g^{h_1})\lambda_0(h_2)+\rho_0^1(h_1h_2)^{-1}\rho_1(g)\big{(}-\lambda_0(h_1h_2)+\lambda_0(h_1)\big{)} \\
											 =\rho_0^1(h_1h_2)^{-1}\rho_1(g)\big{(}\rho_0^0(h_1)\lambda_0(h_2)-\lambda_0(h_1h_2)+\lambda_0(h_1)\big{)}=(\delta'\delta\lambda_0)(h_1,h_2;g)$$ 
and $(\nabla^2\lambda_0)^{0,2}_1=0$. This computation can be easily 
generalized to prove the following:

\lem\label{1stpCol}
$$\bfig
\square/>`<-`<-`>/<1000,450>[C(\G_p^{q+1},V)`C(\G_p^{q+1}\times G,W)`C(\G_p^q,V)`C(\G_p^q\times G,W),;\delta'`\delta`\delta`\delta']
\efig$$
commutes for all $q>0$.
\endlem

\pf 
Let $\omega\in C(\G_p^q,V)$, $\vec{\gamma}=(\gamma_0,...,\gamma_q)^T\in\G_p^{q+1}$ 
and $g\in G$, then
$$\mld \delta'\delta\omega(\vec{\gamma};g) & =\rho_0^1(t_p(\gamma_0)...t_p(\gamma_q))^{-1}\rho_1(g)\delta\omega(\vec{\gamma}) \\
	 =\rho_0^1(t_p(\gamma_0)...t_p(\gamma_q))^{-1}\rho_1(g)\Big{(}\rho_0^0(t_p(\gamma_0))\omega(\delta_{0}\vec{\gamma})+\sum_{j=1}^{q+1}(-1)^{j}\omega(\delta_j\vec{\gamma})\Big{)}, $$
and
$$\mld \delta\delta'\omega(\vec{\gamma};g) & =\delta'\omega(\delta_0\vec{\gamma};g^{t_p(\gamma_0)})+\sum_{j=1}^{q}(-1)^{j}\delta'\omega(\delta_j\vec{\gamma};g)+(-1)^{q+1}\rho_{\G_p}^1(\gamma_q;g)\delta'\omega(\delta_{q+1}\vec{\gamma};g) \\
				 =\rho_0^1(t_p(\gamma_1)...t_p(\gamma_q))^{-1}\rho_1(g^{t_p(\gamma_0)})\omega(\delta_0\vec{\gamma})+ \\
				 \qquad +\sum_{j=1}^{q}(-1)^{j}\rho_0^1(t_p(\gamma_0)...t_p(\gamma_{j-1}\vJoin\gamma_j)...t_p(\gamma_q))^{-1}\rho_1(g)\omega(\delta_j\vec{\gamma})+ \\
				 \qquad\qquad +(-1)^{q+1}\rho_0^1(t_p(\gamma_q))^{-1}\rho_0^1(t_p(\gamma_0)...t_p(\gamma_{q-1}))^{-1}\rho_1(g)\omega(\delta_{q+1}\vec{\gamma}) \\
				 =\rho_0^1(t_p(\gamma_0)...t_p(\gamma_q))^{-1}\rho_1(g)\Big{(}\rho_0^0(t_p(\gamma_0))\omega(\delta_{0}\vec{\gamma})+\sum_{j=1}^{q+1}(-1)^{j}\omega(\delta_j\vec{\gamma})\Big{)}. $$ \epf

Let $(\gamma;f)\in\G\times G$ 
and let $(g,h)\in G\rtimes H$ be the image of $\gamma$ under the 
isomorphism of Remark \ref{Equiv2Gp:Mods}, then
$$\mld (\partial\delta'\lambda_0)(\gamma;f) & =\rho_0^1(i(g))^{-1}(\delta'\lambda_0)(h;f)-(\delta'\lambda_0)(hi(g);f) \\
											 = \rho_0^1(i(g))^{-1}\rho_0^1(h)^{-1}\rho_1(f)\lambda_0(h)-\rho_0^1(hi(g))^{-1}\rho_1(f)\lambda_0(hi(g)) \\
											 =\rho_0^1(hi(g))^{-1}\rho_1(f)\big{(}\lambda_0(h)-\lambda_0(hi(g))\big{)}=(\delta'\partial\lambda_0)(\gamma;f)$$
and $(\nabla^2\lambda_0)^{1,1}_1=0$. This computation can be easily 
generalized to prove the following: 

\lem\label{1stqSquare}
$$\bfig
\square/>`<-`<-`>/<1000,450>[C(\G_p^q\times G,W)`C(\G_{p+1}^q\times G,W)`C(\G_p^q,V)`C(\G_{p+1}^q,V),;\partial`\delta'`\delta'`\partial]
\efig$$
commutes for all $q>1$. \endlem

\pf
We adopt the convention that, for $\gamma_a=\pmatrix{
\gamma_{a0} & \cdots & \gamma_{ap}}\in\G_{p+1}$ and $\pmatrix{
g_{ab} \cr 
h_{ab} }$
is the image of $\gamma_{ab}$ under the isomorphism of Remark 
\ref{Equiv2Gp:Mods}. 

Let $\omega\in C(\G_p^q,V)$, $\vec{\gamma}=(\gamma_1,...,\gamma_q)^T\in\G_{p+1}^q$ 
and $f\in G$, then 
$$\mld \partial\delta'\omega(\vec{\gamma};f) & =\rho_0^1(i(pr_G(\gamma_{10}\vJoin ...\vJoin\gamma_{q0})))^{-1}\delta'\omega(\partial_0\vec{\gamma};f)+\sum_{j=1}^{p+1}(-1)^j\delta'\omega(\partial_j\vec{\gamma};f) \\
          =\rho_0^1(i(pr_G(\gamma_{10}\vJoin ...\vJoin\gamma_{q0})))^{-1}\rho_0^1(t_{p}(\partial_0\gamma_1)...t_{p}(\partial_0\gamma_q))^{-1}\rho_1(f)\omega(\partial_0\vec{\gamma})+ \\
          \qquad +\sum_{j=1}^{p+1}(-1)^j\rho_0^1(t_{p}(\partial_j\gamma_1)...t_{p}(\partial_j\gamma_q))^{-1}\rho_1(f)\omega(\partial_j\vec{\gamma}). $$
          
If $j>0$, $t_p(\partial_j\gamma_a)=t(\gamma_{a0})$; otherwise, 
$t_p(\partial_0\gamma_a)=s(\gamma_{a0})$. Hence, since $s$ and $t$ are 
group homomorphisms,           
$$\mld \partial\delta'\omega(\vec{\gamma};f) & =\rho_0^1(t(\gamma_{10})...t_p(\gamma_{q0}))^{-1}\rho_1(f)\sum_{j=0}^{p+1}(-1)^j\omega(\partial_j\vec{\gamma}) \\
           =\rho_0^1(t_{p+1}(\gamma_1)...t_{p+1}(\gamma_q))^{-1}\rho_1(f)\partial\omega(\vec{\gamma})=\delta'\partial\omega(\vec{\gamma};f). $$ \epf
           
$(\nabla^2\lambda_0)^{1,2}_0$, as it is, does not vanish. Let 
$(\gamma_1,\gamma_2)\in\G^2$ and let $(g_k,h_k)\in G\rtimes H$ be the 
image of $\gamma_k$ under the isomorphism of Remark \ref{Equiv2Gp:Mods} 
for $k\in\lbrace 1,2\rbrace$, then
$$\mld (\nabla^2\lambda_0)^{1,2}_0(\gamma_1,\gamma_2) & =(\partial\delta\lambda_0-\delta\partial\lambda_0)(\gamma_1,\gamma_2) \\
										  = \rho_0^0(h_1)\big{(}\lambda_0(h_2)-\rho_0^0(i(g_1))\lambda_0(h_2)\big{)} 
										  = -\rho_0^0(h_1)\big{(}\phi\circ\rho_1(g_1)\lambda_0(h_2)\big{)}, $$
where the last equality is the first part of Eq. (\ref{2}). Let 
$$\bfig
\morphism<900,0>[\Delta:C^{0,1}_1(\G,\phi)`C^{1,2}_{0}(\G,\phi);]
\efig$$				   
be defined by 
\begin{equation}\label{3.Delta}
\Delta\omega(\gamma_1,\gamma_2):=\rho_0^0(h_1h_2)\circ\phi(\omega(h_2;g_1)),
\end{equation} 
for $\omega\in C(H\times G,W)$ and $(\gamma_1,\gamma_2)\in\G^2$. Here, 
$(g_k,h_k)\in G\rtimes H$ is the image of $\gamma_k$ under the 
isomorphism of Remark \ref{Equiv2Gp:Mods} for $k\in\lbrace 1,2\rbrace$. 
Adding $\Delta$ to Eq. (\ref{preDiff}), makes $(\nabla^2\lambda_0)^{1,2}_0=0$ 
as a consequence of Eq. (\ref{0}) and that for all $(\gamma_1,\gamma_2)\in\G^2$
$$(\Delta\delta'\lambda_0)(\gamma_1,\gamma_2)=\rho_0^0(h_1h_2)\circ\phi((\delta'\lambda_0)(h_2;g_1))=\rho_0^0(h_1h_2)\circ\phi\big{(}\rho_0^1(h_2)^{-1}\rho_1(g_1)\lambda_0(h_2)\big{)}. $$
All other components vanish from the onset.

In $\nabla^2\lambda_1$, 
$(\nabla^2\lambda_1)^{0,0}_3=(\nabla^2\lambda_1)^{2,0}_1=(\nabla^2\lambda_1)^{0,2}_1=0$,
because $\delta_{(1)}$, $\partial$ and $\delta$ are differentials. 
Since by definition, $\partial:C^{0,0}_r(\G,\phi)\to C^{1,0}_r(\G,\phi)$ 
is zero, $(\nabla^2 \lambda_1)^{1,0}_2=0$ trivially. Moreover, 
$\delta_{(1)}$ commutes with $\delta$; hence, 
$(\nabla^2\lambda_1)^{0,1}_2=0$. Indeed, let $(h;g_1,g_2)\in H\times G^2$, 
then
$$\mld (\delta\delta_{(1)}\lambda_1)(h;g_1,g_2) & =(\delta_{(1)}\lambda_1)(g_1^h,g_2^h)-\rho_0^1(h)^{-1}(\delta_{(1)}\lambda_1)(g_1,g_2) \\
											 = \rho_0^1(i(g_1^h))\lambda_1(g_2^h)-\lambda_1((g_1g_2)^h)+\lambda_1(g_1^h)+ \\
						\qquad-\rho_0^1(h)^{-1}\big{(}\rho_0^1(i(g_1))\lambda_1(g_2)-\lambda_1(g_1g_2)+\lambda_1(g_1)\big{)} \\
											 = \rho_0^1(i(g_1^h))\big{(}\lambda_1(g_2^h)-\rho_0^1(h)^{-1}\lambda_1(g_2)\big{)}-(\delta\lambda_1)(h;g_1g_2)+(\delta\lambda_1)(h;g_1) \\
											 =(\delta_{(1)}\delta\lambda_1)(h;g_1,g_2).$$ 
$(\nabla^2\lambda_1)^{1,1}_1$, $(\nabla^2\lambda_1)^{2,1}_0$ and $(\nabla^2\lambda_1)^{1,2}_0$, 
as they are, do not vanish. Let $(\gamma;f)\in\G\times G$ and let 
$(g,h)\in G\rtimes H$ be the image of $\gamma$ under the 
isomorphism of Remark \ref{Equiv2Gp:Mods}, then
$$\mld (\nabla^2\lambda_1)^{1,1}_1(\gamma;f) & =(\delta\cancelto{0}{\partial\lambda_1}-\partial\delta\lambda_1-\delta'\Delta\lambda_1)(\gamma;f) \\
										  = \lambda_1(f^{hi(g)})-\rho_0^1(i(g))^{-1}\lambda_1(f^h)-\rho_0^1(hi(g))^{-1}\rho_1(f)\rho_0^0(h)\phi\lambda_1(g) \\
										  = \lambda_1(f^{hi(g)})-\rho_0^1(i(g))^{-1}\big{(}\lambda_1(f^h)+(\rho_0^1(i(f^h))-I)\lambda_1(g)\big{)} $$
where the last equality follows from Eq. (\ref{3}) combined with 
the second part of Eq. (\ref{2}). Let 
$$\bfig
\morphism<900,0>[\Delta:C^{0,0}_2(\G,\phi)`C^{1,1}_{1}(\G,\phi);]
\efig$$				   
be defined by 
\begin{equation}\label{DeltaOffFront}
\Delta\omega(\gamma;f):=\rho_0^1(i(g)))^{-1}\omega(f^h,g)+\omega(g^{-1},f^hg)-\omega(g^{-1},g),
\end{equation}
for $\omega\in C(G^2,W)$ and $(\gamma;f)\in\G\times G$. Here, 
$(g,h)\in G\rtimes H$ is the image of $\gamma$ under the 
isomorphism of Remark \ref{Equiv2Gp:Mods}. 
Adding $\Delta$ to Eq. (\ref{preDiff}), makes $(\nabla^2\lambda_1)^{1,1}_1=0$ 
as, using the cocycle equation $\delta_{(1)}^2\lambda_1=0$, 
$$\mld (\Delta\delta_{(1)}\lambda_1)(\gamma;f) & =\rho_0^1(i(g)))^{-1}(\delta_{(1)}\lambda_1)(f^h,g)+(\delta_{(1)}\lambda_1)(g^{-1},f^hg)-(\delta_{(1)}\lambda_1)(g^{-1},g) \\ 
												 =(\delta_{(1)}\lambda_1)(g^{-1},f^h)+(\delta_{(1)}\lambda_1)(g^{-1}f^h,g)-(\delta_{(1)}\lambda_1)(g^{-1},g) \\ 
												 = \rho_0^1(i(g))^{-1}\lambda_1(f^h)+\rho_0^1(i(g^{-1}f^h))\lambda_1(g)-\lambda_1(f^{hi(g)})-\rho_0^1(i(g))^{-1}\lambda_1(g) $$
for all $(\gamma;f)\in\G\times G$. 

Let $(g_1,g_2,h)\in G^2\times H\cong\G_2$, then
$$\mld (\nabla^2\lambda_1)^{2,1}_0(g_1,g_2,h) & =-(\partial\Delta\lambda_1+\Delta\cancelto{0}{\partial\lambda_1})(g_1,g_2,h) \\
										  = -\Delta\lambda_1(g_2,h)+\Delta\lambda_1(g_2g_1,h)-\Delta\lambda_1(g_1,hi(g_2)) \\
										  = -\rho_0^0(h)\circ\phi\big{(}\lambda_1(g_2)-\lambda_1(g_2g_1)+\rho_0^1(i(g_2))\lambda_1(g_1)\big{)}. $$
Let 
$$\bfig
\morphism<900,0>[\Delta_{2,1}:C^{0,0}_2(\G,\phi)`C^{2,1}_{0}(\G,\phi);]
\efig$$				   
be defined by 
\begin{equation}\label{Delta21}
\Delta_{2,1}\omega(g_1,g_2,h):=\rho_0^0(h)\circ\phi(\omega(g_2,g_1)),
\end{equation}
for $\omega\in C(G^2,W)$ and $(g_1,g_2,h)\in G^2\times H\cong\G_2$. 
Adding $\Delta_{2,1}$ to Eq. (\ref{preDiff}), makes 
$(\nabla^2\lambda_1)^{2,1}_0=0$ as 
$$(\Delta_{2,1}\delta_{(1)}\lambda_1)(g_1,g_2,h)=\rho_0^0(h)\circ\phi\big{(}\delta_{(1)}\lambda_1(g_2,g_1)\big{)}  $$
for all $(g_1,g_2,h)\in G^2\times H\cong\G_2$. 

Let $(\gamma_1,\gamma_2)\in\G^2$ and let $(g_k,h_k)\in G\rtimes H$ be 
the image of $\gamma_k$ under the isomorphism of Remark 
\ref{Equiv2Gp:Mods} for $k\in\lbrace 1,2\rbrace$, then
$$\mld (\nabla^2\lambda_1)^{1,2}_0(\gamma_1,\gamma_2) & =-(\delta\Delta\lambda_1+\Delta\delta\lambda_1)(\gamma_1,\gamma_2) \\
										  = -\rho_0^0(h_1h_2)\circ\phi\big{(}\rho_0^1(i(g_1^{h_2}))\lambda_1(g_2)-\lambda_1(g_1^{h_2}g_2)+\lambda_1(g_1^{h_2})\big{)} $$									
Let 
$$\bfig
\morphism<900,0>[\Delta_{1,2}:C^{0,0}_2(\G,\phi)`C^{1,2}_{0}(\G,\phi);]
\efig$$				   
be defined by 
\begin{equation}\label{Delta12}
\Delta_{1,2}\omega(\gamma_1,\gamma_2):=\rho_0^0(h_1h_2)\circ\phi(\omega(g_1^{h_2},g_2)),
\end{equation}
for $\omega\in C(G^2,W)$ and $(\gamma_1,\gamma_2)\in\G^2$. Here, 
$(g_k,h_k)\in G\rtimes H$ is the image of $\gamma_k$ under the 
isomorphism of Remark \ref{Equiv2Gp:Mods} for $k\in\lbrace 1,2\rbrace$. 
Adding $\Delta_{1,2}$ to Eq. (\ref{preDiff}), makes 
$(\nabla^2\lambda_1)^{1,2}_0=0$ as 
$$(\Delta_{1,2}\delta_{(1)}\lambda_1)(\gamma_1,\gamma_2)=\rho_0^0(h_1h_2)\circ\phi\big{(}\delta_{(1)}\lambda_1(g_1^{h_2},g_1)\big{)} $$
for all $(\gamma_1,\gamma_2)\in\G^2$. The remaining two components 
vanish from the onset.

Ultimately, the updated differential $\nabla$ is the one appearing in Eq. 
(\ref{preDiff3}). $\nabla$ squares to zero in the truncated complex 
(\ref{trunc3}) whose building blocks we represent schematically:

$$\bfig
\node 000(3050,600)[C^{0,0}_0]
\node 001(2450,400)[C^{0,0}_1]
\node 010(3050,1200)[C^{0,1}_0]
\node 100(2050,600)[C^{1,0}_0]
\node 002(1850,200)[C^{0,0}_2]
\node 020(3050,1800)[C^{0,2}_0]
\node 200(1050,600)[C^{2,0}_0]
\node 110(2050,1200)[C^{1,1}_0]
\node 101(1450,400)[C^{1,0}_1]
\node 011(2450,1000)[C^{0,1}_1]
\node 003(1050,0)[C^{0,0}_3]
\node 030(3050,2400)[C^{0,3}_0]
\node 300(0,600)[C^{3,0}_0]
\node 120(2050,1800)[C^{1,2}_0]
\node 210(1050,1200)[C^{2,1}_0]
\node 102(650,200)[C^{1,0}_2]
\node 201(450,400)[C^{2,0}_1]
\node 012(1850,800)[C^{0,1}_2]
\node 021(2450,1600)[C^{0,2}_1]
\node 111(1450,1000)[C^{1,1}_1]
\node 4(1850,1400)[\circ]        
\arrow|r|[000`001;\delta_{(1)}]
\arrow[000`010;\delta]
\arrow||/@{->}_<>(0.3){\partial}/[000`100;]
\arrow[001`002;]
\arrow[001`011;]
\arrow[001`101;]
\arrow/=>/[001`110;]
\arrow[010`011;]
\arrow[010`020;]
\arrow[010`110;]
\arrow[100`101;]
\arrow[100`110;]
\arrow[100`200;]
\arrow[002`003;]
\arrow[002`012;]
\arrow[002`102;]
\arrow/=>/[002`111;]
\arrow/-->/[002`120;]
\arrow/-->/[002`210;]
\arrow[020`021;]
\arrow[020`030;]
\arrow[020`120;]
\arrow[200`201;]
\arrow[200`210;]
\arrow[200`300;]
\arrow[110`111;]
\arrow[110`120;]
\arrow[110`210;]
\arrow[101`102;]
\arrow[101`111;]
\arrow[101`201;]
\arrow/=>/[101`210;]
\arrow[011`012;]
\arrow[011`021;]
\arrow[011`111;]
\arrow/=>/[011`120;]
\arrow/./[012`4;]         
\arrow/./[021`4;]         
\arrow/./[012`111;]       
\arrow/./[111`210;]       
\arrow/./[101`200;]       
\arrow/./[120`4;]         
\arrow/./[002`100;]       
\arrow/./[002`101;]       
\arrow/./[120`021;]       
\arrow/./[100`001;]       
\efig$$
\begin{equation}\label{trunc}
\textnormal{Fig.}
\end{equation}
In this diagram, the double arrows and the dashed arrows represent 
respectively the first difference maps 
$$\bfig
\morphism<900,0>[\Delta:C^{p,q}_r(\G,\phi)`C^{p+1,q+1}_{r-1}(\G,\phi);]
\efig$$
and, for $(a,b)\in\lbrace(1,2),(2,1)\rbrace$, the second difference maps
$$\bfig
\morphism<1000,0>[\Delta_{a,b}:C^{p,q}_r(\G,\phi)`C^{p+a,q+b}_{r-2}(\G,\phi),;]
\efig$$
which owe the ordinal in their names to the degree difference in the 
$r$-direction.

\subsubsection{The cocycle equations}\label{CocEqns}
The second cohomology group $H^2_\nabla(\G,\phi)$ is in one-to-one 
correspondence with the equivalence classes of split extensions of 
the Lie $2$-group $\G$ by the $2$-vector space $W\to^\phi V$. To 
make this patent, we express an extension in terms of simpler data 
-as in the proof of Proposition \ref{Ind2GpRep}- and determine the 
equations these need to satisfy to build the extension back up.

\proposition\label{Gp2-cocycles}
Let $\rho$ be a $2$-representation of $G\to^iH$ on $W\to^\phi V$ and
$(\check{\varphi},\omega_0,\alpha,\omega_1)\in C(G,V)\oplus C(H^2,V)\oplus C(H\times G,W)\oplus C(G^2,W) $. 
Let $\mathcal{E}_{(\check{\varphi},\omega_0,\alpha,\omega_1)}$ be the 
space in the middle of the sequence 
$$\bfig
\square/>``>`>/<525,500>[1`W`1`V;```]
\square(525,0)/^{ (}->``>`^{ (}->/<750,500>[W`G{}_{\rho_0^1\circ i}\ltimes^{\omega_1}W`V`H{}_{\rho_0^0}\ltimes^{\omega_0}V;j_1`\phi``j_0]
\square(1275,0)/>``>`>/<750,500>[G{}_{\rho_0^1\circ i}\ltimes^{\omega_1}W`G`H{}_{\rho_0^0}\ltimes^{\omega_0}V`H;pr_1`\epsilon`i`pr_1]
\square(2025,0)/>```>/<525,500>[G`1`H`1,;```]
\efig$$
with 
\begin{equation}\label{structHom}
\epsilon(g,w) := (i(g),\phi(w)+\check{\varphi}(g)),
\end{equation}
and
\begin{equation}\label{action}
(g,w)^{(h,v)} := (g^h,\rho^1_0(h)^{-1}(w+\rho_1(g)v)+\alpha(h;g)),\qquad\textnormal{{\em for }}(h,v)\in H\times V. 
\end{equation} 
Then, 
\begin{itemize}
\item[i)] the product on $H{}_{\rho_0^0}\ltimes^{\omega_0}V$ (cf. Eq. 
(\ref{twstdProd})) is associative if and only if $\omega_0$ is a group 
cocycle with respect to $\rho_0^0$; 
\item[ii)] the product on $G{}_{\rho_0^1\circ i}\ltimes^{\omega_1}W$ is 
associative if and only if $\omega_1$ is a group cocycle with respect to 
$\rho_0^1\circ i$; 
\item[iii)] Eq. (\ref{structHom}) defines a Lie group homomorphism if 
and only if for all $g_1,g_2\in G$,
\begin{equation}\label{iii)}
\phi(\omega_1(g_1,g_2))-\omega_0(i(g_1),i(g_2))=\rho^0_0(i(g_1))\check{\varphi}(g_2)-\check{\varphi}(g_1g_2)+\check{\varphi}(g_1);
\end{equation}
\item[iv)] Eq. (\ref{action}) defines a right action if for all 
$h_1,h_2\in H$ and all $g\in G$,
\begin{equation}\label{iv)}
\rho^1_0(h_1h_2)^{-1}\rho_1(g)\omega_0(h_1,h_2)=\rho^1_0(h_2)^{-1}\alpha(h_1;g)-\alpha(h_1h_2;g)+\alpha(h_2;g^{h_1});
\end{equation}
\item[v)] Eq. (\ref{structHom}) defines an equivariant map if and only 
if for all $(h;g)\in H\times G$
\begin{equation}\label{v)}
\check{\varphi}(g^h)-\rho^0_0(h^{-1})\check{\varphi}(g)+\phi(\alpha(h;g))=\rho^0_0(h^{-1})\omega_0(i(g),h)+\omega_0(h^{-1},i(g)h)-\omega_0(h^{-1},h),
\end{equation}
\item[vi)] Peiffer equation holds if and only if for all $g_1,g_2\in G$, 
\begin{equation}\label{vi)}
\rho^1_0(i(g_2))^{-1}\rho_1(g_1)\check{\varphi}(g_2)+\alpha(i(g_2);g_1)=\rho^1_0(i(g_2))^{-1}\omega_1(g_1,g_2)+\omega_1(g_2^{-1},g_1g_2)-\omega_1(g_2^{-1},g_2),
\end{equation}
and
\item[vii)] Eq. (\ref{action}) defines an action by automorphisms if for 
all $h\in H$ and all $g_1,g_2\in G$,
\begin{equation}\label{vii)}
\rho^1_0(h)^{-1}\omega_1(g_1,g_2)-\omega_1(g_1^h,g_2^h)=\rho^1_0(i(g_1^h))\alpha(h;g_2)-\alpha(h;g_1g_2)+\alpha(h;g_1).
\end{equation}
\end{itemize}
Therefore, $\mathcal{E}_{(\check{\varphi},\omega_0,\alpha,\omega_1)}$ is 
a Lie $2$-group extension if and only if $i)$-$vii)$ hold. \endproposition

\rem\label{sanityCheck}
As a consistency check, we prove that, when applied to the trivial 
$2$-representation on $(0)\to\Rr$, Proposition \ref{Gp2-cocycles} is 
equivalent to Lemma \ref{easyExt}. For this case, the only 
non-trivial equations in Proposition \ref{Gp2-cocycles} are 
\begin{itemize}
    \item[i)] $\delta\omega_0=0$,
    \item[ii)] $-\omega_0(i(g_1),i(g_2))=\check{\varphi}(g_2)-\check{\varphi}(g_1g_2)+\check{\varphi}(g_1)$, 
    and
    \item[iii)] $\check{\varphi}(g^h)-\check{\varphi}(g)=\omega_0(i(g),h)+\omega_0(h^{-1},i(g)h)-\omega_0(h^{-1},h)$.
\end{itemize}
Let $(F,f)\in\Omega^2(\G)$ be a $2$-cocycle, $\omega_0:=F$ and 
$\check{\varphi}(g):=f(g,1)$. Then Eq. i) holds right away and  
Eq. ii) coincides with Eq. (\ref{singledOut}). Eq. iii) follows 
from evaluating 
$\delta f-\partial F=0$ at
$$\pmatrix{g & 1 \cr
    h^{-1} & h} : \cancel{f\pmatrix{1\cr h}}-f\pmatrix{g^h\cr 1}+f\pmatrix{g\cr h^{-1}} = F(h^{-1},h)-F(h^{-1}i(g),h), $$
$$\pmatrix{1 & g \cr
    h^{-1} & 1} : f\pmatrix{g\cr 1}-f\pmatrix{g\cr h^{-1}}+\cancel{f\pmatrix{1\cr h^{-1}}} = \cancel{F(h^{-1},1)}-F(h^{-1},i(g)). $$
Adding these expressions together yields
$$f\pmatrix{g\cr 1}-f\pmatrix{g^h\cr 1}=F(h^{-1},h)-F(h^{-1}i(g),h)-F(h^{-1},i(g)) $$
and, since $\delta F=0$, 
$$F(h^{-1}i(g),h)+F(h^{-1},i(g))=F(i(g),h)+F(h^{-1},i(g)h), $$
and the claim follows. Conversely, suppose $\omega_0$ and 
$\check{\varphi}$ verify Eq.'s i)-iii), and let $F:=\omega_0$ and 
$f(g,h):=\omega_0(h,i(g))+\check{\varphi}(g)$. Naturally, $\delta F=0$ 
holds right away. Using successively Eq.'s ii) and i),  
$$\mld f\pmatrix{g_2g_1\cr h} & =\omega_0(h,i(g_2g_1))+\check{\varphi}(g_2g_1) 
           =\omega_0(h,i(g_2g_1))+\omega_0(i(g_2),i(g_1))+\check{\varphi}(g_1)+\check{\varphi}(g_2) \\
           =\omega_0(hi(g_2),i(g_1))+\omega_0(h,i(g_2))+\check{\varphi}(g_1)+\check{\varphi}(g_2)=f\pmatrix{g_1\cr hi(g_2)}+f\pmatrix{g_2\cr h}; $$
hence, $\partial f=0$. As a consequence of $\partial f=0$ and Eq. ii),
$$\mld f\pmatrix{g_1^{h_2}g_2\cr h_1h_2} & =\omega_0(h_1h_2i(g_1^{h_2}),i(g_2))+\check{\varphi}(g_2)+\omega_0(h_1h_2,i(g_1^{h_2}))+\check{\varphi}(g_1^{h_2}) \\
                             =\omega_0(h_1h_2i(g_1^{h_2}),i(g_2))+\check{\varphi}(g_2)+\omega_0(h_1h_2,i(g_1^{h_2}))+\check{\varphi}(g_1)+ \\
                             \qquad +\omega_0(i(g_1),h_2)+\omega_0(h_2^{-1},i(g_1)h_2)-\omega_0(h_2^{-1},h_2). $$
Subtracting this expression from
$$f\pmatrix{g_2\cr h_2}+f\pmatrix{g_1\cr h_1}=\omega_0(h_2,i(g_2))+\check{\varphi}(g_2)+\omega_0(h_1,i(g_1))+\check{\varphi}(g_1), $$
one gets
$$\mld \delta f\pmatrix{g_1 & g_2 \cr
                   h_1 & h_2} & =\omega_0(h_2,i(g_2))+\omega_0(h_1,i(g_1))-\omega_0(h_1i(g_1)h_2,i(g_2))-\omega_0(h_1h_2,i(g_1^{h_2}))+ \\
                    \qquad -\omega_0(i(g_1),h_2)-\omega_0(h_2^{-1},i(g_1)h_2)+\omega_0(h_2^{-1},h_2). $$
Since 
$\delta\omega_0(h_1i(g_1),h_2,i(g_2))-\delta\omega_0(h_1,i(g_1),h_2)=0$,  
$$\omega_0(h_2,i(g_2))+\omega_0(h_1,i(g_1))-\omega_0(h_1i(g_1)h_2,i(g_2))-\omega_0(i(g_1),h_2)=\omega_0(h_1,i(g_1)h_2)-\omega_0(h_1i(g_1),h_2i(g_2)) $$
and since 
$\delta\omega_0(h_1h_2,h_2^{-1},i(g_2)h_2)-\delta\omega_0(h_1,h_2,h_2^{-1})=0$,
$$\omega_0(h_1,i(g_1)h_2)-\omega_0(h_1h_2,i(g_1^{h_2}))-\omega_0(h_2^{-1},i(g_1)h_2)=\omega_0(h_1,h_2)-\omega_0(h_2,h_2^{-1}).$$
Therefore, 
$$\delta f\pmatrix{g_1 & g_2 \cr
                h_1 & h_2} =\omega_0(h_1,h_2)-\omega_0(h_2,h_2^{-1})-\omega_0(h_1i(g_1),h_2i(g_2))+\omega_0(h_2^{-1},h_2)=\partial\omega_0\pmatrix{g_1 & g_2 \cr
                 h_1 & h_2 }, $$
where the last equality follows from $\delta\omega_0(h_2^{-1},h_2,h_2^{-1})=\omega_0(h_2,h_2^{-1})-\omega_0(h_2^{-1},h_2)=0$.
\endrem

Next, we consider equivalent extensions $\mathcal{E}$ and $\mathcal{F}$ 
as in 
$$\bfig
\node 1a(1200,0)[F_0]
\node 2a(0,150)[1]
\node 2b(400,150)[V]
\node 2c(1600,150)[H]
\node 2d(2000,150)[1.]
\node 3a(800,300)[E_0]
\node 4a(1200,575)[F_1]
\node 5a(0,725)[1]
\node 5b(400,725)[W]
\node 5c(1600,725)[G]
\node 5d(2000,725)[1]
\node 6a(800,875)[E_1]
\arrow[2a`2b;]
\arrow[2b`3a;]
\arrow[2b`1a;]
\arrow[3a`2c;]
\arrow[1a`2c;]
\arrow[2c`2d;]
\arrow[5a`5b;]
\arrow[5b`6a;]
\arrow[5b`4a;]
\arrow[6a`5c;]
\arrow[4a`5c;]
\arrow[5c`5d;]
\arrow[5b`2b;]
\arrow|l|[6a`3a;\epsilon]
\arrow|r|[4a`1a;f]
\arrow[5c`2c;]
\arrow/.>/[3a`1a;\psi_0]
\arrow/.>/[6a`4a;\psi_1]
\efig$$

Supposing that the extensions split, one can apply the decomposition of 
the proof of Proposition \ref{Ind2GpRep}. Picking a splitting of either 
extension and composing it with the isomorphism yields a splitting for 
the other. In this manner, the induced $2$-representations of 
Proposition \ref{Ind2GpRep} are identical, and we identify both 
$\mathcal{E}$ and $\mathcal{F}$ with a corresponding (twisted) 
semi-direct product. We write the functor $\psi$ in these coordinates as
$$\psi_k(z,a)=(\psi_k^{Gp}(z,a),\psi_k^{Vec}(z,a)), $$
for $k\in\lbrace 0,1\rbrace$. Both $\psi_0$ and $\psi_1$ respect 
inclusions and projections; hence, $\psi_k(1,a)=(1,a)$ and 
$\psi_k^{Gp}(z,a)=z$. Further using that both $\psi_k$'s are group 
homomorphisms together with the factorization 
$(z,a)=(1,a)\vJoin(z,0)$,
$$\mld \psi_k(z,a) & =\psi_k((1,a)\vJoin(z,0)) 
			 =\psi_k(1,a)\vJoin\psi_k(z,0)\\
			 =(1,a)\vJoin(z,\psi_k^{Vec}(z,0)) 
			 =(z,a+\cancel{\rho_k(1)}(\psi_k^{Vec}(z,0))+\cancel{\omega_k(1,z)}), $$
where $\rho_k$ stands for $\rho_0^0$ when $k=0$ and for 
$\rho_0^1\circ i$ when $k=1$. As a consequence,
$$\psi_k(z,a)=(z,a+\lambda_k(z)), $$
where the maps $\lambda_0:H\to V$ and $\lambda_1:G\to W$ are defined by 
$\psi_k^{Vec}(z,0)$ for $k=0$ and $k=1$ respectively.
\proposition\label{Gp2-coboundaries}
Let $\rho$ be a $2$-representation of $G\to^iH$ on $W\to^\phi V$ and let 
$(\check{\varphi},\omega_0,\alpha,\omega_1)$, 
$(\check{\varphi}',\omega_0',\alpha',\omega_1')$ be a pair of $4$-tuples 
verifying the equations of Proposition \ref{Gp2-cocycles}. Then the 
extensions $\mathcal{E}_{(\check{\varphi},\omega_0,\alpha,\omega_1)}$ 
and $\mathcal{E}_{(\check{\varphi}',\omega_0',\alpha',\omega_1')}$ are 
isomorphic if and only if there are maps $\lambda_0\in C(H,V)$ and 
$\lambda_1\in C(G,W)$ verifying  
\begin{itemize}
\item[i)] $\omega_0-\omega_0'=\delta\lambda_0$. Explicitly, for all 
pairs $h_0,h_1\in H$, 
$$\omega_0(h_0,h_1)-\omega_0'(h_0,h_1)=\rho^0_0(h_0)\lambda_0(h_1)-\lambda_0(h_0h_1)+\lambda_0(h_0). $$ 
\item[ii)] $\omega_1-\omega_1'=\delta\lambda_1$. Explicitly, for all 
pairs $g_0,g_1\in G$, 
$$\omega_1(g_0,g_1)-\omega_1'(g_0,g_1)=\rho^1_0(i(g_0))\lambda_1(g_1)-\lambda_1(g_0g_1)+\lambda_1(g_0). $$
\item[iii)] For all $g\in G$,
\begin{equation}\label{cohPhi}
\check{\varphi}(g)-\check{\varphi}'(g)=\phi(\lambda_1(g))-\lambda_0(i(g)).
\end{equation}
\item[iv)] For all $h\in H$ and all $g\in G$,
\begin{equation}\label{cohAlpha}
\alpha(h;g)-\alpha'(h;g)=\rho^1_0(h)^{-1}(\lambda_1(g)+\rho_1(g)\lambda_0(h))-\lambda_1(g^h).
\end{equation}
\end{itemize}\endproposition

\pf
Items $i)$ and $ii)$ are the usual relations identifying isomorphic 
extensions with cohomologous cocycles of Lie groups. Define 
$\psi_k(z,a):=(z,a+\lambda_k(z))$ for $k\in\lbrace 0,1\rbrace$. Then, 
$$\bfig
\square<850,500>[G{}_{\rho^1_0\circ i}\ltimes^{\omega_1}W`G{}_{\rho^1_0\circ i}\ltimes^{\omega_1'}W`H{}_{\rho^0_0}\ltimes^{\omega_0}V`H{}_{\rho^0_0}\ltimes^{\omega_0'}V;\psi_1```\psi_0]
\efig$$
commutes if and only if 
$$\big{(}i(g),\phi(w+\lambda_1(g))+\check{\varphi}'(g)\big{)}=\big{(}i(g),\phi(w)+\check{\varphi}(g)+\lambda_0(i(g))\big{)}, $$
which is equivalent to Eq. (\ref{cohPhi}). On the other hand, the 
expressions 
$$\mld \psi_1\big{(}(g,w)^{(h,v)}\big{)} & =\psi_1\big{(}g^h,\rho^1_0(h)^{-1}(w+\rho_1(g)v)+\alpha(h;g)\big{)} \\
					   =\big{(}g^h,\rho^1_0(h)^{-1}(w+\rho_1(g)v)+\alpha(h;g)+\lambda_1(g^h)\big{)}, $$
and
$$\mld \psi_1(g,w)^{\psi_0(h,v)} & =\big{(}g,w+\lambda_1(g)\big{)}^{(h,v+\lambda_0(h))} \\
			   =\big{(}g^h,\rho^1_0(h)^{-1}\big{(}w+\lambda_1(g)+\rho_1(g)(v+\lambda_0(h))\big{)}+\alpha'(h;g)\big{)} $$
coincide if and only if Eq. (\ref{cohAlpha}) holds.
\epf

\rem\label{sanity2}
Continuing Remark \ref{sanityCheck}, one can also prove that, when 
appropriately restricted to the trivial $2$-representation on 
$(0)\to\Rr$, Proposition \ref{Gp2-coboundaries} is equivalent to Lemma 
\ref{isoEasyExt}. \endrem
 
If 
$\vec{\omega}=(0,\varphi,\omega_0,\lambda,\alpha,\omega_1)\in C^2_{tot}(\G,\phi)$ 
is a $2$-cocycle, then $(\nabla\vec{\omega})^{2,0}_1=0$ implies 
$\lambda=0$. Furthermore, $(\nabla\vec{\omega})^{2,1}_0=0$ reads 
$$\Delta_{2,1}\omega_1=\partial\varphi ; $$
which, evaluated at an arbitrary element 
$(g_0,g_1,h)\in G^2\times H\cong\G_2$, yields
$$\varphi\pmatrix{g_1g_0 \cr
h} = \varphi\pmatrix{g_1 \cr
h}+\varphi\pmatrix{g_0 \cr
hi(g_1)}-\phi\big{(}\rho_0^1(h)\omega_1(g_1,g_0)\big{)}. $$
In particular, making $g_0=1$, one sees that $\varphi$ vanishes 
identically on the space of units $H$. Defining
\begin{equation}\label{checkOutOfPhi}
\check{\varphi}(g):=\varphi\pmatrix{g \cr
1},
\end{equation}
one is left with a $4$-tuple 
$(\check{\varphi},\omega_0,\alpha,\omega_1)$ as in the statement of 
Proposition \ref{Gp2-cocycles}. In the sequel, we show that this 
$4$-tuple verifies the equations of Proposition \ref{Gp2-cocycles} if 
and only if $\nabla\vec{\omega}=0$, and thus defines an extension.
\begin{itemize}
\item $(\nabla\vec{\omega})^{0,3}_0=\delta\omega_0 =0$; therefore, 
$\omega_0$ is a $2$-cocycle with respect to $\rho_0^0$.  
\item $(\nabla\vec{\omega})^{0,0}_3=\delta_{(1)}\omega_1 =0$; therefore, 
$\omega_1$ is a $2$-cocycle with respect to $\rho_0^1\circ i$. 
\item Evaluating $(\nabla\vec{\omega})^{1,2}_0=0$ at 
$\pmatrix{g_1 & g_2 \cr
1   & 1}\in\G^2$ yields Eq. (\ref{iii)}). 
\item $(\nabla\vec{\omega})^{0,2}_1=0$ is exactly Eq. (\ref{iv)}). 
\item Evaluate $(\nabla\vec{\omega})^{1,2}_0=0$ at
$$\pmatrix{g & 1 \cr
h^{-1} & h}: \omega_0(h^{-1},h)-\omega_0(h^{-1}i(g),h)+\check{\varphi}(g^h)-\varphi\pmatrix{
g \cr
h^{-1}}+\phi(\alpha(h;g)) = 0, $$
and  
$$\pmatrix{1 & g \cr
h^{-1} & 1}: -\omega_0(h^{-1},i(g))-\rho_0^0(h)^{-1}\check{\varphi}(g)+\varphi\pmatrix{g \cr 
h^{-1}} = 0, $$
and consider their sum
$$-\omega_0(h^{-1},h)+\omega_0(h^{-1}i(g),h)+\omega_0(h^{-1},i(g))=\check{\varphi}(g^h)-\rho_0^0(h)^{-1}\check{\varphi}(g)+\phi(\alpha(h;g)). $$
Since $\delta\omega_0=0$,
$$\rho_0^0(h)^{-1}\omega_0(i(g),h)+\omega_0(h^{-1},i(g)h)=\omega_0(h^{-1}i(g),h)+\omega_0(h^{-1},i(g)); $$
hence, by substituting, one gets Eq. (\ref{v)}).
\item Evaluating $(\nabla\vec{\omega})^{1,1}_1=0$ at 
$(\gamma;g_1)\in\G\times G$, where $\gamma=\pmatrix{g_2 \cr
1}\in\G$ yields Eq. (\ref{vi)}). 
\item $(\nabla\vec{\omega})^{0,1}_2=0$ is exactly Eq. (\ref{vii)}).
\end{itemize}

Conversely, let $(\check{\varphi},\omega_0,\alpha,\omega_1)$ be a 
$4$-tuple as in the statement of Proposition \ref{Gp2-cocycles} verifying 
the equations therein. Define
\begin{equation}\label{phiOutOfCheck}
\varphi\pmatrix{g \cr
h}:=\omega_0(h,i(g))+\rho_0^0(h)\check{\varphi}(g) 
\end{equation}
and set $\vec{\omega}=(0,\varphi,\omega_0,0,\alpha,\omega_1)$. We 
proceed to show that $\vec{\omega}$ is a $2$-cocycle. From the previous 
discussion, it suffices to check that $(\nabla\vec{\omega})^{2,1}_0$, 
$(\nabla\vec{\omega})^{1,1}_1$ and $(\nabla\vec{\omega})^{1,2}_0$ vanish. 
Indeed, let $(g_1,g_2,h)\in G^2\times H$, then from the cocycle equation 
$\delta\omega_0=0$ and Eq. (\ref{iii)}), 
$$\omega_0(h,i(g_2g_1))=-\rho_0^0(h)\omega_0(i(g_2),i(g_1))+\omega_0(hi(g_2),i(g_1))+\omega_0(h,i(g_2)); $$
and
$$\mld \varphi\pmatrix{g_2g_1 \cr
h} & = \omega_0(h,i(g_2g_1))+\rho_0^0(h)\check{\varphi}(g_2g_1)\\ 
     = \omega_0(hi(g_2),i(g_1))+\omega_0(h,i(g_2))+\rho_0^0(h)\big{[}\rho_0^1(i(g_2))\check{\varphi}(g_1)+\check{\varphi}(g_2)-\phi(\omega_1(g_2,g_1))\big{]} \\
     =\varphi\pmatrix{g_1 \cr
hi(g_2)}+\varphi\pmatrix{g_2 \cr
h}-\rho_0^0(h)\circ\phi(\omega_1(g_2,g_1)); $$
hence, $(\nabla\vec{\omega})^{2,1}_0=0$. Next, let 
$(\gamma;f)\in\G\times G$ and $(g,h)\in G\rtimes H$ the image of 
$\gamma$ under the isomorphism of Remark \ref{Equiv2Gp:Mods}. Then,
$$\mld \delta'\varphi(\gamma;f) & =\rho_0^1(hi(g))^{-1}\rho_1(f)\varphi(\gamma) 
                                  =\rho_0^1(hi(g))^{-1}\rho_1(f)\big{(}\omega_0(h,i(g))+\rho_0^0(h)\check{\varphi}(g)\big{)} \\
                             =\rho_0^1(hi(g))^{-1}\rho_1(f)\omega_0(h,i(g))+\rho_0^1(i(g))^{-1}\rho_1(f^h)\check{\varphi}(g). $$
Using Eq.'s (\ref{iv)}) and (\ref{vi)}), one gets
$$\mld \delta'\varphi(\gamma;f) & =\rho_0^1(i(g))^{-1}\alpha(h;f)-\alpha(hi(g);f)+\alpha(i(g);f^h)+ \\
                   \qquad\quad -\alpha(i(g),f^h)+\rho_0^1(i(g))^{-1}\omega_1(f^h,g)+\omega_1(g^{-1},f^hg)-\omega_1(g^{-1},g); $$
hence, $(\nabla\vec{\omega})^{1,1}_1=0$. Concludingly, let 
$(\gamma_1,\gamma_2)\in\G^2$ and let $(g_k,h_k)\in G\rtimes H$ be the 
image of $\gamma_k$ under the isomorphism of Remark \ref{Equiv2Gp:Mods} 
for $k\in\lbrace 1,2\rbrace$. Then,
$$\mld \delta\varphi\pmatrix{g_1 & g_2 \cr
h_1 & h_2} & =\rho_0^0(h_1i(g_1))\varphi\pmatrix{
g_2 \cr
h_2}-\varphi\pmatrix{
g_1^{h_2}g_2 \cr
h_1h_2}+\varphi\pmatrix{g_1 \cr
h_1}\\
    = \rho_0^0(h_1i(g_1))\big{[}\omega_0(h_2,i(g_2))+\rho_0^0(h_2)\check{\varphi}(g_2)\big{]}+ \\
      \quad -\big{(}\omega_0(h_1h_2,i(g_1^{h_2}g_2))+\rho_0^0(h_1h_2)\check{\varphi}(g_1^{h_2}g_2)\big{)}+\omega_0(h_1,i(g_1))+\rho_0^0(h_1)\varphi(g_1) $$
Using successively Eq.'s (\ref{iii)}) and (\ref{v)}), one computes
$$\mld \check{\varphi}(g_1^{h_2}g_2) & =\rho_0^0(i(g_1^{h_2}))\check{\varphi}(g_2)+\check{\varphi}(g_1^{h_2})+\omega_0(i(g_1^{h_2}),i(g_2))-\phi(\omega_1(g_1^{h_2},g_2)) \\
  =\rho_0^0(i(g_1^{h_2}))\check{\varphi}(g_2)+\rho_0^0(h_2)^{-1}\omega_0(i(g_1),h_2)+\omega_0(h_2^{-1},i(g_1)h_2)-\omega_0(h_2^{-1},h_2)+ \\
  \qquad +\rho_0^0(h_2)^{-1}\check{\varphi}(g_1)-\phi(\alpha(h_2;g_1))+\omega_0(i(g_1^{h_2}),i(g_2))-\phi(\omega_1(g_1^{h_2},g_2)). $$
Substituting,
$$\mld \delta\varphi\pmatrix{g_1 & g_2 \cr
h_1 & h_2} & =\rho_0^0(h_1h_2)\phi(\alpha(h_2;g_1)+\omega_1(g_1^{h_2},g_2))+R(\omega_0) \\
             =\Delta\alpha\pmatrix{g_0 & g_1 \cr
h_0 & h_1}+\Delta_{1,2}\omega_1\pmatrix{g_0 & g_1 \cr
h_0 & h_1}+R(\omega_0), $$
where 
$$\mld R(\omega_0) & =\rho_0^0(h_1i(g_1))\omega_0(h_2,i(g_2))-\omega_0(h_1h_2,i(g_1^{h_2}g_2))+\omega_0(h_1,i(g_1))-\rho_0^0(h_1)\omega_0(i(g_1),h_2)+ \\
             \qquad -\rho_0^0(h_1h_2)\big{(}\omega_0(h_2^{-1},i(g_1)h_2)-\omega_0(h_2^{-1},h_2)+\omega_0(i(g_1^{h_2}),i(g_2))\big{)}. $$
We claim that $R(\omega_0)=\partial\omega_0\pmatrix{g_1 & g_2 \cr
h_1 & h_2}$, and consequently $(\nabla\vec{\omega})^{1,2}_0=0$. Indeed, 
this claim follows from repeatedly applying the cocycle equation to 
rewrite $R(\omega_0)$. Using $\delta\omega_0=0$ evaluated at 
$(h_2^{-1},i(g_1)h_2,i(g_2))$ yields
$$\mld R(\omega_0) & =\rho_0^0(h_1i(g_1))\omega_0(h_2,i(g_2))-\omega_0(h_1h_2,i(g_1^{h_2}g_2))+\omega_0(h_1,i(g_1))-\rho_0^0(h_1)\omega_0(i(g_1),h_2)+ \\
             \qquad -\rho_0^0(h_1)\omega_0(i(g_1)h_2,i(g_2))-\rho_0^0(h_1h_2)\big{(}\omega_0(h_2^{-1},i(g_1)h_2i(g_2))-\omega_0(h_2^{-1},h_2)\big{)}. $$
Successively, using $\delta\omega_0=0$ evaluated at 
$(h_1h_2,h_2^{-1},i(g_1)h_2,i(g_2))$, $(h_1,i(g_1)h_2,i(g_2))$ 
and $(h_1,i(g_1),h_2)$ yields
$$\mld R(\omega_0) & =\rho_0^0(h_1i(g_1))\omega_0(h_2,i(g_2))-\omega_0(h_1,i(g_1)h_2i(g_2))-\omega_0(h_1h_2,h_2^{-1})+\omega_0(h_1,i(g_1))+ \\
             \qquad -\rho_0^0(h_1)\big{(}\omega_0(i(g_1),h_2)+\omega_0(i(g_1)h_2,i(g_2))\big{)}+\rho_0^0(h_1h_2)\omega_0(h_2^{-1},h_2) \\
             		 =\rho_0^0(h_1i(g_1))\omega_0(h_2,i(g_2))-\omega_0(h_1i(g_1)h_2,i(g_2))-\omega_0(h_1h_2,h_2^{-1})+\omega_0(h_1,i(g_1))+ \\
             \qquad -\omega_0(h_1,i(g_1)h_2)-\rho_0^0(h_1)\omega_0(i(g_1),h_2)+\rho_0^0(h_1h_2)\omega_0(h_2^{-1},h_2) \\
             		 =\rho_0^0(h_1i(g_1))\omega_0(h_2,i(g_2))-\omega_0(h_1i(g_1)h_2,i(g_2))-\omega_0(h_1h_2,h_2^{-1})+ \\
             \qquad -\omega_0(h_1i(g_1),h_2)+\rho_0^0(h_1h_2)\omega_0(h_2^{-1},h_2).  $$
Thus, the claim $R(\omega_0) = -\omega_0(h_1i(g_1),h_2i(g_2))+\omega_0(h_1,h_2)$ 
follows from $\delta\omega_0(h_1i(g_1),h_2,i(g_2))=\delta\omega_0(h_1h_2,h_2^{-1},h_2)=0$.

Furthermore, if  
$$(\vec{\omega})^k=(0,\varphi^k,\omega_0^k,\lambda^k,\alpha^k,\omega_1^k), \qquad\qquad k\in\lbrace 1,2\rbrace $$
are a pair of cohomologous of $2$-cocycles, say
$$(\vec{\omega})^2-(\vec{\omega})^1=\nabla(v,\lambda_0,\lambda_1), $$
then coordinate-wise one has got
$$\omega_0^2-\omega_0^1 = \delta\lambda_0 $$
$$\varphi^2-\varphi^1 = \partial\lambda_0-\delta v+\Delta\lambda_1 \qquad\qquad \alpha^2-\alpha^1 = \delta'\lambda_0-\delta\lambda_1 $$
$$0 = -\partial v=-v \qquad\qquad 0 = -\delta'v \qquad\qquad \omega_1^2-\omega_1^1 = \delta_{(1)}\lambda_1 . $$
Explicitly, for $\gamma=(g,h)\in G\rtimes H\cong\G$ and 
$(h;f)\in H\times G$
\begin{equation}\label{2Coc1}
(\varphi^2-\varphi^1)(\gamma) = \lambda_0(h)-\lambda_0(hi(g))+\phi(\rho_0^1(h)\lambda_1(g))
\end{equation}
\begin{equation}\label{2Coc2}
(\alpha^2-\alpha^1)(h;f) = \rho_0^1(h)^{-1}\rho_1(f)\lambda_0(h)-\big{(}\lambda_1(f^h)-\rho_0^1(h)^{-1}\lambda_1(f)\big{)}.
\end{equation}
Eq. (\ref{2Coc1}) evaluated at $\gamma=(g,1)$ yields Eq. (\ref{cohPhi}) 
with $\check{\varphi}^k$ defined by Eq. (\ref{checkOutOfPhi}), 
and Eq. (\ref{2Coc2}) is exactly Eq. (\ref{cohAlpha}). Conversely, given 
two $4$-tuples as in the statement of Proposition \ref{Gp2-coboundaries} 
whose difference verifies the equations therein and defining a 
corresponding pair of $\varphi^k$'s by Eq. (\ref{phiOutOfCheck}),
$$\mld (\varphi^2-\varphi^1)(\gamma) & =\omega_0^2(h,i(g))+\rho_0^0(h)\varphi^2(g)-\omega_0^1(h,i(g))-\rho_0^0(h)\varphi^1(g) \\    
           =\delta\lambda_0(h,i(g))+\rho_0^0(h)(\phi(\lambda_1(g))-\lambda_0(i(g))) \\
           =\rho_0^0(h)\lambda_0(i(g))-\lambda_0(hi(g))+\lambda_0(h)+\phi(\rho_0^1(h)\lambda_1(g))-\rho_0^0(h)\lambda_0(i(g)) \\
           =\partial\lambda_0(\gamma)+\Delta\lambda_1(\gamma). $$  

Summing this discussion up, we have shown the following:
\thm\label{H2Gp}
$H^2_\nabla(\G,\phi)$, the second degree cohomology of the subcomplex of 
$C_{tot}(\G,\phi)$ that has trivial $(^{p,0}_0)$-coordinate for $p>0$, 
is in one-to-one correspondence with split extensions of the Lie 
$2$-group $\G$ by the $2$-vector space $W\to^\phi V$. \endthm

\rem\label{scope}
Because of the application that we have in mind (see the Introduction), 
Theorem \ref{H2Gp} is stated and proved to classify extensions of $\G$ 
by $2$-vector spaces. However, do notice that the results of this 
section can be effortlessly extended to more general abelian $2$-groups: 
Replacing the $2$-representation by an action of $\G$ on an abelian 
$2$-group (see \cite{Norrie}), one forms a truncated complex using the 
same formulae. The second cohomology of the complex with values in the 
abelian $2$-group $\mathbbm{A}$ is seen to classify extensions by 
$\mathbbm{A}$ simply reinterpreting the arguments we presented. \endrem

\section{Inkling of a larger complex}\label{sec-Ink}

In this section, we study the grid of Section \ref{sec-tetrahedralCx}. 
We prove that for either $q$ or $p$ constant, there is a double complex. 
More precisely, we prove that $\partial$ and $\delta_{(1)}$ do commute, 
yielding a sequence of double complexes that we call the $q$-pages, and 
show that $\delta$ and $\delta_{(1)}$ also commute, thus yielding a 
sequence of double complexes that we call the $p$-pages. We prove that 
there is a general formula for the {\em difference maps} $\Delta$ 
that measure the difference between $\partial\delta$ and 
$\delta\partial$. In fact, using $\Delta$ we show that the two 
dimensional grid $r=0$ commutes up to isomorphism in the $2$-vector 
space (cf. Proposition \ref{FrontPage-upToIso}), and the successive 
$r$-\textit{pages} commute up to homotopy (cf. Proposition \ref{r-page}). 
Since $\partial$ and $\delta$ fail to commute on the nose, the total 
differential $\nabla$ (cf. Eq. (\ref{preDiff})) does not, in general, 
square to zero. In Section \ref{sec-tetrahedralCx}, it is explained how 
this defect is partially solved by adding the difference maps; however, 
the updated $\nabla$ still fails to square to zero. In a similar fashion, 
the coordinates that fail to vanish, do so because of some differentials 
that do not commute, but that do commute up to higher homotopies. 
Accordingly, we call these homotopies  {\em higher difference maps} and 
note them $\Delta_{a,b}$ (cf. Eq.'s (\ref{Delta21}) and (\ref{Delta12})). 
We provide explicit formulas for some families of higher difference maps 
and explain what the grid is lacking to form a complex.

Ultimately, lacking additional structure such as a product or a more 
insightful definition of the higher difference maps, the proofs of the 
relations we present boil down to unfortunately long and rather 
unenlightening computations.

\rem\label{expectations}
In this section, we collect evidence suggesting there is a full triple 
complex analog to that in \cite{Lie2AlgCoh}. Eventually, one could 
collapse one dimension in the triple complex, by taking, {\em e.g.}, the 
total complex of the $p$-pages; however, as opposed to the case of 
\cite{Lie2AlgCoh}, the remaining object is not a double complex 
because of the higher difference maps $\Delta_{a,b}$ with $a>1$. \endrem

Throughout, let $\G$ be a Lie $2$-group with associated crossed module 
$G\to^iH$ and let $\rho$ be a $2$-representation of $\G$ on the 
$2$-vector space $W\to^\phi V$. Using the notation conventions laid down 
in Subsection \ref{subsec-rep}, we think of $\rho$ as a triple 
$(\rho_0^0,\rho_0^1;\rho_1)$. Also, we take the isomorphisms of Remarks 
\ref{Equiv2Gp:Mods} and \ref{2GpNerve} to be fixed and we abuse notation 
and often treat them as equalities. Recall the notation from Section 
\ref{sec-tetrahedralCx}: We write
$$\bfig
\morphism<850,0>[\partial:C^{p,q}_r(\G,\phi)`C^{p+1,q}_r(\G,\phi);]
\efig$$
for the differential in the $p$-direction (cf. \ref{p-dir}),
$$\bfig
\morphism<850,0>[\delta:C^{p,q}_r(\G,\phi)`C^{p,q+1}_r(\G,\phi);]
\efig$$
for the differential in the $q$-direction (cf. \ref{q-dir}), and
$$\bfig
\morphism<900,0>[\delta_{(1)}:C^{p,q}_r(\G,\phi)`C^{p,q}_{r+1}(\G,\phi);]
\efig$$
for the differential in the $r$-direction (cf. \ref{r-dir}).

\subsection{Commuting differentials}\label{pq-pag} In this subsection, 
we prove that for constant $q$ the two dimensional grid of maps is a 
double complex that we call the $q$-{\em page}. Analogously, for 
constant $p$, there is a double complex that we call the $p$-{\em page}.

\subsubsection*{$q$-pages} When $q=0$, there is a double complex all of 
whose columns are equal and the intertwining maps are either zero or the 
identity, thus commuting trivially. For $q>0$, there is a double complex 
as well.

\proposition\label{q-pages}
For each $q>0$, 
$$\bfig
\node 1a(0,1350)[\vdots]
\node 1b(850,1350)[\vdots]
\node 1c(1800,1350)[\vdots]
\node 1d(2400,1350)[]
\node 2a(0,1000)[C(H^q\times G^2,W)]
\node 2b(850,1000)[C(\G^q\times G^2,W)]
\node 2c(1800,1000)[C(\G_2^q\times G^2,W)]
\node 2d(2400,1000)[\cdots]
\node 3a(0,500)[C(H^q\times G,W)]
\node 3b(850,500)[C(\G^q\times G,W)]
\node 3c(1800,500)[C(\G_2^q\times G,W)]
\node 3d(2400,500)[\cdots]
\node 4a(0,0)[C(H^q,V)]
\node 4b(850,0)[C(\G^q,V)]
\node 4c(1800,0)[C(\G_2^q,V)]
\node 4d(2400,0)[\cdots]
\arrow|l|[4a`3a;\delta']
\arrow[4a`4b;\partial]
\arrow|l|[3a`2a;\delta_{(1)}]
\arrow[3a`3b;\partial]
\arrow|l|[4b`3b;\delta']
\arrow[4b`4c;\partial]
\arrow[2a`1a;]
\arrow[2a`2b;\partial]
\arrow|l|[3b`2b;\delta_{(1)}]
\arrow[3b`3c;\partial]
\arrow|l|[4c`3c;\delta']
\arrow[4c`4d;]
\arrow[2b`1b;]
\arrow[2b`2c;\partial]
\arrow|l|[3c`2c;\delta_{(1)}]
\arrow[3c`3d;]
\arrow[2c`1c;]
\arrow[2c`2d;]
\efig$$
is a double complex. \endproposition
 
\pf
Each row is a complex and due to Lemma \ref{r-cx}, so is each column. 
Moreover, due to Lemma \ref{1stqSquare}, we just need to prove that 
$\partial$ commutes with $\delta_{(1)}$.

Let $\omega\in C(\G_p^q\times G^r,V)$, 
$\vec{\gamma}=(\gamma_1,...,\gamma_q)^T\in\G_{p+1}^q$ and 
$\vec{f}=(f_0,...,f_r)\in G^{r+1}$. Adopting the convention of Lemma 
\ref{1stqSquare} for $\gamma_b\in\G_{p+1}$, 
$$ \delta_{(1)}\partial\omega(\vec{\gamma};\vec{f})=\rho_0^1(i(f_0^{t_{p+1}(\gamma_1)...t_{p+1}(\gamma_q)}))\partial\omega(\vec{\gamma};\delta_0\vec{f})+\sum_{k=1}^{r+1}(-1)^k\partial\omega(\vec{\gamma};\delta_k\vec{f}) $$
$$\mld & =\rho_0^1(i(f_0^{t_{p+1}(\gamma_1)...t_{p+1}(\gamma_q)}))\Big{(}\rho_0^1(i(pr_G(\gamma_{10}\vJoin ...\vJoin\gamma_{q0})))^{-1}\omega(\partial_0\vec{\gamma};\delta_0\vec{f})+\sum_{j=1}^{p+1}(-1)^j\omega(\partial_j\vec{\gamma};\delta_0\vec{f})\Big{)}+ \\
           \qquad +\sum_{k=1}^{r+1}(-1)^k\Big{(}\rho_0^1(i(pr_G(\gamma_{10}\vJoin ...\vJoin\gamma_{q0})))^{-1}\omega(\partial_0\vec{\gamma};\delta_k\vec{f})+\sum_{j=1}^{p+1}(-1)^j\omega(\partial_j\vec{\gamma};\delta_k\vec{f})\Big{)}. $$
Now, given that 
$$\mld \rho_0^1\Big{(}i\big{(}f_0^{t_{p+1}(\gamma_1)...t_{p+1}(\gamma_q)}\big{)}\Big{)} & =\rho_0^1\Big{(}i\big{(}f_0^{t(\gamma_{10})...t(\gamma_{q0})}\big{)}\Big{)} 
     =\rho_0^1\Big{(}i\big{(}f_0^{t(\gamma_{10}\vJoin ...\vJoin\gamma_{q0})}\big{)}\Big{)} \\
     =\rho_0^1\Big{(}i\big{(}f_0^{h_{10}...h_{q0}i(pr_G(\gamma_{10}\vJoin ...\vJoin\gamma_{q0}))}\big{)}\Big{)} \\
     =\rho_0^1(i(pr_G(\gamma_{10}\vJoin ...\vJoin\gamma_{q0})))^{-1}\rho_0^1(i(f_0^{h_{10}...h_{q0}}))\rho_0^1(i(pr_G(\gamma_{10}\vJoin ...\vJoin\gamma_{q0}))), $$
one can group terms to get 
$$\mld \delta_{(1)}\partial\omega(\vec{\gamma};\vec{f}) & =\rho_0^1(i(pr_G(\gamma_{10}\vJoin ...\vJoin\gamma_{q0})))^{-1}\Big{(}\rho_0^1(i(f_0^{h_{10}...h_{q0}}))\omega(\partial_0\vec{\gamma};\delta_0\vec{f})+\sum_{k=1}^{r+1}(-1)^k\omega(\partial_0\vec{\gamma};\delta_k\vec{f})\Big{)} \\
       \qquad +\sum_{j=1}^{p+1}(-1)^j\Big{(}\rho_0^1(i(f_0^{t_{p+1}(\gamma_1)...t_{p+1}(\gamma_q)}))\omega(\partial_j\vec{\gamma};\delta_0\vec{f})+\sum_{k=1}^{r+1}(-1)^k\omega(\partial_j\vec{\gamma};\delta_k\vec{f})\Big{)}. $$
If $j>0$, $t_p(\partial_j\gamma_b)=t(\gamma_{b0})=t_{p+1}(\gamma_b)$; 
otherwise, $t_p(\partial_0\gamma_b)=t(\gamma_{b1})=s(\gamma_{b0})=h_{b0}$. 
Hence, 
$$\delta_{(1)}\partial\omega(\vec{\gamma};\vec{f})=\rho_0^1(i(pr_G(\gamma_{10}\vJoin ...\vJoin\gamma_{q0})))^{-1}\delta_{(1)}\omega(\partial_0\vec{\gamma};\vec{f})+\sum_{j=1}^{p+1}(-1)^j\delta_{(1)}\omega(\partial_j\vec{\gamma};\vec{f})= \partial\delta_{(1)}\omega(\vec{\gamma};\vec{f}). $$ \epf

\subsubsection*{$p$-pages} We start by outlining some general facts 
about actions of Lie groups and representations. 

\lem\label{aDblGpd}
Let $X$ and $Y$ be Lie groups, and let $Y$ act on the right on $X$ by 
Lie group automorphisms. Then 
\begin{equation}\label{aDbl}
\bfig
\square/@{>}@<3pt>`@{>}@<3pt>`@{>}@<3pt>`@{>}@<3pt>/[Y\ltimes X`Y`X`\ast,;```]
\square/@{>}@<-3pt>`@{>}@<-3pt>`@{>}@<-3pt>`@{>}@<-3pt>/[Y\ltimes X`Y`X`\ast,;```]
\efig
\end{equation}
where the top groupoid is a bundle of Lie groups and the left groupoid 
is the (right!) action groupoid for the action of $Y$ on $X$, i.e, 
\begin{eqnarray*}
\Lf{s}(y;x)=x^{y},
\end{eqnarray*}
for $x\in X$ and $y\in Y$, is a double Lie groupoid. \endlem

\rem\label{OffbeatNot} 
A word on the odd choice of notation: we opted out of calling our groups 
$G$ and $H$, so to avoid any concurrence and to help the reader stray 
away from believing that there is a Lie $2$-group somewhere in this 
diagram. There is none. In fact, though the action of $Y$ on $X$ is by 
automorphisms, there is no crossed module around. \endrem

\pf[of Lemma \ref{aDblGpd}]
Given each of the sides of the square (\ref{aDbl}) is clearly a Lie 
groupoid, to prove that the array (\ref{aDbl}) is a double Lie groupoid, 
we show that the horizontal structural maps are groupoid morphisms. 
First, the square commutes as all compositions land in $\ast$. Both 
$\Tp{s}=\Tp{t}$ are functors: for $x\in X$ and $y_1,y_2\in Y$, 
$$\Tp{s}(1;x)=1\quad\textnormal{and}\quad\Tp{s}((y_1;x)\vJoin(y_2;x^{y_1}))=\Tp{s}(y_1y_2;x)=y_1y_2=\Tp{s}(y_1;x)\Tp{s}(y_2;x^{y_1}). $$
$\Tp{u}$ is a functor as well: for $y\in Y$, 
$$\Lf{s}(\Tp{u}(y))=\Lf{s}(y;1)=1^{y}=1=u(\ast) \qquad\textnormal{and}\qquad \Lf{t}(\Tp{u}(y))=\Lf{t}(y;1)=1=u(\ast), $$
so it is well-defined, and it respects units and the multiplication:
$$\Tp{u}(1)=(1;1)\quad\textnormal{and}\quad\Tp{u}(y_1y_2)=(y_1y_2;1)=(y_1;1)\vJoin(y_2;1^{y_1})=\Tp{u}(y_1)\vJoin\Tp{u}(y_2). $$
To conclude, 
$(Y\ltimes X)_{\Tp{s}}\times_{\Tp{t}}(Y\ltimes X)\cong Y\ltimes X^2\two X^2$ 
is the right action groupoid for the diagonal action, we prove that 
$\Tp{m}$ is a Lie groupoid homomorphism: for $x_1,x_2\in X$ and $y\in Y$,  
$$\Lf{s}(\Tp{m}(y;x_1,x_2))=\Lf{s}(y;x_1x_2)=(x_1x_2)^y=x_1^yx_2^y=m(\Lf{s}^2(y;x_1,x_2)) $$
and
$$\Lf{t}(\Tp{m}(y;x_1,x_2))=\Lf{t}(y;x_1x_2)=x_1x_2=m(\Lf{t}^2(y;x_1,x_2)).$$
Furthermore, $\Tp{m}$ is also compatible with units and the 
multiplication, as $\Tp{m}(1;x_1,x_2)=(1;x_1x_2)$ and 
$((y_1;x_1)\vJoin(y_2;x_1^{y_1}))\Join ((y_1;x_2)\vJoin(y_2;x_2^{y_1}))$ 
equals
$$\mld (y_1y_2;x_1)\Join (y_1y_2;x_2) & =\Tp{m}(y_1y_2;x_1,x_2)=(y_1y_2;x_1x_2) \\
                                  =(y_1;x_1x_2)\vJoin(y_2;(x_1x_2)^{y_1}) 
                                  =\Tp{m}(y_1;x_1,x_2)\vJoin\Tp{m}(y_2;x_1^{y_1},x_2^{y_1}). $$
\epf
                                  
\example\label{TheDoubleOfGL} 
Let $W$ be a vector space and let $X=Y=GL(W)$ along with the right 
action by conjugation on itself. Then 
$$\bfig
\square/@{>}@<3pt>`@{>}@<3pt>`@{>}@<3pt>`@{>}@<3pt>/<900,450>[GL(W)\ltimes GL(W)`GL(W)`GL(W)`\ast;```]
\square/@{>}@<-3pt>`@{>}@<-3pt>`@{>}@<-3pt>`@{>}@<-3pt>/<900,450>[GL(W)\ltimes GL(W)`GL(W)`GL(W)`\ast;```]
\efig$$
is a double Lie groupoid. \endexample

\lem\label{p-pageModel}
Let $X$ and $Y$ be as in the statement of Lemma \ref{aDblGpd} and let 
\begin{equation}\label{dblGpdMap}
Y\ltimes X\to<350>GL(W)\ltimes GL(W):(y;x)\to/|->/<350>(\rho_Y(y),\rho_X(x))
\end{equation} 
be a map of double Lie groupoids from the double Lie groupoid of Lemma 
\ref{aDblGpd} to the one in Example \ref{TheDoubleOfGL}. Then, for every 
pair of integers $q,r\geq 0$, there are representations of the Lie group 
bundles
\begin{equation}\label{GpBdles}
Y^q\times X\two Y^q
\end{equation}
on $Y^q\times W\to Y^q$ and of the right transformation groupoids 
\begin{equation}\label{ActGpds}
Y\ltimes X^r\two X^r
\end{equation}
on $X^r\times W \to X^r$ that make the grid
$$\bfig
\node 1a(0,1350)[\vdots]
\node 1b(850,1350)[\vdots]
\node 1c(1800,1350)[\vdots]
\node 1d(2400,1350)[]
\node 2a(0,1000)[C(Y^2,W)]
\node 2b(850,1000)[C(Y^2\times X,W)]
\node 2c(1800,1000)[C(Y^2\times X^2,W)]
\node 2d(2400,1000)[\cdots]
\node 3a(0,500)[C(Y,W)]
\node 3b(850,500)[C(Y\times X,W)]
\node 3c(1800,500)[C(Y\times X^2,W)]
\node 3d(2400,500)[\cdots]
\node 4a(0,0)[W]
\node 4b(850,0)[C(X,W)]
\node 4c(1800,0)[C(X^2,W)]
\node 4d(2400,0)[\cdots]
\arrow|l|[4a`3a;\delta]
\arrow[4a`4b;\delta']
\arrow|l|[3a`2a;\delta]
\arrow[3a`3b;\delta']
\arrow|l|[4b`3b;\delta]
\arrow[4b`4c;\delta']
\arrow[2a`1a;]
\arrow[2a`2b;\delta']
\arrow|l|[3b`2b;\delta]
\arrow[3b`3c;\delta']
\arrow|l|[4c`3c;\delta]
\arrow[4c`4d;]
\arrow[2b`1b;]
\arrow[2b`2c;\delta']
\arrow|l|[3c`2c;\delta]
\arrow[3c`3d;]
\arrow[2c`1c;]
\arrow[2c`2d;]
\efig$$
whose rows and columns are Lie groupoid cochain complexes taking values 
in these representations into a double complex. \endlem

\pf
Since (\ref{dblGpdMap}) is a map of double Lie groupoids, its 
restrictions to the bottom and right groupoids give respectively 
representations $\rho_X$ of $X$ and $\rho_Y$ of $Y$, both on $W$. These 
are the representations for the group bundle over a point $X\two\ast$ 
and for the trivial transformation groupoid $Y\ltimes\ast\two\ast$. For 
each $q>0$, the representation of the group bundle (\ref{GpBdles}) is the 
pull-back of $\rho_X$ along the groupoid homomorphism $\Lf{s}_q$. In 
symbols, define the representation $\rho_X^q$ of the group bundle 
(\ref{GpBdles}) on $(s_Y)_q^*W=Y^q\times W\to Y^q$ by
$$\rho_X^q(y_1,...,y_q;x):=\Lf{s}_q^*\rho_X(y_1,...,y_q;x)=\rho_X(x^{y_1...y_q}) $$
for $(y_1,...,y_q;x)\in Y^q\times X\cong(Y\ltimes X)^{(q)}$. Analogously, 
for each $r>0$, the representation of the transformation groupoid 
(\ref{ActGpds}) is the pull-back of $\rho_Y$ along the groupoid 
homomorphisms $\Tp{t}_r$. In symbols, define the right representation 
$\rho_Y^r$ of the transformation groupoid (\ref{ActGpds}) on 
$t_r^*W=X^r\times W\to X^r$ by 
$$\rho_Y^r(y;x_1,...,x_r):=(\iota_Y\circ\Tp{t}_r)^*\rho_Y(y;x_1,...,x_r)=\rho_Y(y)^{-1} $$
for $(y;x_1,...,x_r)\in Y\times X^r\cong(Y\times X)^{(r)}$.

If $y\in Y$ and $\vec{x}=(x_1,...,x_r)\in X^r$, we adopt the convention 
that $(\vec{x})^{y}:=(x_1^y,...,x_r^y)$. 

We proceed to check that, for fixed $(q,r)$, the spaces of cochains of 
the groupoids (\ref{GpBdles}) and (\ref{ActGpds}) with respect to 
$\rho_X^q$ and $\rho_Y^r$ concur. On the one hand,
$$(Y^q\times X)^{(r)}=\lbrace (\vec{y}_1,x_1;...;\vec{y}_r,x_r)\in(Y^q\times X)^r:\Tp{s}(\vec{y}_j,x_j)=\vec{y}_j=\vec{y}_{j+1}=\Tp{t}(\vec{y}_{j+1},x_{j+1})\rbrace ; $$
therefore, $(Y^q\times X)^{(r)}\cong Y^q\times X^r$, where the 
diffeomorphism is obviously given by 
$(\vec{y},x_1;...;\vec{y},x_r)\to/|->/(\vec{y};x_1,...,x_r)$. On the 
other hand,  
$$(Y\ltimes X^r)^{(q)}=\lbrace (y_1,\vec{x}_1;...;y_q,\vec{x}_q)\in(Y\ltimes X^r)^{q}:\Lf{s}(y_j,\vec{x}_j)=(\vec{x}_j)^{y_j}=\vec{x}_{j+1}=\Lf{t}(y_{j+1},\vec{x}_{j+1})\rbrace ; $$
therefore, $(Y\ltimes X^r)^{(q)}\cong Y^q\times X^r$, where the 
diffeomorphism is naturally given by 
$(y_1,\vec{x};y_2,(\vec{x})^{y_1};...;y_q,(\vec{x})^{y_1...y_{q-1}})\to/|->/(y_1,...,y_q;\vec{x})$. 
Since the representations $\rho_X$ and $\rho_Y$ take values on a trivial 
vector bundle, the pull-back bundle along any groupoid homomorphism 
remains trivial and its sections are but smooth functions to $W$; thus, 
$$C^r(Y^q\times X;Y^q\times W)=C(Y^q\times X^r,W)=C^q(Y\ltimes X^r;X^r\times W).$$

We use the diffeomorphisms of the latter discussion and the face maps 
$\delta_j$ of the Lie group $Y$ and $\delta'_k$ of the Lie group $X$ to 
rewrite the face maps of the groupoids (\ref{GpBdles}) and 
(\ref{ActGpds}): For 
$\vec{y}=(y_0,...,y_q)\in Y^{q+1}$ and 
$\vec{x}=(x_0,...,x_r)\in X^{r+1}$,    
$$\delta_j(\vec{y};\vec{x})=\cases{
    (\delta_0\vec{y};(\vec{x})^{y_0}) & if $j=0$ \cr
    (\delta_{j}\vec{y};\vec{x})       & otherwise } \quad\textnormal{and}\quad 
    \delta'_k(\vec{y};\vec{x})=(\vec{y};\delta'_k\vec{x}). $$

We are left to prove that the generic square
\begin{equation}\label{generic p-square}
\bfig
\square/>`<-`<-`>/<1200,450>[C(Y^{q+1}\times X^r,W)`C(Y^{q+1}\times X^{r+1},W)`C(Y^q\times X^r,W)`C(Y^q\times X^{r+1},W);\delta'`\delta`\delta`\delta']
\efig
\end{equation}
commutes. Indeed, let $\omega\in C(Y^q\times X^r,W)$ and $\vec{y}$ and 
$\vec{x}$ be as above. Then, 
$$\delta'\delta\omega(\vec{y};\vec{x})=\rho_X^{q+1}(\vec{y};x_0)\delta\omega(\vec{y};\delta'_0\vec{x})+\sum_{k=1}^{r+1}(-1)^{k}\delta\omega(\vec{y};\delta'_k\vec{x}), $$
while
$$\delta\delta'\omega(\vec{y};\vec{x})=\delta'\omega(\delta_0\vec{y};(\vec{x})^{y_0})+\sum_{j=1}^{q}(-1)^{j}\delta'\omega(\delta_j\vec{y};\vec{x})+(-1)^{q+1}\rho_Y^{r+1}(y_q;\vec{x})\delta'\omega(\delta_{q+1}\vec{y};\vec{x}). $$
We expand further to make the common terms evident:
$$\mld \delta'\delta\omega(\vec{y};\vec{x}) & =\rho_X^{q+1}(\vec{y};x_0)\Big{[}\omega(\delta_0\vec{y};(\delta'_0\vec{x})^{y_0})+\sum_{j=1}^{q}(-1)^{j}\omega(\delta_j\vec{y};\delta'_0\vec{x})+(-1)^{q+1}\rho_Y^r(y_q;\delta'_0\vec{x})\omega(\delta_{q+1}\vec{y};\delta'_0\vec{x})\Big{]}+ \\
									  \quad +\sum_{k=1}^{r+1}(-1)^{k}\Big{[}\omega(\delta_0\vec{y};(\delta'_k\vec{x})^{y_0})+\sum_{j=1}^{q}(-1)^{j}\omega(\delta_j\vec{y};\delta'_k\vec{x})+(-1)^{q+1}\rho_Y^r(y_q;\delta'_k\vec{x})\omega(\delta_{q+1}\vec{y};\delta'_k\vec{x})\Big{]} $$
and
$$\mld \delta\delta'\omega(\vec{y};\vec{x}) & =\rho_X^q(\delta_0\vec{y};x_0^{y_0})\omega(\delta_0\vec{y};\delta'_0(\vec{x})^{y_0})+\sum_{k=1}^{r+1}(-1)^{k}\omega(\delta_0\vec{y};\delta'_k(\vec{x})^{y_0})+ \\
									  \qquad +\sum_{j=1}^{q}(-1)^{j}\Big{[}\rho_X^q(\delta_j\vec{y};x_0)\omega(\delta_j\vec{y};\delta'_0\vec{x})+\sum_{k=1}^{r+1}(-1)^{k}\omega(\delta_j\vec{y};\delta'_k\vec{x})\Big{]}+ \\
									  \quad\qquad +(-1)^{q+1}\rho_Y^{r+1}(y_q;\vec{x})\Big{[}\rho_X^q(\delta_{q+1}\vec{y};x_0)\omega(\delta_{q+1}\vec{y};\delta'_0\vec{x})+\sum_{k=1}^{r+1}(-1)^{k}\omega(\delta_{q+1}\vec{y};\delta'_k\vec{x})\Big{]}. $$
The equality follows from the following identities: Firstly, one 
obviously has 
$$(\delta'_k\vec{x})^{y_0}=\delta'_k(\vec{x})^{y_0}. $$
Secondly, 
$$\rho_Y^r(y_q;\delta'_k\vec{x})=\rho_Y(y_q)^{-1}=\rho_Y^{r+1}(y_q;\vec{x}) $$
and 
$$\rho_X^{q+1}(\vec{y};x_0)=\rho_X(x_0^{y_0...y_q})=\rho_X((x_0^{y_0})^{y_1...y_q})=\rho_X^{q}(\delta_{0}\vec{y};x_0^{y_0}). $$
Also, for all values $0<j\leq q$,
$$\rho_X^{q+1}(\vec{y};x_0)=\rho_X(x_0^{y_0...y_q})=\rho_X(x_0^{y_0...(y_{j-1}y_j)...y_q})=\rho_X^{q}(\delta_{j}\vec{y};x_0). $$
Finally, as $\rho_Y\times\rho_X$ is a double Lie groupoid map
$$\rho_X(x^y)=\rho_X(\Lf{s}(y;x))=\Lf{s}(\rho_Y(y),\rho_X(x))=\rho_Y(y)^{-1}\rho_X(x)\rho_Y(y); $$
thereby implying,
$$\rho_X^{q+1}(\vec{y};x_0)=\rho_X(x_0^{y_0...y_q})=\rho_X((x_0^{y_0...y_{q-1}})^{y_q})=\rho_Y(y_q)^{-1}\rho_X(x_0^{y_0...y_{q-1}})\rho_Y(y_q), $$
which is
$$\rho_X^{q+1}(\vec{y};x_0)\rho_Y^r(y_q;\delta'_0\vec{x})=\rho_Y^{r+1}(y_q;\vec{x})\rho_X^q(\delta_{q+1}\vec{y};x_0) $$
and so, the commutativity of the square (\ref{generic p-square}) follows.

\epf

We use Lemma \ref{p-pageModel} to help proving that the $p$-pages are 
double complexes. For any given $p$, the right action of $\G_p$ on $G$ 
is by automorphisms, as it is defined by the pull-back of the action of 
$H$ along $t_p$, that is, for $g\in G$ and $\gamma\in\G_p$,
$$g^\gamma:=g^{t_p(\gamma)}.$$
Furthermore, the map 
$$(\gamma;g)\to/|->/<350>(\rho_{\G_p}(\gamma),\rho_G(g)):=(\rho_0^1(t_p(\gamma)),\rho_0^1(i(g))) $$
whose components are pull-back representations, verifies the hypothesis 
of Lemma \ref{p-pageModel}. Indeed, after looking at the compatibility 
with the whole structure, one sees that it suffices to prove 
$$\rho_G(\Lf{s}(\gamma;g))=\Lf{s}(\rho_{\G_p}(\gamma),\rho_G(g)). $$
This follows easily from the equivariance of the crossed module map $i$,
$$\rho_G(\Lf{s}(\gamma;g))=\rho_0^1(i(g^{t_p(\gamma)})) 
	 =\rho_0^1(t_p(\gamma)^{-1}i(g)t_p(\gamma)) 
	 =\rho_{\G_p}(\gamma)^{-1}\rho_G(g)\rho_{\G_p}(\gamma). $$
In fact, the $p$-pages coincide by definition with the outcome of Lemma 
\ref{p-pageModel}, but for one caveat; the first column of the $p$-pages 
consists of the $0$th degrees of the complexes in the $r$-direction, 
which are modified (see Lemmas \ref{q=0r-cx} and \ref{r-cx}).

\proposition\label{p-pages}
For each $p$, 
$$\bfig
\node 1a(0,1350)[\vdots]
\node 1b(850,1350)[\vdots]
\node 1c(1800,1350)[\vdots]
\node 1d(2400,1350)[]
\node 2a(0,1000)[C(\G_p^2,V)]
\node 2b(850,1000)[C(\G_p^2\times G,W)]
\node 2c(1800,1000)[C(\G_p^2\times G^2,W)]
\node 2d(2400,1000)[\cdots]
\node 3a(0,500)[C(\G_p,V)]
\node 3b(850,500)[C(\G_p\times G,W)]
\node 3c(1800,500)[C(\G_p\times G^2,W)]
\node 3d(2400,500)[\cdots]
\node 4a(0,0)[V]
\node 4b(850,0)[C(G,W)]
\node 4c(1800,0)[C(G^2,W)]
\node 4d(2400,0)[\cdots]
\arrow|l|[4a`3a;\delta]
\arrow[4a`4b;\delta']
\arrow|l|[3a`2a;\delta]
\arrow[3a`3b;\delta']
\arrow|l|[4b`3b;\delta]
\arrow[4b`4c;\delta_{(1)}]
\arrow[2a`1a;]
\arrow[2a`2b;\delta']
\arrow|l|[3b`2b;\delta]
\arrow[3b`3c;\delta_{(1)}]
\arrow|l|[4c`3c;\delta]
\arrow[4c`4d;]
\arrow[2b`1b;]
\arrow[2b`2c;\delta_{(1)}]
\arrow|l|[3c`2c;\delta]
\arrow[3c`3d;]
\arrow[2c`1c;]
\arrow[2c`2d;]
\efig$$
is a double complex. \endproposition

\pf
Due to Lemmas \ref{q=0r-cx} and \ref{r-cx}, each row is a complex, and 
clearly so is each column. Lemma \ref{p-pageModel} implies that, 
disregarding the first column of squares, we have got a double complex. 
In order to finish the proof, one needs to check that the generic square 
in the first column commutes. This is precisely the content of Lemma 
\ref{1stSquare}, Corollary \ref{1st pSquare} and Lemma \ref{1stpCol}.
					 
\epf

\subsection{Difference maps}\label{DiffMaps} In contrast with Subsection 
\ref{pq-pag}, when $r$ is left constant, the resulting $r$-page fails to 
be a double complex. Nonetheless, as it is briefly mentioned at the 
beginning of this section, the {\em front page}, i.e., the $r$-page for 
$r=0$ commutes up to isomorphism in the $2$-vector space, and, when 
$r>0$, there is a commutation relation up to homotopy. In this 
subsection, we make these comments precise and prove them. 

\subsubsection*{The front page} Let $\omega\in C(\G_p^q,V)$ and 
$\vec{\gamma}\in\G_{p+1}^{q+1}$, then 
$\delta\partial\omega(\vec{\gamma}),\partial\delta\omega(\vec{\gamma})\in V$. 
Recall that in the $2$-vector space $W\to^\phi V$, $v_1,v_2\in V$ belong 
to the same orbit if there exists a $w\in W$ such that $v_2=v_1+\phi(w)$. 
When we say that the front page commutes up to isomorphism, we mean that 
$\delta\partial\omega(\vec{\gamma})$ and $\partial\delta\omega(\vec{\gamma})$ 
belong to the same orbit of $V$. The element in $W$ that realizes the 
isomorphism, is coherently defined to be the composition of the map 
$\delta'$ from the complex in the $r$-direction and certain maps that we 
note $\Delta$ and refer to as the {\em difference maps}.

We start by setting notation and defining the difference maps. Throughout 
this subsection, $\vec{\gamma}\in\G_{p}^{q}$ has components
\begin{equation}\label{gammaMatrix}
\vec{\gamma}=\pmatrix{\gamma_1 \cr \vdots \cr \gamma_q}=\pmatrix{
    \gamma_{11} & \gamma_{12} & ... & \gamma_{1p} \cr
    \gamma_{21} & \gamma_{22} & ... & \gamma_{2p} \cr
   \vdots  & \vdots &     & \vdots \cr
    \gamma_{q1} & \gamma_{q2} & ... & \gamma_{qp}}=\pmatrix{
    g_{11} & g_{12} & ... & g_{1p} & h_1 \cr
    g_{21} & g_{22} & ... & g_{2p} & h_2\cr
   \vdots  & \vdots &     & \vdots & \vdots \cr
    g_{q1} & g_{q2} & ... & g_{qp} & h_q},
\end{equation}
where the last equality is a notation abuse corresponding to the 
row-wise isomorphism $\G_{p}\cong G^{p}\times H$ from Remark 
\ref{2GpNerve}. Here, $(g_{ab},h_{ab})$ is the image of $\gamma_{ab}$ 
under the isomorphism of Remark \ref{Equiv2Gp:Mods} for all values of 
$a$ and $b$. We abbreviate $\partial_0\delta_0\vec{\gamma}$ by regarding 
$\vec{\gamma}$ as a matrix and using ``minor'' notation, {\it i.e.},
$$\vec{\gamma}_{1,1}:=\partial_0\delta_0\vec{\gamma}=\pmatrix{
    \gamma_{22} & \gamma_{23} & ... & \gamma_{2q} \cr
    \gamma_{32} & \gamma_{33} & ... & \gamma_{3q} \cr
   \vdots  & \vdots &     & \vdots \cr
    \gamma_{p2} & \gamma_{p3} & ... & \gamma_{pq}}. $$
    
The difference maps
$$\Delta:C(\G_p^q\times G,W)\to<350> C(\G_{p+1}^{q+1},V)$$ 
are defined by
$$\Delta\omega(\vec{\gamma}):=\rho_0^1(t_p(\partial_0\gamma_1)...t_p(\partial_0\gamma_{q+1}))\circ\phi\big{(}\omega(\vec{\gamma}_{1,1};g_{11})\big{)}, $$
for $\omega\in C(\G_p^q\times G,W)$ and $\vec{\gamma}\in\G_{p+1}^{q+1}$. 
Clearly, when $q=0$, one drops the $\vec{\gamma}_{1,1}$-entry (compare 
with Eq.'s (\ref{1.Delta}), (\ref{2.Delta}) and (\ref{3.Delta})).

\lem\label{FrontPage-upToIsoq=0}
Let $v\in C^{p,0}_0(\G,\phi)=V$, then 
\begin{equation}\label{upToIsoq=0}
\delta\partial v=\partial\delta v+\Delta\delta' v \in C(\G_{p+1},V).
\end{equation} \endlem

\pf
First, we compute the first term on the right hand side of Eq. 
(\ref{upToIsoq=0}). Let $\gamma\in\G_{p+1}$ and let 
$(g_0,...,g_p;h)\in G^{p+1}\times H$ be its corresponding image under 
the isomorphism of Remark \ref{2GpNerve}. Then,
$$\mld (\partial\delta v)(\gamma) & =(\delta v)(\partial_0\gamma) 
	+\sum_{j=1}^{p}(-1)^{j+1}(\delta v)(\partial_j\gamma)+(-1)^{p+1}(\delta v)(\partial_{p+1}\gamma) \\
		 =\rho_0^0(hi(g_p...g_1))v-v+\sum_{j=1}^{p}(-1)^{j+1}(\rho_0^0(hi(g_p...g_0))v-v)+ \\
		 \qquad\quad +(-1)^{p+1}(\rho_0^0(hi(g_p)i(g_{p-1}...g_0))v-v), $$
and 
$$(\partial\delta v)(\gamma)=\cases{\rho_0^0(hi(g_p...g_1))v-v & if $p$ is odd \cr 
				  \rho_0^0(hi(g_p...g_1))(v-\rho_0^0(i(g_0))v) & otherwise. }$$
On the other hand, 
$$(\delta\partial v)(\gamma)=\cases{\rho_0^0(hi(g_p...g_0))v-v & if $p$ is odd \cr
														   0 & otherwise; }$$ 
hence, either way, 
$$(\delta\partial v-\partial\delta v)(\gamma)=\rho_0^0(hi(g_p...g_1))(\rho_0^0(i(g_0))v-v) $$
and the result follows from the first part of Eq. (\ref{2}). 
 
\epf

\proposition\label{FrontPage-upToIso}
Let $\omega\in C(\G_p^q,V)$, then 
\begin{equation}\label{upToIso}
\delta\partial\omega=\partial\delta\omega+\Delta\delta'\omega \in C(\G_{p+1}^{q+1},V).
\end{equation} \endproposition

\pf
Eq. (\ref{upToIso}) holds essentially due to the commutativity of 
$\delta_j$ and $\partial_k$ for all values of $(j,k)$, one just need to 
take care of the representations that appear at $(0,0)$. Let 
$\vec{\gamma}\in\G_{p+1}^{q+1}$, then
$$\mld \delta\partial\omega(\vec{\gamma}) & = \rho_0^0(t_{p+1}(\gamma_1))\partial\omega(\delta_0\vec{\gamma})+\sum_{j=1}^{q+1}(-1)^j\partial\omega(\delta_j\vec{\gamma}) \\
		= \rho_0^0(t_{p+1}(\gamma_1))\sum_{k=0}^{p+1}(-1)^k\omega(\partial_k\delta_0\vec{\gamma})+\sum_{j=1}^{q+1}\sum_{k=0}^{p+1}(-1)^{j+k}\omega(\partial_k\delta_j\vec{\gamma}). $$
On the other hand,
$$\mld \partial\delta\omega(\vec{\gamma}) & = \sum_{k=0}^{p+1}(-1)^k\delta\omega(\partial_k\vec{\gamma}) \\
	    = \sum_{k=0}^{p+1}(-1)^k\Big{(}\rho_0^0(t_p((\partial_k\vec{\gamma})_1))\omega(\delta_0\partial_k\vec{\gamma})+\sum_{j=1}^{q+1}(-1)^{j}\omega(\delta_j\partial_k\vec{\gamma})\Big{)}.  $$
As stated, the double sums in the above expressions coincide. 
By definition, $(\partial_k\vec{\gamma})_1=\partial_k\gamma_1$. 
In accordance with Remark \ref{2GpNerve}, this corresponds to 
$(g_{11},...,g_{1k-2},g_{1k}g_{1k-1},g_{1k+1},...,g_{1p+1},h_1)$, 
when $0\lt k\leq p$, and to $(g_{11},...,g_{1p},h_0i(g_{0p+1}))$, 
when $k=p+1$. Hence, for all $0\lt k\leq p+1$,
$$t_p((\partial_k\vec{\gamma})_1)=h_1i(g_{1p+1}...g_{11})=t_{p+1}(\gamma_1). $$
Thus, using the first part of Eq. (\ref{2}), one computes
$$\mld (\delta\partial\omega-\partial\delta\omega)(\vec{\gamma}) & =\rho_0^0(t_{p+1}(\gamma_1))\omega(\vec{\gamma}_{1,1})-\rho_0^0(t_p(\partial_0\gamma_1))\omega(\vec{\gamma}_{1,1}) \\
	    =\rho_0^0(h_1i(g_{1p+1}...g_{12}))\Big{[}\rho_0^0(i(g_{11}))-I\Big{]}\omega(\vec{\gamma}_{1,1}) 
		=\rho_0^0(t_p(\partial_0\gamma_1))\circ\phi\Big{(}\rho_1(g_{11})\omega(\vec{\gamma}_{1,1})\Big{)}. $$
													
\epf

\rem\label{upToIsoIndeed}
As claimed, Lemma \ref{FrontPage-upToIsoq=0} and Proposition 
\ref{FrontPage-upToIso} are interpreted as saying that, if 
$\omega\in C(\G_p^q,V)$ and $\vec{\gamma}\in\G_{p+1}^{q+1}$, 
$\delta\partial\omega(\vec{\gamma})$ and 
$\partial\delta\omega(\vec{\gamma})$ lie on the same orbit of $V$. 
Indeed, their difference lies in the image of the difference map 
$\Delta$, which, using Eq. (\ref{0}), lies in the image of the structural 
map $\phi$ of the $2$-vector space. \endrem

\subsubsection*{$r$-pages} If $\omega\in C(\G_p^q\times G^r,W)$ and 
$(\vec{\gamma};\vec{f})\in\G_{p+1}^{q+1}\times G^r$, then 
$\delta\partial\omega(\vec{\gamma};\vec{f}),\partial\delta\omega(\vec{\gamma};\vec{f})\in W$. 
One cannot expect results analogous to Lemma \ref{FrontPage-upToIsoq=0} 
and Proposition \ref{FrontPage-upToIso}, because there are no orbits in 
$W$. We prove instead that the compositions $\delta\circ\partial$ and 
$\partial\circ\delta$ are homotopic when regarded as maps between 
$r$-complexes. In what seems an overlap of notation, we call the 
homotopies $\Delta$ and refer to them as difference maps too. 

We need to introduce further notation conventions to define the difference 
maps. Let $\vec{f}=(f_1,...,f_r)\in G^r$, then for any pair of integers 
$1\leq a<b\leq r$, define
$$\vec{f}_{[a,b]}:=(f_a,f_{a+1},...,f_{b-1},f_b)  \qquad\textnormal{and}\qquad\  \vec{f}_{[a,b)}:=(f_a,f_{a+1},...,f_{b-2},f_{b-1}). $$
With this shorthand, for $r>1$ and $0\lt n\lt r$, we define 
$$c_{2n-1},c_{2n}:G^{r-1}\times \G \to<350> G^r $$
by
\begin{equation}\label{c(11)2n-1}
c_{2n-1}(\vec{f};\gamma):=\Big{(}\big{(}\vec{f}_{[1,r-n)}\big{)}^{hi(g)},g^{-1},\big{(}\vec{f}_{[r-n,r-2]}\big{)}^h,f_{r-1}^hg\Big{)}
\end{equation}
and
\begin{equation}\label{c(11)2n}
c_{2n}(\vec{f};\gamma):=\Big{(}\big{(}\vec{f}_{[1,r-n)}\big{)}^{hi(g)},g^{-1},\big{(}\vec{f}_{[r-n,r-2]}\big{)}^h,g\Big{)},
\end{equation}
where $\vec{f}=(f_1,...,f_{r-1})\in G^{r-1}$, $\gamma\in\G$ and 
$(g,h)\in G\rtimes H$ is the image of $\gamma$ under the isomorphism of 
Remark \ref{Equiv2Gp:Mods}.
 
Let $p\geq 0$ and $q=0$, then the difference maps  
$$\Delta:C(G^{r},W)\to<350> C(\G_{p+1}\times G^{r-1},W)$$ 
are defined by
$$\Delta\omega(\gamma;\vec{f}) = \rho_0^1(i(g_0))^{-1}\omega((\vec{f})^{h_0},g_0)+ 
  \sum_{n=1}^{r-1}(-1)^{n+1}\big{(}\omega(c_{2n-1}(\vec{f};\gamma_0))-\omega(c_{2n}(\vec{f};\gamma_0))\big{)}, $$
for $\omega\in C(G^{r},W)$, $\vec{f}=(f_1,...,f_{r-1})\in G^{r-1}$ and 
$\gamma=(\gamma_0,...,\gamma_p)\in\G_{p+1}$. Here, 
$(g_0,h_0)\in G\rtimes H$ is the image of $\gamma_0$ under the 
isomorphism of Remark \ref{Equiv2Gp:Mods}. 

When $q>0$, the difference maps
$$\Delta:C(\G_p^q\times G^{r},W)\to<350> C(\G_{p+1}^{q+1}\times G^{r-1},W)$$ 
are defined by
$$\mld \Delta\omega(\vec{\gamma};\vec{f}) & =\rho_0^1(i(pr_G(\gamma_{21}\vJoin ...\vJoin\gamma_{(q+1)1})))^{-1}\Big{[}\rho_0^1(i(g_{11}^{h_{21}...h_{(q+1)1}}))^{-1}\omega(\vec{\gamma}_{1,1};(\vec{f})^{h_{11}},g_{11})+ \\
 \qquad +\sum_{n=1}^{r-1}(-1)^{n+1}\big{(}\omega(\vec{\gamma}_{1,1};c_{2n-1}(\vec{f};\gamma_{11}))-\omega(\vec{\gamma}_{1,1};c_{2n}(\vec{f};\gamma_{11}))\big{)}\Big{]}, $$
for $\omega\in C(\G_p^q\times G^{r+1},W)$, $\vec{f}\in G^{r-1}$ and 
$\vec{\gamma}\in\G_{p+1}^{q+1}$ is as in Eq. (\ref{gammaMatrix}).

Lemma \ref{q=0r=1-page} and Proposition \ref{r=1-page} below justify the 
choice of notation.

\lemma\label{q=0r=1-page}
Let $\omega\in C^{p,0}_1(\G,\phi)=C(G,W)$, then 
\begin{equation}\label{q=0r=1UpToHom}
(\partial\delta-\delta\partial)\omega=(\Delta\delta_{(1)}-\delta'\Delta)\omega\in C(\G_{p+1}\times G,W).
\end{equation}\endlem

\pf
Let $f\in G$, $\gamma\in\G_{p+1}$ and 
$(g_0,...,g_p;h)\in G^{p+1}\times H$ its corresponding image under the 
isomorphism of Remark \ref{2GpNerve}. Then,
$$\delta\partial\omega(\gamma;f)=\cases{\omega(f^{t_{p+1}(\gamma)})-\rho_0^1(t_{p+1}(\gamma))^{-1}\omega(f) & if $p$ is odd \cr
0 & otherwise} $$
and 
$$\mld \partial\delta\omega(\gamma;f) & =\rho_0^1(i(g_0))^{-1}\delta\omega(\partial_0\gamma;f)+\sum_{j=1}^{p+1}(-1)^j\delta\omega(\partial_j\gamma;f) \\
	=\rho_0^1(i(g_0))^{-1}\Big{[}\omega(f^{t_p(\partial_0\gamma)})-\rho_0^1(t_p(\partial_0\gamma))^{-1}\omega(f)\Big{]}+ \\
 \qquad\qquad +\sum_{j=1}^{p+1}(-1)^j\Big{[}\omega(f^{t_p(\partial_j\gamma)})-\rho_0^1(t_p(\partial_j\gamma))^{-1}\omega(f)\Big{]}. $$
As in the proof of Proposition \ref{FrontPage-upToIso}, 
$t_p(\partial_0\gamma)=hi(g_p...g_1)$ and, for $0\lt j\leq p+1$, 
$t_p(\partial_j\gamma)=hi(g_p...g_0)=t_{p+1}(\gamma)$; therefore,
$$\partial\delta\omega(\gamma;f)=\cases{\rho_0^1(i(g_0))^{-1}\omega(f^{t_p(\partial_0\gamma)})-\rho_0^1(t_{p+1}(\gamma))^{-1}\omega(f) & if $p$ is odd \cr
\rho_0^1(i(g_0))^{-1}\omega(f^{t_p(\partial_0\gamma)})-\omega(f^{t_{p+1}(\gamma)}) & otherwise, } $$
and, either way,
$$(\partial\delta-\delta\partial)\omega(\gamma;f)=\rho_0^1(i(g_0))^{-1}\omega(f^{t_p(\partial_0\gamma)})-\omega(f^{t_{p+1}(\gamma)}). $$
On the other hand,
$$\mld \delta'\Delta\omega(\gamma;f) & =\rho_0^1(t_{p+1}(\gamma))^{-1}\rho_1(f)\Delta\omega(\gamma) \\
	=\rho_0^1(t_p(\partial_0\gamma)i(g_0))^{-1}\rho_1(f)\rho_0^0(t_p(\partial_0\gamma))\phi\big{(}\omega(g_0)\big{)}=\rho_0^1(i(g_0))^{-1}\rho_1(f^{t_p(\partial_0\gamma)})\phi\big{(}\omega(g_0)\big{)} $$
and 
$$\mld \Delta\delta_{(1)}\omega(\gamma;f) & =\rho_0^1(i(g_0))^{-1}\delta_{(1)}\omega(f^{h_0},g_0)+\delta_{(1)}\omega(g_0^{-1},f^{h_0}g_0)-\delta_{(1)}\omega(g_0^{-1},g_0) \\
	=\rho_0^1(i(g_0))^{-1}\Big{[}\rho_0^1(i(f^{h_0}))\omega(g_0)-\omega(f^{h_0}g_0)+\omega(f^{h_0})\Big{]}+ \\
	\qquad + \Big{[}\rho_0^1(i(g_0))^{-1}\omega(f^{h_0}g_0)-\omega(g_0^{-1}f^{h_0}g_0)+\omega(g_0^{-1})\Big{]} \\
	\qquad\qquad -\Big{[}\rho_0^1(i(g_0))^{-1}\omega(g_0)-\cancel{\omega(g_0^{-1}g_0)}+\omega(g_0^{-1})\Big{]} \\
	=\rho_0^1(i(g_0))^{-1}\Big{[}\rho_0^1(i(f^{h_0}))\omega(g_0)+\omega(f^{h_0})-\omega(g_0)\Big{]}-\omega(f^{h_0i(g_0)}) \\
	=\rho_0^1(i(g_0))^{-1}\Big{[}\rho_1(f^{h_0})\phi\big{(}\omega(g_0)\big{)}+\omega(f^{h_0})\Big{]}-\omega(f^{h_0i(g_0)}), $$
where the last equality follows from the second half of Eq. (\ref{2}). 
In so, given that $t_p(\partial_0\gamma)=h_0$ and 
$t_{p+1}(\gamma)=h_0i(g_0)$, one computes the difference to be
$$(\Delta\delta_{(1)}-\delta'\Delta)\omega(\gamma;f)=\rho_0^1(i(g_0))^{-1}\omega(f^{t_p(\partial_0\gamma)})-\omega(f^{t_{p+1}(\gamma)}), $$
and the result follows.

\epf

\proposition\label{r=1-page}
Let $\omega\in C(\G_p^q\times G,W)$, then
\begin{equation}\label{r=1UpToHom}
(\partial\delta-\delta\partial)\omega=(\Delta\delta_{(1)}-\delta'\Delta)\omega\in C(\G_{p+1}^{q+1}\times G,W).
\end{equation}\endproposition

\pf
Let $\vec{\gamma}\in\G_{p+1}^{q+1}$ be as in Eq. (\ref{gammaMatrix}) and 
$f\in G$. To ease up notation, we introduce the shorthand
\begin{equation}\label{p-repShorthand}
\rho^{q}(\vec{\gamma}_{\bullet 1}):=\rho_0^1(i(pr_G(\gamma_{11}\vJoin...\vJoin\gamma_{q1})))^{-1}.
\end{equation} 
We compute the left hand side of Eq. (\ref{r=1UpToHom}). On the one hand,
$$\mld \partial\delta\omega(\vec{\gamma};f) & =\rho^{q+1}(\vec{\gamma}_{\bullet 1})\delta\omega(\partial_0\vec{\gamma};f)+\sum_{j=1}^{p+1}(-1)^j\delta\omega(\partial_j\vec{\gamma};f) \\
         =\rho^{q+1}(\vec{\gamma}_{\bullet 1})\Big{(}\omega(\vec{\gamma}_{1,1};f^{t_p(\partial_0\gamma_1)})+\sum_{k=1}^q(-1)^k\omega(\delta_k\partial_0\vec{\gamma};f)+ \\
         \qquad\qquad\qquad\qquad\qquad +(-1)^{q+1}\rho_0^1(t_p(\partial_0\gamma_{q+1}))^{-1}\omega(\delta_{q+1}\partial_0\vec{\gamma};f)\Big{)}+ \\
         \quad +\sum_{j=1}^{p+1}(-1)^j\Big{(}\omega(\delta_0\partial_j\vec{\gamma};f^{t_p(\partial_j\gamma_1)})+\sum_{k=1}^q(-1)^k\omega(\delta_k\partial_j\vec{\gamma};f)+ \\
         \qquad\qquad\qquad\qquad\qquad +(-1)^{q+1}\rho_0^1(t_p(\partial_j\gamma_{q+1}))^{-1}\omega(\delta_{q+1}\partial_j\vec{\gamma};f)\Big{)}; $$
while on the other,
$$\mld \delta\partial\omega(\vec{\gamma};f) & =\partial\omega(\delta_0\vec{\gamma};f^{t_{p+1}(\gamma_1)})+\sum_{k=1}^q(-1)^k\partial\omega(\delta_k\vec{\gamma};f)+(-1)^{q+1}\rho_0^1(t_{p+1}(\gamma_{q+1}))^{-1}\omega(\delta_{q+1}\vec{\gamma};f) \\
       =\rho^{q}((\delta_0\vec{\gamma})_{\bullet 1})\omega(\vec{\gamma}_{1,1};f^{t_{p+1}(\gamma_1)})+\sum_{j=1}^{p+1}(-1)^j\omega(\partial_j\delta_0\vec{\gamma};f^{t_{p+1}(\gamma_1)})+ \\
       \qquad +\sum_{k=1}^q(-1)^k\Big{(}\rho^{q}((\delta_k\vec{\gamma})_{\bullet 1})\omega(\partial_0\delta_k\vec{\gamma};f)+\sum_{j=1}^{p+1}(-1)^j\omega(\partial_j\delta_k\vec{\gamma};f)\Big{)}+ $$
$$ +(-1)^{q+1}\rho_0^1(t_{p+1}(\gamma_{q+1}))^{-1}\Big{(}\rho^{q}((\delta_{q+1}\vec{\gamma})_{\bullet 1})\omega(\partial_0\delta_{q+1}\vec{\gamma};f)+\sum_{j=1}^{p+1}(-1)^j\omega(\partial_j\delta_{q+1}\vec{\gamma};f)\Big{)}. $$
We claim that taking the difference yields
$$ (\partial\delta-\delta\partial)\omega(\vec{\gamma};f)=\rho^{q+1}(\vec{\gamma}_{\bullet 1})\omega(\vec{\gamma}_{1,1};f^{t_p(\partial_0\gamma_1)})-\rho^{q}((\delta_0\vec{\gamma})_{\bullet 1})\omega(\vec{\gamma}_{1,1};f^{t_{p+1}(\gamma_1)}). $$
This follows from the commutativity of all simplicial maps 
$\delta_k\partial_j=\partial_j\delta_k$, together with the following 
identities:
\begin{itemize}
  \item $\rho^{q+1}(\vec{\gamma}_{\bullet 1})=\rho^{q}((\delta_k\vec{\gamma})_{\bullet 1})$ 
  for $1\leq k\leq q$. Indeed, for the ranging values of $k$, 
    $$\rho^{q}((\delta_k\vec{\gamma})_{\bullet 1})=\rho_0^1(i(pr_G(\gamma_{11}\vJoin...\vJoin(\gamma_{(k-1)1}\vJoin\gamma_{k1})\vJoin...\vJoin\gamma_{(q+1)1})))^{-1}=\rho^{q+1}(\vec{\gamma}_{\bullet 1}). $$
  \item $\rho^{q+1}(\vec{\gamma}_{\bullet 1})\rho_0^1(t_p(\partial_0\gamma_{q+1}))^{-1}=\rho_0^1(t_{p+1}(\gamma_{q+1}))^{-1}\rho^{q}((\delta_{q+1}\vec{\gamma})_{\bullet 1})$. Indeed, using Lemma 
  \ref{multiprods}, we write 
    $$\mld \rho^{q+1}(\vec{\gamma}_{\bullet 1}) & =\rho_0^1(i(g_{11}^{h_{21}...h_{(q+1)1}}g_{21}^{h_{31}...h_{(q+1)1}}...g_{(q-1)1}^{h_{q1}h_{(q+1)1}}g_{q1}^{h_{(q+1)1}}g_{(q+1)1}))^{-1} \\
    		=\rho_0^1(i(g_{(q+1)1}))^{-1}\rho_0^1\big{(}i\big{(}(g_{11}^{h_{21}...h_{q1}}g_{21}^{h_{31}...h_{q1}}...g_{(q-1)1}^{h_{q1}}g_{q1})^{h_{(q+1)1}}\big{)}\big{)}^{-1} \\
    		=\rho_0^1(h_{(q+1)1}i(g_{(q+1)1}))^{-1}\rho_0^1\big{(}i(g_{11}^{h_{21}...h_{q1}}g_{21}^{h_{31}...h_{q1}}...g_{(q-1)1}^{h_{q1}}g_{q1})\big{)}^{-1}\rho_0^1(h_{(q+1)1}),  $$
	and by definition $t_p(\partial_0\gamma_{q+1})=t(\gamma_{(q+1)2})=s(\gamma_{(q+1)1})=h_{(q+1)1}$.
    \item For $1\leq j\leq q+1$, 
    $t_p(\partial_j\gamma_b)=t_{p+1}(\gamma_b)$. 
\end{itemize}
Using Lemma \ref{multiprods} once more,   
$$\mld \rho^{q+1}(\vec{\gamma}_{\bullet 1}) & =\rho_0^1(i(g_{21}^{h_{31}...h_{(q+1)1}}...g_{(q-1)1}^{h_{q1}h_{(q+1)1}}g_{q1}^{h_{(q+1)1}}g_{(q+1)1}))^{-1}\rho_0^1(i(g_{11}^{h_{21}...h_{(q+1)1}}))^{-1} \\
	   =\rho^q((\delta_0\vec{\gamma})_{\bullet 1})\rho_0^1(i(g_{11}^{h_{21}...h_{(q+1)1}}))^{-1}, $$
and we rewrite the difference as
$$(\partial\delta-\delta\partial)\omega(\vec{\gamma};f)=\rho^q((\delta_0\vec{\gamma})_{\bullet 1})\big{[}\rho_0^1(i(g_{11}^{h_{21}...h_{(q+1)1}}))^{-1}\omega(\vec{\gamma}_{1,1};f^{h_{11}})-\omega(\vec{\gamma}_{1,1};f^{h_{11}i(g_{11})})\big{]}. $$
We turn to compute the right hand side of Eq. (\ref{r=1UpToHom}),
$$\mld \delta'\Delta\omega(\vec{\gamma};f) & =\rho_0^1(t_{p+1}(\gamma_1)...t_{p+1}(\gamma_{q+1}))^{-1}\rho_1(f)\Delta\omega(\vec{\gamma}) \\
       =\rho_0^1(t(\gamma_{11}\vJoin...\vJoin\gamma_{(q+1)1}))^{-1}\rho_1(f)\rho_0^0(t_p(\partial_0\gamma_1)...t_p(\partial_0\gamma_{q+1}))\phi(\omega(\vec{\gamma}_{1,1};g_{11})) \\
       =\rho_0^1(h_{11}...h_{(q+1)1}i(pr_G(\gamma_{11}\vJoin...\vJoin\gamma_{(q+1)1})))^{-1}\rho_1(f)\rho_0^0(h_{11}...h_{(q+1)1})\phi(\omega(\vec{\gamma}_{1,1};g_{11})) \\
       =\rho^{q+1}(\vec{\gamma}_{\bullet 1})\rho_1(f^{h_{11}...h_{(q+1)1}})\phi(\omega(\vec{\gamma}_{1,1};g_{11})).  $$
Adding all terms, factoring $\rho^{q+1}(\vec{\gamma}_{\bullet 1})$ as 
before and using the second half of Eq. (\ref{2}) yields 
$$\mld (\delta\partial-\partial\delta-\delta'\Delta & )\omega(\vec{\gamma};f) 
      =\rho^q((\delta_0\vec{\gamma})_{\bullet 1})\Big{[}\omega(\vec{\gamma}_{1,1};f^{h_{11}i(g_{11})})+ \\
       -\rho_0^1(i(g_{11}^{h_{21}...h_{(q+1)1}}))^{-1}\big{(}\omega(\vec{\gamma}_{1,1};f^{h_{11}})
       +[\rho_0^1(i(f^{h_{11}...h_{(q+1)1}}))-I]\omega(\vec{\gamma}_{1,1};g_{11})\big{)}\Big{]}. $$
Adding and subtracting 
$\rho_0^1(i(g_{11}^{h_{21}...h_{(q+1)1}}))^{-1}\omega(\vec{\gamma}_{1,1};f^{h_{11}}g_{11})$ 
and $\omega(\vec{\gamma}_{1,1};g_{11}^{-1})$,
$$\mld (\delta\partial-\partial\delta-\delta'\Delta)\omega(\vec{\gamma};f)= &  
      \rho^q((\delta_0\vec{\gamma})_{\bullet 1})\Big{[}\delta_{(1)}\omega(\vec{\gamma}_{1,1};g_{11}^{-1},g_{11})
      -\delta_{(1)}\omega(\vec{\gamma}_{1,1};g_{11}^{-1},f^{h_{11}}g_{11})+ \\
      -\rho_0^1(i(g_{11}^{h_{21}...h_{(q+1)1}}))^{-1}\delta_{(1)}\omega(\vec{\gamma}_{1,1};f^{h_{11}},g_{11})\Big{]} 
       =-\Delta\delta_{(1)}\omega(\vec{\gamma};f). $$
\epf

\proposition\label{q=0r-page}
Let $r>1$ and $\omega\in C^{p,0}_r(\G,\phi)=C(G^r,W)$, then
\begin{equation}\label{q=0r-UpToHom}
(-1)^r(\delta\partial-\partial\delta)\omega=(\Delta\delta_{(1)}-\delta_{(1)}\Delta)\omega\in C(\G_{p+1}\times G^r,W).	
\end{equation} \endproposition

\pf
We compute the left hand side of Eq. (\ref{q=0r-UpToHom}). Let 
$\vec{f}=(f_1,...,f_r)\in G^r$, 
$\gamma=(\gamma_0,...,\gamma_p)\in\G_{p+1}$ and $(g_0,h_0)\in G\rtimes H$ 
be the image of $\gamma_0$ under the isomorphism of Remark 
\ref{Equiv2Gp:Mods}. Then, 
$$\delta\partial\omega(\gamma;\vec{f})=\cases{\omega\big{(}(\vec{f})^{h_0i(g_0)}\big{)}-\rho_0^1(h_0i(g_0))^{-1}\omega(\vec{f})                & if $p$ is odd \cr
 0 & otherwise },$$
and
$$\mld \partial\delta\omega(\gamma;\vec{f}) & =\rho_0^1(i(g_0))^{-1}\delta\omega(\partial_0\gamma;\vec{f})+\sum_{j=1}^{p+1}(-1)^j\delta\omega(\partial_j\gamma;\vec{f}) \\
     =\rho_0^1(i(g_0))^{-1}\Big{[}\omega\big{(}(\vec{f})^{h_0}\big{)}-\rho_0^1(h_0)^{-1}\omega(\vec{f})\Big{]}+ \\
     \qquad\qquad +\sum_{j=1}^{p+1}(-1)^j\Big{[}\omega\big{(}(\vec{f})^{h_0i(g_0)}\big{)}-\rho_0^1(h_0i(g_0))^{-1}\omega(\vec{f})\Big{]}. $$
As in the proof of Lemma \ref{q=0r=1-page}, one concludes 
$$\partial\delta\omega(\gamma;\vec{f})=\cases{\rho_0^1(i(g_0))^{-1}\Big{[}\omega\big{(}(\vec{f})^{h_0}\big{)}-\rho_0^1(h_0)^{-1}\omega(\vec{f})\Big{]} & if $p$ is odd \cr
 \rho_0^1(i(g_0))^{-1}\omega\big{(}(\vec{f})^{h_0}\big{)}-\omega\big{(}(\vec{f})^{h_0i(g_0)}\big{)} & otherwise}, $$ 
and in both cases
\begin{equation}\label{theLHS r-pag}
(\delta\partial-\partial\delta)\omega(\gamma;\vec{f})=\omega\big{(}(\vec{f})^{h_0i(g_0)}\big{)}-\rho_0^1(i(g_0))^{-1}\omega\big{(}(\vec{f})^{h_0}\big{)}.
\end{equation}

Turning to the left hand side of Eq. (\ref{q=0r-UpToHom}), on the one 
hand, 
$$\Delta\delta_{(1)}\omega(\gamma;\vec{f})=\rho_0^1(i(g_0))^{-1}\delta_{(1)}\omega((\vec{f})^{h_0},g_0)+\sum_{n=1}^r(-1)^{n+1}T_n, $$ 
where 
$$T_n:=\delta_{(1)}\omega(c_{2n-1}(\vec{f};\gamma_0))-\delta_{(1)}\omega(c_{2n}(\vec{f};\gamma_0)), $$
and, on the other,
\begin{equation}\label{RHS2 r-pag}
\delta_{(1)}\Delta\omega(\gamma;\vec{f})=\rho_0^1(i(f_1^{h_0i(g_0)}))\Delta\omega(\gamma;\delta_0\vec{f})+\sum_{k=1}^r(-1)^k\Delta\omega(\gamma;\delta_k\vec{f}).
\end{equation}
Rearranging Eq. (\ref{RHS2 r-pag}),
$$\delta_{(1)}\Delta\omega(\gamma;\vec{f})=S_0+\sum_{n=1}^{r-1}(-1)^{n+1}S_n, $$
where
$$S_0:=\rho_0^1(i(f_1^{h_0i(g_0)}))\rho_0^1(i(g_0))^{-1}\omega((\delta_0\vec{f})^{h_0},g_0)+\sum_{k=1}^r(-1)^k\rho_0^1(i(g_0))^{-1}\omega((\delta_k\vec{f})^{h_0},g_0) $$
and
$$\mld S_n:=\rho_0^1(i(f_1^{h_0i(g_0)})) & \big{(}\omega(c_{2n-1}(\delta_0\vec{f};\gamma_0))-\omega(c_{2n}(\delta_0\vec{f};\gamma_0)\big{)}+ \\
		+\sum_{k=1}^r(-1)^k\big{(}\omega(c_{2n-1}(\delta_k\vec{f};\gamma_0))-\omega(c_{2n}(\delta_k\vec{f};\gamma_0)\big{)}. $$
Expanding further the first term of $\Delta\delta_{(1)}\omega$, 
$$\mld \rho_0^1(i(g_0))^{-1} & \delta_{(1)}\omega((\vec{f})^{h_0},g_0)=\rho_0^1(i(g_0))^{-1}\Big{(}\rho_0^1(i(f_1^{h_0}))\omega(\delta_0(\vec{f})^{h_0},g_0)+ \\
	+\sum_{k=1}^{r-1}(-1)^k\omega(\delta_k(\vec{f})^{h_0},g_0)+(-1)^r\omega(\delta_r(\vec{f})^{h_0},f_r^{h_0}g_0)+(-1)^{r+1}\omega((\vec{f})^{h_0})\Big{)}, $$
it is made patent that, if 
$\epsilon_0:=\rho_0^1(i(g_0))^{-1}\delta_{(1)}\omega((\vec{f})^{h_0},g_0)-S_0$, 
then
$$\epsilon_0=(-1)^r\rho_0^1(i(g_0))^{-1}\Big{(}\omega(\delta_r(\vec{f})^{h_0},f_r^{h_0}g_0)-\omega((\vec{f})^{h_0})+\omega((\delta_r\vec{f})^{h_0},g_0)\Big{)}, $$ 
as $(\delta_k\vec{f})^{h_0}=\delta_k(\vec{f})^{h_0}$ and 
$g_0^{-1}f_1^{h_0}=f_1^{h_0i(g_0)}g_0^{-1}$. 

Since one can equivalently write
$$c_{2r-1}(\vec{f};\gamma_0)=(g_0^{-1},(\delta_r\vec{f})^{h_0},f_r^{h_0}g_0)\quad\textnormal{ and }\quad 
c_{2r}(\vec{f};\gamma_0)=(g_0^{-1},(\delta_r\vec{f})^{h_0},g_0), $$ 
$T_r$ gets expanded as
$$\mld T_r= & \rho_0^1(i(g_0))^{-1}\big{(}\omega((\delta_r\vec{f})^{h_0},f_r^{h_0}g_0)-\omega((\delta_r\vec{f})^{h_0},g_0)\big{)}+ \\
		-\big{(}\omega(g_0^{-1}f_1^{h_0},\delta_0(\delta_r\vec{f})^{h_0},f_r^{h_0}g_0)-\omega(g_0^{-1}f_1^{h_0},\delta_0(\delta_r\vec{f})^{h_0},g_0)\big{)}+ \\
	\quad +\sum_{k=1}^{r-2}(-1)^{k+1}\big{(}\omega(g_0^{-1},\delta_k(\delta_r\vec{f})^{h_0},f_r^{h_0}g_0)-\omega(g_0^{-1},\delta_k(\delta_r\vec{f})^{h_0},g_0)\big{)}+ \\
	\qquad +(-1)^{r}\big{(}\omega(g_0^{-1},\delta_{r-1}(\delta_r\vec{f})^{h_0},(f_{r-1}f_r)^{h_0}g_0)-\omega(g_0^{-1},\delta_{r-1}(\delta_r\vec{f})^{h_0},f_{r-1}^{h_0}g_0)\big{)}+ \\
	\quad\qquad +(-1)^{r+1}\big{(}\cancel{\omega(g_0^{-1},(\delta_r\vec{f})^{h_0})}-\cancel{\omega(g_0^{-1},(\delta_r\vec{f})^{h_0})}\big{)}. $$ 
On the other hand, 
$$\mld S_{r-1}= & \rho_0^1(i(f_1^{h_0i(g_0)}))^{-1}\big{(}\omega(g_0^{-1},\delta_{r-1}(\delta_0\vec{f})^{h_0},f_r^{h_0}g_0)-\omega(g_0^{-1},\delta_{r-1}(\delta_0\vec{f})^{h_0},g_0)\big{)}+ \\
 	+\sum_{k=1}^{r-2}(-1)^{k}\big{(}\omega(g_0^{-1},\delta_{r-1}(\delta_k\vec{f})^{h_0},f_r^{h_0}g_0)-\omega(g_0^{-1},\delta_{r-1}(\delta_k\vec{f})^{h_0},g_0)\big{)}+ \\
 	\quad +(-1)^{r-1}\big{(}\omega(g_0^{-1},\delta_{r-1}(\delta_{r-1}\vec{f})^{h_0},(f_{r-1}f_r)^{h_0}g_0)-\cancel{\omega(g_0^{-1},\delta_{r-1}(\delta_{r-1}\vec{f})^{h_0},g_0)}\big{)} \\
 	\qquad +(-1)^r\big{(}\omega(g_0^{-1},\delta_{r-1}(\delta_{r}\vec{f})^{h_0},f_{r-1}^{h_0}g_0)-\cancel{\omega(g_0^{-1},\delta_{r-1}(\delta_{r}\vec{f})^{h_0},g_0)}\big{)}; $$
hence, updating the difference 
$\epsilon_1:=\epsilon_0+(-1)^{r+1}T_r-(-1)^rS_{r-1}$ becomes 
$$\mld \epsilon_1= & (-1)^{r+1}\Big{(}\rho_0^1(i(g_0))^{-1}\omega((\vec{f})^{h_0})-\omega(g_0^{-1}f_1^{h_0},\delta_0(\delta_r\vec{f})^{h_0},f_r^{h_0}g_0)+\omega(g_0^{-1}f_1^{h_0},\delta_0(\delta_r\vec{f})^{h_0},g_0)+ \\
	+\rho_0^1(i(f_1^{h_0i(g_0)}))^{-1}\big{(}\omega(g_0^{-1},\delta_{r-1}(\delta_0\vec{f})^{h_0},f_r^{h_0}g_0)-\omega(g_0^{-1},\delta_{r-1}(\delta_0\vec{f})^{h_0},g_0)\big{)}\Big{)}. $$
In general, for $2\leq n\lt r$,
$$\mld T_n= & \rho_0^1(i(f_1^{h_0i(g_0)}))^{-1}\big{(}\omega((\vec{f}_{[2,r-n]})^{h_0i(g_0)},g_0^{-1},(\vec{f}_{[r-n+1,r)})^{h_0},f_r^{h_0}g_0)+ \\
\qquad\qquad\qquad\qquad\qquad\qquad-\omega((\vec{f}_{[2,r-n]})^{h_0i(g_0)},g_0^{-1},(\vec{f}_{[r-n+1,r)})^{h_0},g_0)\big{)}+ \\
		+\sum_{k=1}^{r-(n+1)}(-1)^k\big{(}\omega(\delta_k(\vec{f}_{[1,r-n]})^{h_0i(g_0)},g_0^{-1},(\vec{f}_{[r-n+1,r)})^{h_0},f_r^{h_0}g_0)+ \\
\qquad\qquad\qquad\qquad\qquad\qquad-\omega(\delta_k(\vec{f}_{[1,r-n]})^{h_0i(g_0)},g_0^{-1},(\vec{f}_{[r-n+1,r)})^{h_0},g_0)\big{)}+ \\
	+(-1)^{r-n}\big{(}\omega((\vec{f}_{[1,r-n)})^{h_0i(g_0)},g_0^{-1}f_{r-n}^{h_0},(\vec{f}_{[r-n+1,r)})^{h_0},f_r^{h_0}g_0)+ \\
\qquad\qquad\qquad\qquad\qquad\qquad-\omega((\vec{f}_{[1,r-n)})^{h_0i(g_0)},g_0^{-1}f_{r-n}^{h_0},(\vec{f}_{[r-n+1,r)})^{h_0},g_0)\big{)}+ \\
	+(-1)^{r-n+1}\big{(}\omega((\vec{f}_{[1,r-n]})^{h_0i(g_0)},g_0^{-1}f_{r-n+1}^{h_0},(\vec{f}_{[r-n+2,r)})^{h_0},f_r^{h_0}g_0)+ \\
\qquad\qquad\qquad\qquad\qquad\qquad-\omega((\vec{f}_{[1,r-n]})^{h_0i(g_0)},g_0^{-1}f_{r-n+1}^{h_0},(\vec{f}_{[r-n+2,r)})^{h_0},g_0)\big{)}+ \\
	+\sum_{k=1}^{n-2}(-1)^{r-n+1+k}\big{(}\omega((\vec{f}_{[1,r-n]})^{h_0i(g_0)},g_0^{-1},\delta_k(\vec{f}_{[r-n+1,r)})^{h_0},f_r^{h_0}g_0)+ \\
\qquad\qquad\qquad\qquad\qquad\qquad-\omega((\vec{f}_{[1,r-n]})^{h_0i(g_0)},g_0^{-1},\delta_k(\vec{f}_{[r-n+1,r)})^{h_0},g_0)\big{)}+ \\
	+(-1)^{r}\big{(}\omega((\vec{f}_{[1,r-n]})^{h_0i(g_0)},g_0^{-1},(\vec{f}_{[r-n+1,r-1)})^{h_0},(f_{r-1}f_r)^{h_0}g_0)+ \\
\qquad\qquad\qquad\qquad\qquad\qquad-\omega((\vec{f}_{[1,r-n]})^{h_0i(g_0)},g_0^{-1},(\vec{f}_{[r-n+1,r-1)})^{h_0},f_{r-1}^{h_0}g_0)\big{)}+ \\
	+(-1)^{r+1}\big{(}\cancel{\omega((\vec{f}_{[1,r-n]})^{h_0i(g_0)},g_0^{-1},(\vec{f}_{[r-n+1,r)})^{h_0})}-\cancel{\omega((\vec{f}_{[1,r-n]})^{h_0i(g_0)},g_0^{-1},(\vec{f}_{[r-n+1,r)})^{h_0})}\big{)}, $$
and, for $1\leq n\leq r-2$, 
$$\mld S_n= & \rho_0^1(i(f_1^{h_0i(g_0)}))^{-1}\big{(}\omega((\vec{f}_{[2,r-n]})^{h_0i(g_0)},g_0^{-1},(\vec{f}_{[r-n+1,r)})^{h_0},f_r^{h_0}g_0)+ \\
\qquad\qquad\qquad\qquad\qquad\qquad-\omega((\vec{f}_{[2,r-n]})^{h_0i(g_0)},g_0^{-1},(\vec{f}_{[r-n+1,r)})^{h_0},g_0)\big{)}+ \\
 	+\sum_{k=1}^{r-(n+1)}(-1)^k\big{(}\omega(\delta_k(\vec{f}_{[1,r-n]})^{h_0i(g_0)},g_0^{-1},(\vec{f}_{[r-n+1,r)})^{h_0},f_r^{h_0}g_0)+ \\
\qquad\qquad\qquad\qquad\qquad\qquad-\omega(\delta_k(\vec{f}_{[1,r-n]})^{h_0i(g_0)},g_0^{-1},(\vec{f}_{[r-n+1,r)})^{h_0},g_0)\big{)}+ \\
 	+\sum_{k=r-n}^{r-2}(-1)^k\big{(}\omega((\vec{f}_{[1,r-n)})^{h_0i(g_0)},g_0^{-1},(\delta_k\vec{f}_{[r-n,r)})^{h_0},f_r^{h_0}g_0)+ \\
\qquad\qquad\qquad\qquad\qquad\qquad-\omega((\vec{f}_{[1,r-n]})^{h_0i(g_0)},g_0^{-1},(\delta_k\vec{f}_{[r-n,r)})^{h_0},g_0)\big{)}+ \\
	+(-1)^{r-1}\big{(}\omega((\vec{f}_{[1,r-n)})^{h_0i(g_0)},g_0^{-1},(\vec{f}_{[r-n,r-1)})^{h_0},(f_{r-1}f_r)^{h_0}g_0)+ \\
\qquad\qquad\qquad\qquad\qquad\qquad-\cancel{\omega((\vec{f}_{[1,r-n)})^{h_0i(g_0)},g_0^{-1},(\vec{f}_{[r-n,r-1)})^{h_0},g_0)}\big{)}+ \\
	+(-1)^r\big{(}\omega((\vec{f}_{[1,r-n)})^{h_0i(g_0)},g_0^{-1},(\vec{f}_{[r-n,r-1)})^{h_0},f_{r-1}^{h_0}g_0)+ \\
\qquad\qquad\qquad\qquad\qquad\qquad-\cancel{\omega((\vec{f}_{[1,r-n)})^{h_0i(g_0)},g_0^{-1},(\vec{f}_{[r-n,r-1)})^{h_0},g_0)}\big{)}. $$
Thus, defining inductively 
$\epsilon_{n+1}:=\epsilon_n+(-1)^{r+1-n}(T_{r-n}+S_{r-(n+1)})$,
$$\mld \epsilon_{r-1}= & (-1)^{r+1}\rho_0^1(i(g_0))^{-1}\omega((\vec{f})^{h_0})+ \\
	-\Big{(}\rho_0^1(i(f_1^{h_0i(g_0)}))^{-1}\big{(}\omega((\vec{f}_{[2,r)})^{h_0i(g_0)},g_0^{-1},f_r^{h_0}g_0)-\omega((\vec{f}_{[2,r)})^{h_0i(g_0)},g_0^{-1},g_0)\big{)}+ \\
 	+\sum_{k=1}^{r-2}(-1)^k\big{(}\omega(\delta_k(\vec{f}_{[1,r)})^{h_0i(g_0)},g_0^{-1},f_r^{h_0}g_0)-\omega(\delta_k(\vec{f}_{[1,r)})^{h_0i(g_0)},g_0^{-1},g_0)\big{)}+ \\
	+(-1)^{r-1}\big{(}\omega((\vec{f}_{[1,r-1)})^{h_0i(g_0)},g_0^{-1}f_{r-1}^{h_0},f_r^{h_0}g_0)-\omega((\vec{f}_{[1,r-1)})^{h_0i(g_0)},g_0^{-1}f_{r-1}^{h_0},g_0)\big{)}\Big{)}. $$
Naturally, 
$(\Delta\delta_{(1)}-\delta'\Delta)\omega(\gamma;\vec{f})=\epsilon_{r-1}+T_1$;
therefore,
$$\mld (\Delta\delta_{(1)}-\delta'\Delta)\omega(\gamma;\vec{f})= & (-1)^{r+1}\rho_0^1(i(g_0))^{-1}\omega((\vec{f})^{h_0})+ \\
	+(-1)^{r}\Big{(}\omega((\delta_r\vec{f})^{h_0i(g_0)},g_0^{-1}f_r^{h_0}g_0)-\cancelto{0}{\omega((\delta_r\vec{f})^{h_0i(g_0)},g_0^{-1}g_0)}\Big{)} $$
and the result follows.

\epf

\proposition\label{r-page}
Let $r>1$ and $\omega\in C(\G_p^q\times G^r,W)$, then 
\begin{equation}\label{r-UpToHom}
(-1)^r(\delta\partial-\partial\delta)\omega=(\Delta\delta_{(1)}-\delta_{(1)}\Delta)\omega\in C(\G_{p+1}^{q+1}\times G^r,W)
\end{equation} \endproposition

\pf
Eq. (\ref{r-UpToHom}) follows from an argument analog to the one in 
Proposition \ref{q=0r-page} after noticing that the left hand side takes 
the form of Eq. (\ref{theLHS r-pag}) (cf. Eq. (\ref{goal r-pag})). Let 
$\vec{\gamma}\in\G_{p+1}^{q+1}$ be as in Eq. (\ref{gammaMatrix}) and 
$\vec{f}=(f_1,...,f_r)\in G^r$, then
$$\mld \delta\partial\omega(\vec{\gamma};\vec{f}) & =\partial\omega(\delta_0\vec{\gamma};(\vec{f})^{h_{11}i(g_{11})})+ \\
\qquad +\sum_{k=1}^{q}(-1)^k\partial\omega(\delta_k\vec{\gamma};\vec{f})+(-1)^{q+1}\rho_0^1(h_{(q+1)1}i(g_{(q+1)1}))^{-1}\partial\omega(\delta_{q+1}\vec{\gamma};\vec{f}). $$
Assuming the convention of Eq. (\ref{p-repShorthand}),
$$ \mld \delta\partial\omega(\vec{\gamma};\vec{f}) & =\rho^{q}((\delta_0\vec{\gamma})_{\bullet 1})\omega(\vec{\gamma}_{1,1};(\vec{f})^{h_{11}i(g_{11})})+\sum_{j=1}^{p+1}(-1)^j\omega(\partial_j\delta_0\vec{\gamma};(\vec{f})^{h_{11}i(g_{11})})+ \\
        \quad +\sum_{k=1}^{q}(-1)^k\Big{[}\rho^{q+1}(\vec{\gamma}_{\bullet 1})\omega(\partial_0\delta_k\vec{\gamma};\vec{f})+\sum_{j=1}^{p+1}(-1)^j\omega(\partial_j\delta_k\vec{\gamma};\vec{f})\Big{]}+ $$
$$ +(-1)^{q+1}\rho_0^1(h_{(q+1)1}i(g_{(q+1)1}))^{-1}\Big{[}\rho^q((\delta_{q+1}\vec{\gamma})_{\bullet 1})\omega(\vec{\gamma}_{0,q};\vec{f})+\sum_{j=1}^{p+1}(-1)^j\omega(\partial_j\delta_{q+1}\vec{\gamma};\vec{f})\Big{]}, $$
and  
$$\mld \partial & \delta\omega(\vec{\gamma};\vec{f})=\rho^{q+1}(\vec{\gamma}_{\bullet 1})\delta\omega(\partial_0\vec{\gamma};\vec{f})+\sum_{j=1}^{p+1}(-1)^j\delta\omega(\partial_j\vec{\gamma};\vec{f}) \\
   	=\rho^{q+1}(\vec{\gamma}_{\bullet 1})\Big{[}\omega(\vec{\gamma}_{1,1};(\vec{f})^{h_{11}})+\sum_{k=1}^{q}(-1)^k\omega(\delta_k\partial_0\vec{\gamma};\vec{f})+(-1)^{q+1}\rho_0^1(h_{(q+1)1})^{-1}\omega(\vec{\gamma}_{q+1,1};\vec{f})\Big{]}+\\
    \quad +\sum_{j=1}^{p+1}(-1)^j\Big{[}\omega(\delta_0\partial_j\vec{\gamma};(\vec{f})^{h_{11}i(g_{11})})+ \\
    \qquad\qquad +\sum_{k=1}^{q}(-1)^k\omega(\delta_k\partial_j\vec{\gamma};\vec{f})+(-1)^{q+1}\rho_0^1(h_{(q+1)1}i(g_{(q+1)1}))^{-1}\omega(\delta_{q+1}\partial_j\vec{\gamma};\vec{f})\Big{]}. $$
Aside from the visible terms in common, using the following identity for 
$y\in H$ and $x,z\in G$:
$$\rho_0^1(yi(x))^{-1}\rho_0^1(i(z))^{-1}=\rho_0^1(i(z)yi(x))^{-1}=\rho_0^1(yi(z^y)i(x)))^{-1}=\rho_0^1(i(z^yx))^{-1}\rho_0^1(y)^{-1}, $$
applied to $y=h_{(q+1)1}$, $x=g_{(q+1)1}$ and 
$z=pr_G(\gamma_{11}\vJoin ...\vJoin\gamma_{q1})=g_{11}^{h_{21}...h_{q1}}g_{21}^{h_{31}...h_{q1}}...g_{q-1}^{h_{q1}}g_{q1}$ 
(cf. Lemma \ref{multiprods}), the difference becomes
\begin{equation}\label{goal r-pag}
(\delta\partial-\partial\delta)\omega(\vec{\gamma};\vec{f})=\rho^{q+1}((\delta_0\vec{\gamma})_{\bullet 1})\Big{(}\omega(\vec{\gamma}_{1,1};(\vec{f})^{h_{11}i(g_{11})})-\rho_0^1(i(g_{11}^{h_{21}...h_{q1}}))^{-1}\omega(\vec{\gamma}_{1,1};(\vec{f})^{h_{11}})\Big{)}.
\end{equation}

\epf

\subsection{Higher Difference maps and outlook}\label{HDiff} 
Let $\omega\in C^{p,q}_r(\G,\phi)$. In Section \ref{sec-tetrahedralCx}, 
it is explained that if the total differential is defined by Eq. 
(\ref{preDiff}), $\nabla^2\omega$ is not necessarily zero; more 
specifically, in general, 
$$(\nabla^2\omega)^{p+1,q+1}_r=\partial\delta\omega-\delta\partial\omega\neq 0. $$
We can sum up Subsection \ref{DiffMaps} as saying that, 
if the total differential is redefined by 
\begin{equation}\label{preDiff2}
\nabla=(-1)^p\Big{(}\delta_{(1)}+\partial+\Delta+(-1)^r\delta\Big{)},
\end{equation} 
$(\nabla^2\omega)^{p+1,q+1}_r$ is ensured to vanish. However, adding 
the difference maps does not yet imply that $\nabla^2\omega$ vanishes; 
indeed, $(\nabla^2\omega)^{p+2,q+1}_{r-1}$ and 
$(\nabla^2\omega)^{p+1,q+2}_{r-1}$ need not be zero, as $\Delta$ does 
not necessarily commute with neither $\partial$ nor $\delta$. In 
Subsection \ref{H^2}, these non-vanishing coordinates of 
$\nabla^2\omega$ are studied in the case $p+q+r=1$, and it is explained 
that one can ensure again $\nabla^2\omega=0$ by adding the second 
difference maps $\Delta_{2,1}$ and $\Delta_{1,2}$ (cf. Eq.'s 
(\ref{Delta21}) and (\ref{Delta12})). Thus, redefining the total 
differential by Eq. (\ref{preDiff3}) forces 
$(\nabla^2\omega)^{p+2,q+1}_{r-1}$ and $(\nabla^2\omega)^{p+1,q+2}_{r-1}$ 
to be zero, but, in general, as with the first difference maps, the 
$\nabla$ of Eq. (\ref{preDiff3}) does not yet square to zero. Adding 
higher difference maps creates new non-vanishing coordinates, which, in 
turn, can be made to zero by adding further higher difference maps. 
Ultimately, if the total differential is defined by 
\begin{equation}\label{FinalDiff}
\nabla=(-1)^p\Big{(}\delta_{(1)}+\sum_{a+b>0}(-1)^{(a+1)(r+b+1)}\Delta_{a,b}\Big{)}
\end{equation}
for some 
$$\bfig
\morphism<900,0>[\Delta_{a,b}:C^{p,q}_r(\G,\phi)`C^{p+a,q+b}_{r+1-(a+b)}(\G,\phi);]
\efig$$
where we set $\Delta_{1,0}:=\partial$, $\Delta_{0,1}:=\delta$, 
$\Delta_{1,1}:=\Delta$ and $\Delta_{a,0}=\Delta_{0,b}=0$ whenever 
$a,b>1$, the following is a rephrasing of the relation $\nabla^2=0$:

\thm\label{rePhrase}
Let $\omega\in C^{p,q}_r(\G,\phi)$. Then, $\nabla^2\omega=0$ if and only 
if 
\begin{itemize}
\item[i)] $\delta_{(1)}^2\omega=0$; 
\item[ii)] for all $0\leq n\leq r+1$ and $0\leq m\leq n$,
\begin{equation}\label{umbrella}
\sum_{0\lt i+j\lt n}(-1)^{i(i-n)+(i+1)(j+1)}\Delta_{n-m-i,m-j}\circ\Delta_{i,j}\omega=(-1)^{r}[\delta_{(1)},\Delta_{n-m,m}]\omega,
\end{equation}
where $[\cdot,\cdot]$ stands for the commutator of operators; and 
\item[iii)] for all $0\lt m\leq r+1$, 
\begin{equation}\label{wall}
\sum_{0\lt i+j\lt r+2}(-1)^{i(i-r)+(i+1)(j+1)}\Delta_{r+2-m-i,m-j}\circ\Delta_{i,j}\omega=(-1)^{r+1}\Delta_{r+2-m,m}\circ\delta_{(1)}\omega.
\end{equation}
\end{itemize}
Here, we assume the convention that $\Delta_{a,b}=0$ whenever  
$a\lt 0$ or $b\lt 0$.
\endthm

As it is briefly mentioned above, in Section \ref{sec-tetrahedralCx}, 
the necessary higher difference maps for the cases $p+q+r\leq 1$ are 
defined and they are shown to verify the relations of Theorem 
\ref{rePhrase}; moreover, that $\delta_{(1)}$ squares to zero is 
established (cf. Lemmas \ref{q=0r-cx} and \ref{r-cx}), as well as the 
relations of Eq. (\ref{umbrella}) in the cases when $n\leq 2$ (cf. 
Subsections \ref{p-dir}, \ref{q-dir}, \ref{pq-pag}, \ref{DiffMaps}). 
Although as of the writing of this paper, we are unable to provide a 
general formula for the higher difference maps $\Delta_{a,b}$, we 
present evidence that suggests that one can find such maps ultimately 
turning $(C_{tot}(\G,\phi),\nabla)$ into a complex. In particular, with 
the formulas we include in the appendix, the complex (\ref{trunc3}) is 
extended up to degree $5$ (see \ref{leq5}).

We devote the remainder of this section to define the necessary higher 
difference maps to prove that Eq. (\ref{wall}) holds in general for 
$m\in\lbrace 1,r+1\rbrace$. 

Let
$$\Delta_{r,1}:C^{p,q}_r(\G,\phi)\to<350> C^{p+r,q+1}_0(\G,\phi) $$
be defined by
\begin{equation}\label{FrontHighDiff r,1}
\Delta_{r,1}\omega(\vec{\gamma}):=\rho_0^0(t_{p}(\partial_0^r\gamma_1)...t_{p}(\partial_0^r\gamma_{q+1}))\circ\phi\big{(}\omega(\partial_0^{r}\delta_0\vec{\gamma};g_{1r},g_{1r-1},...,g_{12},g_{11})\big{)},
\end{equation}
for $\omega\in C(\G_p^q\times G^r,W)$ and 
$\vec{\gamma}\in\G_{p+r}^{q+1}$ as in Eq. (\ref{gammaMatrix}).

\proposition\label{Wall m=1}
Let $\omega\in C(\G_p^q\times G^r,W)$, then
\begin{equation}\label{Eq r+2,1}
((-1)^{r+1}\Delta_{r,1}\circ\partial+\partial\circ\Delta_{r,1})\omega=(-1)^{r+1}\Delta_{r+1,1}\circ\delta_{(1)}\omega . 
\end{equation}\endproposition

\pf
To prove Eq. (\ref{Eq r+2,1}), one needs to consider two separate 
cases: $q=0$ and $q>0$. 

For $q=0$, let $\gamma\in\G_{p+r+1}$ and let 
$(g_0,...,g_{p+r};h)\in G^{p+r+1}\times H$ be its image under the 
isomorphism of Remark \ref{2GpNerve}, then
$$\partial\Delta_{r,1}\omega(\gamma)=\sum_{j=0}^{p+r+1}(-1)^j\Delta_{r,1}\omega(\partial_j\gamma). $$
Using Eq. (\ref{2GpFaceMaps}), we compute
\begin{equation}\label{Del rFaceMaps}
\partial_0^r\partial_j\gamma =\cases{
    (g_{r+1},...,g_{p+r};h)            & if $0\leq j\leq r$  \cr
    (g_r,...,g_jg_{j-1},...,g_{p+r};h) & if $r\lt j\leq r+p$ \cr
    (g_r,...,g_{p+r-1};hi(g_{p+r}))    & if $j=r+p+1$, }
\end{equation}
thus yielding
$$\mld & \partial\Delta_{r,1}\omega(\gamma) =\phi\Bigg{[}\rho_0^1(hi(g_{r+p}...g_{r+1}))\Big{(}\omega(g_r,...,g_1)-\omega(g_r,...,g_2,g_1g_0)+...+ \\
               \quad +(-1)^{r-1}\omega(g_r,g_{r-2}g_{r-1},...,g_1,g_0)+(-1)^{r}\omega(g_rg_{r-1},g_{r-2}...,g_1,g_0)\Big{)}+ \\
               \qquad +\sum_{j=r+1}^{r+p}(-1)^j\rho_0^1(hi(g_{r+p}...g_{r}))\omega(g_{r-1},...,g_0)+ \\
               \qquad\quad +(-1)^{r+p+1}\rho_0^1(hi(g_{r+p})i(g_{r+p-1}...g_{r}))\omega(g_{r-1},...,g_0)\Bigg{]}. $$

Since the terms in the sum have got no dependence on $j$ and are equal 
to the last term,   
$$\partial\Delta_{r,1}\omega(\gamma)=\cases{
    (-1)^{r+1}\phi\Big{[}\rho_0^1(hi(g_{p+r}...g_{r+1}))\Big{(}\delta_{(1)}\omega(g_{r},...,g_0)-\rho_0^1(i(g_r))\omega(g_{r-1},...,g_0)\Big{)}\Big{]} & if $p$ is odd, \cr
    (-1)^{r+1}\phi\Big{(}\rho_0^1(hi(g_{p+r}...g_{r+1}))\delta_{(1)}\omega(g_{r},...,g_0)\Big{)} & otherwise. } $$
On the other hand, 
$$\Delta_{r,1}\partial\omega(\gamma)=\cases{\rho_0^0(hi(g_{p+r}...g_r))\circ\phi(\omega(g_{r-1},...,g_0)) & if $p$ is odd, \cr
     0 & otherwise; } $$
thus, in either case, 
\begin{equation}
(\partial\circ\Delta_{r,1}\omega+(-1)^{r+1}\Delta_{r,1}\circ\partial\omega)(\gamma)=(-1)^{r+1}(\Delta_{r+1,1}\circ\delta_{(1)}\omega)(\gamma).
\end{equation}

For $q\geq 1$, let $\vec{\gamma}\in\G_{p+r+1}^{q+1}$ be as in Eq. 
(\ref{gammaMatrix}), then
$$\mld \partial\Delta_{r,1}\omega(\vec{\gamma}) & =\sum_{j=0}^{p+r+1}(-1)^j\Delta_{r,1}\omega(\partial_j\vec{\gamma}) \\
      =\phi\Big{[}\rho_0^1(t_{p}(\partial_0^{r+1}\gamma_1)...t_{p}(\partial_0^{r+1}\gamma_{q+1}))\omega(\partial_0^r\delta_0(\partial_0\vec{\gamma});g_{1r+1},...,g_{12})+ \\
       \qquad -\rho_0^1(t_{p}(\partial_0^r\partial_1\gamma_1)...t_{p}(\partial_0^r\partial_1\gamma_{q+1}))\omega(\partial_0^r\delta_0(\partial_1\vec{\gamma});g_{1r+1},...,g_{13},g_{12}g_{11})+...+ $$
$$\mld \quad +(-1) & ^{r-1}\rho_0^1(t_{p}(\partial_0^{r}\partial_{r-1}\gamma_1)...t_{p}(\partial_0^{r}\partial_{r-1}\gamma_{q+1}))\omega(\partial_0^r\delta_0(\partial_{r-1}\vec{\gamma});g_{1r+1},g_{1r}g_{1r-1},...,g_{12},g_{11})+ \\
       +(-1)^{r}\rho_0^1(t_{p}(\partial_0^{r}\partial_{r}\gamma_1)...t_{p}(\partial_0^{r}\partial_{r}\gamma_{q+1}))\omega(\partial_0^r\delta_0(\partial_{r}\vec{\gamma});g_{1r+1}g_{1r},g_{1r-1}...,g_{12},g_{11})+ \\
       \quad +\sum_{j=r+1}^{p+r+1}(-1)^j\rho_0^1(t_{p}(\partial_0^{r}\partial_j\gamma_1)...t_{p}(\partial_0^{r}\partial_j\gamma_{q+1}))\omega(\partial_0^r\delta_0(\partial_j\vec{\gamma})_{0,0};g_{1r},...,g_{11})\Big{]}. $$
Using Eq. (\ref{Del rFaceMaps}) along with 
$h_{ab}=h_{ab+1}i(g_{ab+1})=h_{ap+r}i(g_{ap+r}...g_{ab+1})$, one gets 
$$\mld \partial\Delta_{r,1}\omega(\vec{\gamma})=(-1)^{r+1}\phi & \Bigg{[}\rho_0^1(h_{1r+1}...h_{(q+1)(r+1)})\Big{(}\delta_{(1)}\omega(\partial_0^{r+1}\delta_0\vec{\gamma};g_{1r+1},...,g_{11})+ \\
	\qquad -\rho_0^1(i(g_{1r+1}^{h_{2r+1}...h_{(q+1)(r+1)}}))\omega(\partial_0^{r+1}\delta_0\vec{\gamma};g_{1r},...,g_{11})\Big{)}+ \\
    \qquad\quad -\sum_{j=r+1}^{p+r+1}(-1)^{j-r}\rho_0^1(h_{1r}...h_{(q+1)r})\omega(\partial_0^r\delta_0(\partial_j\vec{\gamma});g_{1r},...,g_{11})\Bigg{]}. $$
Notice that, as opposed to the case $q=0$, due to the explicit 
dependence on $j$, the terms in the sum do not cancel one another. 

On the other hand,
$$\mld  & \Delta_{r,1}\partial\omega(\vec{\gamma})=\rho_0^0(t_{p+1}(\partial_0^r\gamma_1)...t_{p+1}(\partial_0^r\gamma_{q+1}))\circ\phi\Big{(}\partial\omega(\partial_0^r\delta_0\vec{\gamma};g_{1r},g_{1r-1},...,g_{12},g_{11})\Big{)} \\
        =\rho_0^0(h_{1r}...h_{(q+1)r})\circ\phi\Big{(}\rho_0^1(i(pr_G(\gamma_{2r+1}\vJoin ...\vJoin\gamma_{(q+1)(r+1)})))^{-1}\omega(\partial_0^{r+1}\delta_0\vec{\gamma};g_{1r},...,g_{11})+ \\
        \qquad +\sum_{j=1}^{p+1}(-1)^j\omega(\partial_j\partial_0^r\delta_0\vec{\gamma};g_{1r},g_{1r-1},...,g_{12},g_{11})\big{)}\Big{)}; $$
thus, by means of $\partial_0^r\partial_j=\partial_{j-r}\partial_0^r$ 
for 
$j\geq r+1$ and  
$$h_{1r}...h_{(q+1)r}=h_{1r+1}...h_{(q+1)(r+1)}i(pr_G(\gamma_{1r+1}\vJoin ...\vJoin\gamma_{(q+1)(r+1)})), $$
multiplying by the factor of $(-1)^{r+1}$, it follows that
$$\mld (\partial\Delta_{r,1}\omega+(-1)^{r+1}\Delta_{r,1}\partial\omega)(\vec{\gamma}) & =(-1)^{r+1}\rho_0^0(h_{1r+1}...h_{(q+1)(r+1)})\circ\phi\Big{(}\delta_{(1)}\omega(\partial_0^{r+1}\delta_0\vec{\gamma};g_{1r+1},...,g_{11})\Big{)} \\
	=(-1)^{r+1}(\Delta_{r+1,1}\circ\delta_{(1)}\omega)(\vec{\gamma}), $$
as desired. 

\epf

Let
$$\Delta_{1,r}:C^{p,q}_r(\G,\phi)\to<350> C^{p+1,q+r}_0(\G,\phi) $$
be defined by
\begin{equation}\label{FrontHighDiff 1,r}
\Delta_{1,r}\omega(\vec{\gamma}):=\rho_0^0(t_{p}(\partial_0\gamma_1)...t_{p}(\partial_0\gamma_{q+r}))\circ\phi\big{(}\omega(\partial_0\delta^{r}_0\vec{\gamma};g_{11}^{h_{21}...h_{r1}},g_{21}^{h_{31}...h_{r1}},...,g_{(r-1)1}^{h_{r1}},g_{r1})\big{)},
\end{equation}
for $\omega\in C(\G_p^q\times G^r,W)$ and 
$\vec{\gamma}\in\G_{p+1}^{q+r}$ as in Eq. (\ref{gammaMatrix}). 

\proposition\label{Wall m=r+1}
Let $\omega\in C(\G_p^q\times G^r,W)$, then
\begin{equation}\label{Eq r+2,r+1} 
(\Delta_{1,r}\circ\delta+(-1)^{r+1}\delta\circ\Delta_{1,r})\omega=(-1)^{r+1}\Delta_{1,r+1}\circ\delta_{(1)}\omega .
\end{equation}\endproposition

\pf
Let $\vec{\gamma}\in\G_{p+1}^{q+r+1}$ be as in Eq. (\ref{gammaMatrix}). 
Then,

$$\mld  & \delta\Delta_{1,r}\omega(\vec{\gamma})=\rho_0^0(t_{p+1}(\gamma_1))\Delta_{1,r}\omega(\delta_0\vec{\gamma})+\sum_{j=1}^{q+r+1}(-1)^j\Delta_{1,r}\omega(\delta_j\vec{\gamma}) \\
        =\rho_0^0(h_{11}i(g_{11}))\phi\Big{[}\rho_0^1(t_{p}(\partial_0\gamma_2)...t_{p}(\partial_0\gamma_{q+r+1}))\omega(\partial_0\delta_0^r(\delta_0\vec{\gamma});g_{21}^{h_{31}...h_{(r+1)1}},...,g_{r1}^{h_{(r+1)1}},g_{(r+1)1})+ \\
        \quad -\rho_0^1(t_{p}(\partial_0(\gamma_1\vJoin\gamma_2))...t_{p}(\partial_0\gamma_{q+r+1}))\omega(\partial_0\delta_0^r(\delta_1\vec{\gamma});(g_{11}^{h_{21}}g_{21})^{h_{31}...h_{(r+1)1}},...,g_{r1}^{h_{(r+1)1}},g_{(r+1)1})+ \\
        +\sum_{j=2}^{r}(-1)^j\rho_0^1(t_{p}(\partial_0\gamma_1)...t_{p}(\partial_0(\gamma_j\vJoin\gamma_{j+1}))...t_{p}(\partial_0\gamma_{q+r+1}))\omega(\partial_0\delta_0^r(\delta_j\vec{\gamma});...,(g_{j1}^{h_{(j+1)1}}g_{(j+1)1})^{h_{(j+2)1}...h_{(r+1)1}},...)+ \\
        +\sum_{j=r+1}^{r+q}(-1)^j\rho_0^1\rho_0^1(t_{p}(\partial_0\gamma_1)...t_{p}(\partial_0(\gamma_j\vJoin\gamma_{j+1}))...t_{p}(\partial_0\gamma_{q+r+1}))\omega(\partial_0\delta_0^r(\delta_j\vec{\gamma});g_{11}^{h_{21}...h_{r1}},...,g_{(r-1)1}^{h_{r1}},g_{r1})+ \\
        \quad\qquad +(-1)^{r+q+1}\rho_0^1(t_{p}(\partial_0\gamma_1)...t_{p}(\partial_0\gamma_{q+r}))\omega(\partial_0\delta_0^r(\delta_{r+q+1}\vec{\gamma});g_{11}^{h_{21}...h_{r1}},...,g_{(r-1)1}^{h_{r1}},g_{r1})\Big{]} \\
        =\rho_0^0(h_{11}...h_{(q+r+1)1})\phi\Big{[}\rho_0^1(i(g_{11}^{h_{21}...h_{(q+r+1)1}}))\omega(\partial_0\delta_0^{r+1}\vec{\gamma};g_{21}^{h_{31}...h_{(r+1)1}},...,g_{r1}^{h_{(r+1)1}},g_{(r+1)1})+ \\
        \qquad +\sum_{k=1}^r(-1)^k\omega(\partial_0\delta_0^{r+1}\vec{\gamma};\delta_k(g_{11}^{h_{21}...h_{(r+1)1}},...,g_{r1}^{h_{(r+1)1}},g_{(r+1)1}))+ \\
        \qquad\quad +\sum_{j=r+1}^{r+q}(-1)^j\omega(\partial_0\delta_0^r(\delta_j\vec{\gamma});g_{11}^{h_{21}...h_{r1}},...,g_{(r-1)1}^{h_{r1}},g_{r1})) \\
        \qquad\qquad +(-1)^{r+q+1}\rho_0^1(h_{(q+r+1)1})^{-1}\omega(\partial_0\delta_0^r(\delta_{r+q+1}\vec{\gamma});g_{11}^{h_{21}...h_{r1}},...,g_{(r-1)1}^{h_{r1}},g_{r1})\Big{]} $$
and        
$$\mld \Delta_{1,r}\delta\omega(\vec{\gamma}) & =\rho_0^1(t_{p}(\partial_0\gamma_1)...t_{p}(\partial_0\gamma_{q+r+1}))\circ\phi\Big{(}\delta\omega(\partial_0\delta_0^r\vec{\gamma};g_{11}^{h_{21}...h_{r1}},...,g_{(r-1)1}^{h_{r1}},g_{r1})\Big{)} \\
        =\rho_0^0(h_{11}...h_{(q+r+1)1})\phi\Big{[}\omega(\partial_0\delta_0^{r+1}\vec{\gamma};(g_{11}^{h_{21}...h_{r1}})^{h_{(r+1)1}},...,(g_{(r-1)1}^{h_{r1}})^{h_{(r+1)1}},g_{r1}^{h_{(r+1)1}}))+ \\
        \quad +\sum_{j=1}^{q}(-1)^j\omega(\delta_j\partial_0\delta_0^r\vec{\gamma};g_{11}^{h_{21}...h_{r1}},...,g_{(r-1)1}^{h_{r1}},g_{r1})+ \\
        \qquad +(-1)^{q+1}\rho_0^1(t_p(\partial_0\gamma_{q+r+1}))^{-1}\omega(\delta_{q+1}\partial_0\delta_0^r\vec{\gamma};g_{11}^{h_{21}...h_{r1}},...,g_{(r-1)1}^{h_{r1}},g_{r1})\big{)}\Big{]}. $$
Since $\delta_0^r\delta_j=\delta_{j-r}\delta_0^r$ for $j\geq r+1$ and 
$t_p(\partial_0\gamma_{q+r+1})=h_{(q+r+1)1}$, it follows that
$$\mld (\delta\Delta_{1,r}+(-1)^{r+1}\Delta_{1,r}\delta)\omega(\vec{\gamma}) & =\rho_0^0(h_{11}...h_{(q+r+1)1})\circ\phi\Big{(}\delta_{(1)}\omega(\partial_0\delta_0^{r+1}\vec{\gamma};g_{11}^{h_{21}...h_{(r+1)1}},...,g_{r1}^{h_{(r+1)1}},g_{(r+1)1})\Big{)} \\
	=(\Delta_{1,r+1}\circ\delta_{(1)}\omega)(\vec{\gamma}), $$
as desired.       

\epf

\rem\label{final} Propositions \ref{Wall m=1} and \ref{Wall m=r+1} are 
the simplest of relations in Theorem \ref{rePhrase}, both due to the 
definition of the higher difference maps involved and to the number of 
terms in Eq. (\ref{wall}). Continuing Remark \ref{upToIsoIndeed}, these 
results prove that though the difference maps do not commute in general 
with $\partial$ and $\delta$, they do so up to isomorphism in the 
$2$-vector space. \endrem

\section*{Appendix}\label{sec-App}

We devote this appendix to give the rather cumbersome formulas for  
some families of higher difference maps and all maps necessary to extend 
(\ref{trunc3}) up to degree $5$.

We introduce further notation in order to abbreviate the formulas. 
Associated to $\vec{\gamma}\in\G_{p+a}^{q+b}$ as in Eq. 
(\ref{gammaMatrix}) for $a,b>0$, define the sub-matrices
\begin{itemize}
\item $\vec{\gamma}^{a,b}_\sqcap:=(\gamma_{ij})_{1\leq i\leq b;1\leq j\leq a}\in\G_a^b$, 
and
\item $\vec{\gamma}^{a,b}_\sqcup:=(\gamma_{ij})_{b+1\leq i\leq q+b;1\leq j\leq a}\in\G_a^q.$
\end{itemize}
Also, define $\parallel\vec{\gamma}\parallel$ to be the 
{\em full product} of the coordinates of $\vec{\gamma}$, i.e. 
$$\parallel\vec{\gamma}\parallel:=\partial_1^{p+a-1}\delta_1^{q+b-1}\vec{\gamma}=(\gamma_{11}\Join ...\Join\gamma_{1p+a})\vJoin ...\vJoin(\gamma_{(q+b)1}\Join ...\Join\gamma_{(q+b)(p+a)}), $$
and $\parallel\vec{\gamma}\parallel_G:=pr_G(\parallel\vec{\gamma}\parallel)$.  

\subsection*{Front page}\label{To r=0}
Each difference map landing on the front page 
\[ \Delta_{a,b}:C^{p,q}_{a+b-1}(\G,\phi)\to<350> C^{p+a,q+b}_0(\G,\phi) \]
is defined by
\[
\Delta_{a,b}\omega(\vec{\gamma})=\rho_0^0\Big{(}s(\parallel\vec{\gamma}^{a,b}_\sqcap\parallel\vJoin\parallel\vec{\gamma}^{a,b}_\sqcup\parallel)\Big{)}\circ\phi\Bigg{[}\sum_{\alpha\in I_{a,b}}\zeta(\alpha)\omega\Big{(}\partial_0^a\delta_0^b\vec{\gamma};c_{a,b}^\alpha(\vec{\gamma}^{a,b}_\sqcap)\Big{)}\Bigg{]},
\]
for $\omega\in C(\G_p^q\times G^{a+b-1},W)$ and 
$\vec{\gamma}\in\G_{p+a}^{q+b}$. Here, $I_{a,b}$ is a set of indices, 
$\zeta(\alpha)$ is a sign and 
$\lbrace c_{a,b}^\alpha\rbrace_{\alpha\in I_{a,b}}$ is a collection of 
maps 
\[
c_{a,b}^\alpha:\G_a^b\to<350> G^{a+b-1}.
\] 

\subsection*{Otherwise}\label{To r>0}
Each difference map landing off the front page 
\[ \Delta_{a,b}:C^{p,q}_r\to<350> C^{p+a,q+b}_{r+1-(a+b)} \]
with $a+b<r+1$ is defined by
$$ \mld \Delta_{a,b}\omega(\vec{\gamma};\vec{f})= & \rho_0^1(i(\parallel\vec{\gamma}^{a,b}_\sqcup\parallel_G))^{-1}\Bigg{[}\sum_{\beta\in J_{a,b}(r)}\zeta(\beta)\omega\Big{(}\partial_0^a\delta_0^b\vec{\gamma};c^{a,b}_\beta(r)\big{(}\vec{f};\vec{\gamma}^{a,b}_\sqcap\big{)}\Big{)}+ \\
 \quad +\rho_0^1(i(\parallel\vec{\gamma}^{a,b}_\sqcap\parallel_G^{s(\parallel\vec{\gamma}^{a,b}_\sqcup\parallel)}))^{-1}\sum_{\alpha\in I_{a,b}}\zeta(\alpha)\omega\Big{(}\partial_0^a\delta_0^b\vec{\gamma};(\vec{f})^{s(\parallel\vec{\gamma}^{a,b}_\sqcap\parallel)},c^\alpha_{a,b}\big{(}\vec{\gamma}^{a,b}_\sqcap\big{)}\Big{)}\Bigg{]}, $$
for $\omega\in C(\G_p^q\times G^r,W)$, $\vec{\gamma}\in\G_{p+a}^{q+b}$ 
and $\vec{f}\in G^{r+1-(a+b)}$.
Again, $J_{a,b}(r)$ is a set of indices, $\zeta(\beta)$ stands for the 
sign of the index $\beta$ and 
$\lbrace c^{a,b}_\beta(r)\rbrace_{\beta\in J_{a,b}(r)}$ is a set of maps
\[
c^{a,b}_\beta(r):G^{r+1-(a+b)}\times\G_a^b\to<350> G^r.
\] 

We limit ourselves to specifying the index set $I_{a,b}$, the values of 
$c_{a,b}^\alpha$ and the signs $\zeta(\alpha)$ by writing them as the 
formal polynomial 
$p_{a,b}=\sum_{\alpha\in I_{a,b}}\zeta(\alpha)c_{a,b}^\alpha$. 
Analogously, we use the formal polynomial 
$p_{a,b}^{(r)}=\sum_{\beta\in J_{a,b}(r)}\zeta(\beta)c^{a,b}_\beta(r)$ 
to indicate $J_{a,b}(r)$, $c^{a,b}_\beta(r)$ and $\zeta(\beta)$. 

$$\mld p_{2,2}\pmatrix{\gamma_{11} & \gamma_{12}\cr
						\gamma_{21} & \gamma_{22}}= & (g_{12}^{h_{22}}g_{22},g_{11}^{h_{21}},g_{21})-((g_{12}g_{11})^{h_{22}},g_{22},g_{21})+(g_{12}^{h_{22}},g_{11}^{h_{22}},g_{22})+ \\
 +(g_{12}^{h_{22}}g_{22},g_{22}^{-1},g_{11}^{h_{22}}g_{22})-(g_{12}^{h_{22}}g_{22},g_{22}^{-1},g_{22}) $$

Using these coordinates to define $\Delta_{2,2}$ implies
\begin{equation}\label{D22 1st}
(\Delta_{2,1}\circ\delta-\Delta_{1,2}\circ\partial-\Delta\circ\Delta-\partial\circ\Delta_{1,2}+\delta\circ\Delta_{2,1})\omega=-\Delta_{2,2}\circ\delta_{(1)}\omega
\end{equation} 
for all $\omega\in C^{p,q}_2(\G,\phi)$.

We point out that there is a certain recurrence in $p_{2,2}$. Indeed, 
one can use $p_{1,2}$, $p_{2,1}$ and the coordinates of the first 
difference map $\Delta$, $c^{1,1}_{2n-1}(2),c^{1,1}_{2n}(2)$ (cf. 
Eq.'s (\ref{c(11)2n-1}) and (\ref{c(11)2n})) to recast $p_{2,2}$ as 
$$\mld p_{2,2}(\vec{\gamma})= & \Big{(}\parallel\partial_0\vec{\gamma}\parallel_G,c_{1,2}(\partial_2\vec{\gamma})\Big{)}-\Big{(}\parallel\vec{\gamma}^{2,1}_\sqcap\parallel_G^{s(\parallel\delta_0\vec{\gamma}\parallel)},c_{2,1}(\delta_0\vec{\gamma})\Big{)}+ \\
	+\Big{(}\big{(}c_{2,1}(\vec{\gamma}^{1,1}_\sqcap)\big{)}^{s(\partial_0\delta_0\vec{\gamma})},c_{1,1}(\partial_0\delta_0\vec{\gamma})\Big{)}+\Big{(}\parallel\partial_0\vec{\gamma}\parallel_G,p_{1,1}^{(2)}(c_{1,1}(\vec{\gamma}^{1,1}_\sqcap);\partial_0\delta_0\vec{\gamma})\Big{)}, $$
where $\vec{\gamma}\in\G_2^2$ is as in Eq. (\ref{gammaMatrix}).	
Using this scheme, one can write more easily the seventeen terms 
in $p_{3,2}$ and the fifteen terms in $p_{2,3}$; however, in turn, 
these are written in part using $p_{2,1}^{(3)}$ and $p_{1,2}^{(3)}$.

$$\mld p_{2,1}^{(3)}\big{(}f;\pmatrix{\gamma_{11}&\gamma_{12}}\big{)}= & ((g_{12}g_{11})^{-1},f^{h_{12}}g_{12},g_{11})-((g_{12}g_{11})^{-1},g_{12},g_{11})+ \\
	-(g_{11}^{-1},g_{12}^{-1},f^{h_{12}}g_{12})+(g_{11}^{-1},g_{12}^{-1},g_{12})-(g_{11}^{-1},f^{h_{11}},g_{11}). $$

$$\mld p_{1,2}^{(3)}\big{(}f;\pmatrix{\gamma_{11}\cr\gamma_{21}}\big{)}= & ((g_{11}^{h_{21}}g_{21})^{-1},(f^{h_{11}}g_{11})^{h_{21}},g_{21})-((g_{11}^{h_{21}}g_{21})^{-1},g_{11}^{h_{21}},g_{21})+ \\
	-(g_{21}^{-1},(g_{11}^{h_{21}})^{-1},(f^{h_{11}}g_{11})^{h_{21}})+(g_{21}^{-1},(g_{11}^{h_{21}})^{-1},g_{11}^{h_{21}})-(g_{21}^{-1},f^{h_{11}i(g_{11})h_{21}},g_{21}). $$

Inductively, for $\vec{f}=(f_1,...,f_{r-1})\in G$ and 
$\vec{\gamma}\in\G_2$ as in Eq. (\ref{gammaMatrix}),
$$\mld p_{2,1}^{(r+1)}(\vec{f};\vec{\gamma})= & \Big{(}f_1^{h_{11}i(g_{11})},p_{2,1}^{(r)}(\delta_0\vec{f};\vec{\gamma})\Big{)}
	+(-1)^{r+1}\Bigg{[}\Big{(}g_{11}^{-1},p_{1,1}^{(r)}(\vec{f};\vec{\gamma})\Big{)}+\Big{(}g_{11}^{-1},(\vec{f})^{h_{11}},g_{11}\Big{)}\Bigg{]}, $$
and for $\vec{\gamma}\in\G^2$ as in Eq. (\ref{gammaMatrix}),
$$\mld p_{1,2}^{(r+1)}(\vec{f};\vec{\gamma})= & \Big{(}f_1^{h_{11}i(g_{11})h_{21}i(g_{21})},p_{1,2}^{(r)}(\delta_0\vec{f};\vec{\gamma})\Big{)}+ \\
	\quad +(-1)^{r+1}\Bigg{[}\Big{(}g_{21}^{-1},\big{(}p_{1,1}^{(r)}(\vec{f};\vec{\gamma})\big{)}^{h_{21}}\Big{)}+\Big{(}g_{21}^{-1},(\vec{f})^{h_{11}i(g_{11})h_{21}},g_{21}\Big{)}\Bigg{]}, $$

Using these coordinates to define $\Delta_{2,1}$ and 
$\Delta_{1,2}$ implies that for all $\omega\in C^{p,q}_r(\G,\phi)$,
\begin{equation}\label{D21 upToHom}
(\Delta\circ\partial+\partial\circ\Delta)\omega=(-1)^r(\delta_{(1)}\circ\Delta_{2,1}-\Delta_{2,1}\circ\delta_{(1)})\omega,
\end{equation}
and 
\begin{equation}\label{D12 upToHom}
(\Delta\circ\delta+\delta\circ\Delta)\omega=(-1)^r(\delta_{(1)}\circ\Delta_{1,2}-\Delta_{1,2}\circ\delta_{(1)})\omega.
\end{equation}
Eq.'s (\ref{D22 1st}), (\ref{D21 upToHom}) 
and (\ref{D12 upToHom}) were the missing relations to 
imply the following result:

\thm\label{leq4} 
The composition
\begin{equation}\label{trunc4}
C^2_{tot}(\G,\phi)\to<350>^\nabla C^3_{tot}(\G,\phi)\to<350>^\nabla C^4_{tot}(\G,\phi)
\end{equation}
is identically zero. \endthm

For $\vec{\gamma}\in\G_3^2$ as in Eq. (\ref{gammaMatrix}),
$$\mld p_{3,2}(\vec{\gamma})= & \Big{(}\parallel\partial_0^2\vec{\gamma}\parallel_G,c_{2,2}(\vec{\gamma}^{2,2}_\sqcap)\Big{)}+\Big{(}\parallel\vec{\gamma}^{3,1}_\sqcap\parallel_G^{s(\parallel\delta_0\vec{\gamma}\parallel)},c_{3,1}(\delta_0\vec{\gamma})\Big{)}+ \\
	+\Big{(}\big{(}c_{3,1}(\vec{\gamma}^{3,1}_\sqcap)\big{)}^{s(\partial_0^2\delta_0\vec{\gamma})},c_{1,1}(\partial_0^2\delta_0\vec{\gamma})\Big{)}+\Big{(}\parallel\partial_0^2\vec{\gamma}\parallel_G,p_{1,1}^{(3)}(c_{2,1}(\vec{\gamma}^{2,1}_\sqcap);\partial_0^2\delta_0\vec{\gamma})\Big{)}+ \\
	+\Big{(}\delta_0\big{(}c_{3,1}(\vec{\gamma}^{3,1}_\sqcap)\big{)}^{s(\partial_0\delta_0\vec{\gamma})},c_{2,1}(\partial_0\delta_0\vec{\gamma})\Big{)}+\Big{(}\parallel\partial_0\vec{\gamma}\parallel_G,p_{2,1}^{(3)}(c_{1,1}(\vec{\gamma}^{1,1}_\sqcap);\partial_0\delta_0\vec{\gamma})\Big{)}. $$
For $\vec{\gamma}\in\G_2^3$ as in Eq. (\ref{gammaMatrix}), 
$$\mld p_{2,3}(\vec{\gamma})= & \Big{(}\parallel\partial_0\vec{\gamma}\parallel_G,c_{1,3}(\vec{\gamma}^{1,3}_\sqcap)\Big{)}-\Big{(}\parallel\vec{\gamma}^{2,1}_\sqcap\parallel_G^{s(\parallel\delta_0\vec{\gamma}\parallel)},c_{2,2}(\delta_0\vec{\gamma})\Big{)}+ \\
	+\Big{(}\big{(}c_{2,2}^1(\vec{\gamma}^{2,2}_\sqcap)\big{)}^{s(\partial_0\delta_0^2\vec{\gamma})},c_{1,1}(\partial_0\delta_0^2\vec{\gamma})\Big{)}+\Big{(}\parallel\partial_0\vec{\gamma}\parallel_G,p_{1,1}^{(3)}(c_{1,2}(\vec{\gamma}^{1,2}_\sqcap);\partial_0\delta_0^2\vec{\gamma})\Big{)}+ \\
	+\Big{(}\big{(}c_{2,1}(\vec{\gamma}^{2,1}_\sqcap)\big{)}^{s(\partial_0\delta_0\vec{\gamma})},c_{1,2}(\partial_0\delta_0\vec{\gamma})\Big{)}+\Big{(}\parallel\partial_0\vec{\gamma}\parallel_G,p_{2,1}^{(3)}(c_{1,1}(\vec{\gamma}^{1,1}_\sqcap);\partial_0\delta_0\vec{\gamma})\Big{)}, $$
where, for clarity, $c_{2,2}^1(\vec{\gamma}^{2,2}_\sqcap)=(g_{12}^{h_{22}}g_{22},g_{11}^{h_{21}},g_{21})$.

Using these coordinates to define $\Delta_{3,2}$ and 
$\Delta_{2,3}$ implies that for all $\omega\in C^{p,q}_3(\G,\phi)$,
\begin{equation}\label{D32}
(\Delta_{3,1}\circ\delta+\Delta_{2,2}\circ\partial+\Delta_{2,1}\circ\Delta+\Delta\circ\Delta_{2,1}-\partial\circ\Delta_{2,2}+\delta\circ\Delta_{3,1})\omega=\Delta_{3,2}\circ\delta_{(1)}\omega
\end{equation}
and
\begin{equation}\label{D23}
(\Delta_{2,2}\circ\delta+\Delta_{1,3}\circ\partial+\Delta_{1,2}\circ\Delta+\Delta\circ\Delta_{1,2}+\partial\circ\Delta_{1,3}-\delta\circ\Delta_{2,2})\omega=\Delta_{2,3}\circ\delta_{(1)}\omega 
\end{equation}

We conclude defining the necessary higher difference 
maps to extend (\ref{trunc4}) to degree $5$. 
$$\mld p_{3,1}^{(4)}\big{(}f;\pmatrix{\gamma_{11}&\gamma_{12}&\gamma_{13}} & \big{)}=\Big{(}g_{11}^{-1},p_{2,1}^{(3)}\big{(}f;\pmatrix{\gamma_{12}&\gamma_{13}}\big{)}\Big{)}-((g_{12}g_{11})^{-1},f^{h_{12}},g_{12},g_{11})+ \\
+((g_{13}g_{12}g_{11})^{-1},f^{h_{13}}g_{13},g_{12},g_{11})-((g_{13}g_{12}g_{11})^{-1},g_{13},g_{12},g_{11}). $$

$$\mld p_{1,3}^{(4)} & \big{(}f;\pmatrix{\gamma_{11}\cr\gamma_{21}\cr\gamma_{31}}\big{)}=\Big{(}g_{31}^{-1},\big{(}p_{1,2}^{(3)}(f;\pmatrix{\gamma_{12}\cr\gamma_{13}})\big{)}^{h_{31}}\Big{)}-((g_{21}^{h_{31}}g_{31})^{-1},f^{h_{11}i(g_{11})h_{21}h_{31}},g_{21}^{h_{31}},g_{31})+ \\
+((g_{11}^{h_{21}h_{31}}g_{21}^{h_{31}}g_{31})^{-1},(f^{h_{11}}g_{11})^{h_{21}h_{31}},g_{21}^{h_{31}},g_{31})-((g_{11}^{h_{21}h_{31}}g_{21}^{h_{31}}g_{31})^{-1},g_{11}^{h_{21}h_{31}},g_{21}^{h_{31}},g_{31}). $$

Using these coordinates to define $\Delta_{3,1}$ and 
$\Delta_{1,3}$ implies 
\begin{equation}\label{D31 upToHom}
(\partial\circ\Delta_{2,1}-\Delta_{2,1}\circ\partial)\omega=(\Delta_{3,1}\circ\delta_{(1)}-\delta'\circ\Delta_{3,1})\omega
\end{equation}
and
\begin{equation}\label{D13 upToHom}
(\Delta_{1,2}\circ\delta-\delta\circ\Delta_{1,2})\omega=(\Delta_{1,3}\circ\delta_{(1)}-\delta'\circ\Delta_{1,3})\omega, 
\end{equation}
for all $\omega\in C^{p,q}_3(\G,\phi)$. 

Finally,
$$\mld p & _{2,2}^{(4)}\big{(}f;\vec{\gamma}\big{)}=\Big{(}\parallel\vec{\gamma}^{1,2}_\sqcap\parallel_G^{-1},f^{h_{11}h_{21}},c_{1,2}(\vec{\gamma}^{1,2}_\sqcap)\Big{)}-\Big{(}\parallel\delta_0\vec{\gamma}\parallel_G^{-1},f^{t(\gamma_{11})h_{21}},c_{2,1}(\delta_0\vec{\gamma})\Big{)}+ \\
	\qquad-\Big{(}\parallel\vec{\gamma}^{1,2}_\sqcap\parallel_G^{-1},p_{1,2}^{(3)}(f;\partial_0\vec{\gamma})\Big{)}+\Big{(}\parallel\delta_0\vec{\gamma}\parallel_G^{-1},\big{(}p_{2,1}^{(3)}(\vec{\gamma}^{2,1}_\sqcap)\big{)}^{s(\parallel\vec{\gamma}^{1,2}_\sqcap\parallel_G)}\Big{)}+ \\
	\qquad-\Big{(}\parallel\vec{\gamma}^{1,2}_\sqcap\parallel_G^{-1},p_{1,1}^{(3)}(f^{h_{11}},g_{11};\partial_0\delta_0\vec{\gamma})\Big{)}-\Big{(}\parallel\delta_0\vec{\gamma}\parallel_G^{-1},\big{(}p_{1,1}^{(2)}(f;\gamma_{11})\big{)}^{h_{22}},g_{22}\Big{)}+ \\	
	+\Big{(}\parallel\vec{\gamma}\parallel_G^{-1},\delta_1\big{(}f^{h_{12}},p_{2,2}(\vec{\gamma})\big{)}\Big{)}-\Big{(}\parallel\vec{\gamma}\parallel_G^{-1},p_{2,2}(\vec{\gamma})\Big{)}-\big{(}(g_{22}g_{11}^{h_{21}}g_{21})^{-1},f^{h_{11}h_{22}},g_{11}^{h_{22}},g_{22}\big{)}+ \\
	\qquad+\big{(}g_{21}^{-1},g_{22}^{-1},(g_{11}^{h_{22}})^{-1},f^{h_{11}h_{22}}g_{11}^{h_{22}}g_{22}\big{)}-\big{(}g_{21}^{-1},g_{22}^{-1},(g_{11}^{h_{22}})^{-1},g_{11}^{h_{22}}g_{22}\big{)}+ \\
	\qquad-\big{(}g_{21}^{-1},(g_{11}^{h_{21}})^{-1},g_{21}^{-1},f^{h_{11}h_{22}}g_{11}^{h_{22}}g_{22}\big{)}+\big{(}g_{21}^{-1},(g_{11}^{h_{21}})^{-1},g_{21}^{-1},g_{11}^{h_{22}}g_{22}\big{)} $$
for $f\in G$ and 
$\vec{\gamma}=\pmatrix{\gamma_{11}&\gamma_{12}\cr\gamma_{21}&\gamma_{22}}\in\G_2^2$.	

Using the latter to define $\Delta_{2,2}$ implies 
that for all $\omega\in C^{p,q}_3(\G,\phi)$, 
\begin{equation}\label{D22 2nd}
(\Delta_{2,1}\circ\delta-\Delta_{1,2}\circ\partial-\Delta\circ\Delta-\partial\circ\Delta_{1,2}+\delta\circ\Delta_{2,1})\omega=(\Delta_{2,2}\circ\delta_{(1)}-\delta'\circ\Delta_{2,2})\omega;
\end{equation} 
thus, together with Eq.'s (\ref{D32}), (\ref{D23}), 
(\ref{D31 upToHom}) and (\ref{D13 upToHom}), the 
following holds:

\thm\label{leq5} 
The composition
\begin{equation}\label{trunc5}
C^3_{tot}(\G,\phi)\to<350>^\nabla C^4_{tot}(\G,\phi)\to<350>^\nabla C^5_{tot}(\G,\phi)
\end{equation}
is identically zero. \endthm

\refs

\bibitem [Whitehead, 1949]{Whitehead} J. H. C. Whitehead, Combinatorial 
homotopy II. Bull. Amer. Math. Soc. {\em 55} (1949), 453-496.

\bibitem [Loday, 1982]{Loday} J.-L. Loday, Spaces with finitely many  
non-trivial homotopy groups. J. Pure Appl. Algebra {\em 24} (1982), 
179-202.

\bibitem [Norrie, 1987]{Norrie} K. Norrie, Crossed modules and analogues 
of group theorems. Dissertation (1987), King's College, University of 
London.

\bibitem[Baez-Crans, 2004]{Lie2Alg} J. C. Baez and A. S. Crans, 
Higher-dimensional algebra VI: Lie 2-algebras. Theory Appl. Categ. 
{\em 12} (2004), 492-538.

\bibitem[Baez-Lauda, 2004]{Lie2Gps} J. C. Baez  and  A. Lauda, 
Higher-dimensional algebra 5: 2-groups. Theory Appl. Categ. {\em 12} 
(2004), 423-491.

\bibitem[Ellis, 1992]{Ellis} G. J. Ellis, Homology of 2-types. J. London 
Math. Soc. {\em 46} (2) (1992), 1-27.

\bibitem[Ginot-Xu, 2009]{2Coh} G. Ginot and P. Xu, Cohomology of Lie 
2-groups. L'Enseignement Math\'ematique {\em 55} (2) (2009), 1-24.

\bibitem[Wockel, 2011]{IntInfDim} C. Wockel, Categorified central 
extensions, \'etale Lie 2-groups and Lie's Third Theorem for locally 
exponential Lie algebras. Adv. Math. {\em 228} (2011), 2218-2257.

\bibitem[Tseng-Zhu, 2006]{StackyInt} H. H. Tseng and C. Zhu, Integrating 
Lie algebroids via stacks. Compositio Mathematica {\em 142} (2006), 
251-270. 

\bibitem[Sheng-Zhu, 2012]{ZhuInt2Alg} Y. Sheng and C. Zhu, Integration 
of Lie 2-algebras and their morphisms. Lett. Math. Phys. {\em 102} (2) 
(2012), 223-244.
 
\bibitem[vanEst, 1955]{VanEst} W. T. van Est, Une application d'une 
m\'ethode de Cartan-Leray. Proc. Kon. Ned. Akad. {\em 58} (1955) 542-544.

\bibitem[Crainic, 2003]{VanEstC} M. Crainic, Differentiable and 
algebroid cohomology, van Est isomorphisms, and characteristic classes. 
Coment. Math. Helv. Vol.~78 {\em 4} (2003) 681-721.

\bibitem[vanEst-Korthagen, 1964]{NonInt} W. T. van Est and J. Th. Korthagen, 
Non-enlargible Lie algebras. Indag. Math. {\em 26} (1964) 15-31. 

\bibitem[Bursztyn-Cabrera-delHoyo, 2016]{BCD} H. Bursztyn, A. Cabrera 
and M. del Hoyo, Vector bundles over Lie Groupoids and algebroids. Adv. 
Math. {\em 290} (2016) 163-207.

\bibitem[Stefanini, 2008]{LucaPhD} L. Stefanini, On morphic actions and 
integrability of LA-Groupoids. Dissertation (2008), Universit\"at 
Z\"urich Available at http://user.math.uzh.ch/cattaneo/stefanini.pdf  

\bibitem[Neeb, 2002]{Neeb} K.-H. Neeb, Central extensions of 
infinite-dimensional Lie groups. Ann. Inst. Fourier, Grenoble. Vol.~52 
{\em 5} (2002) 1365-1442.

\bibitem[Arias-Abad{\&}Crainic, 2012]{AAC} C. Arias Abad and M. Crainic, 
Representations up to homotopy of Lie algebroids. J. Reine Angew. Math. 
{\em 663} (2012) 91-126.

\bibitem[Arias-Abad{\&}Schatz]{CamFlor} C. Arias Abad and F. Schatz, 
Deformations of Lie brackets and representations up to homotopy. 
Indagationes Mathematicae. {\em 22} (1) (2011), 27-54.

\bibitem[Angulo1]{Lie2AlgCoh} C. Angulo, A new cohomology theory for 
strict Lie 2-algebras. {\em To appear}.

\bibitem[Baues, 1996]{Baues1} H.-J. Baues, Non-abelian extensions and 
homotopies. Non-abelian extensions and homotopies. K-theory {\em 10} (2) 
(1996), 107-133.

\bibitem[Vietes, 1999]{Vietes} A. M. Vietes, Extensiones abelianas, 
cruzadas y 2-extensiones cruzadas de m\'odulos cruzados. Dissertation 
(1999), Universidad de Vigo.

\bibitem[Carrasco-Cegarra-Grandjean, 2002]{(Co)Ho} P. Carrasco, A. M. 
Cegarra and A. R.- Grandjean, (Co)Homology of crossed modules. J. Pure 
Appl. Algebra {\em 168} (2002), 147-176.

\bibitem[Grandjean-Ladra, 1994]{GrandJLadra} M. Ladra and A. R.- 
Grandjean, Crossed modules and homology. J. Pure Appl. Algebra {\em 95} 
(1994), 41-55. 

\bibitem[Angulo2]{2VanEst} C. Angulo, A cohomological proof for the 
integration of Lie 2-algebras. {\em In progress}.
 
\bibitem[Sheng-Zhu, 2012]{IntSubLin2} Y. Sheng and C. Zhu, Integration 
of semidirect product Lie 2-algebras. Int. J. Geom. Methods Mod. Phys. 
{\em Vol.~9} {\bf 5} (2012), 1250043.

\bibitem[Gracia-Saz{\&}Mehta, 2016]{VB&Reps} A. Gracia-Saz and R. Mehta, 
VB-groupoids and representations theory of Lie groupoids. J. Symplect. 
Geom. {\em 15} (3) (2017), 741-783.

\endrefs

\end{document}